%%%The Alternating Sign Matrix Conj. by D. Zeilberger
%%%%%Version of July 31 1995 %%%%%%%%%%%%%%%%%%%%%%%%%%%%%%%
%%%Outputs 84 pages

%begin macros
\baselineskip=14pt
\parskip=10pt
\def\epsilon{\varepsilon}
\def\B{{\cal B}}
\def\S{{\cal S}}
\def\halmos{\hbox{\vrule height0.15cm width0.01cm\vbox{\hrule height
 0.01cm width0.2cm \vskip0.15cm \hrule height 0.01cm width0.2cm}\vrule
 height0.15cm width 0.01cm}}
\font\eightrm=cmr8  
\font\eighttt=cmtt8
\magnification=\magstephalf

\parindent=0pt
\overfullrule=0in
%end macros
\centerline
{\bf PROOF
OF THE ALTERNATING SIGN MATRIX CONJECTURE \footnote{$^1$}
{\eightrm  \raggedright
To appear in Electronic J. of Combinatorics (Foata's 60th Birthday issue).
Version of July 31, 1995;
original version written December 1992. The Maple package
ROBBINS, accompanying this paper, can be downloaded from the
www address in footnote 2 below.
}
} 
\bigskip
\centerline{ {\it Doron ZEILBERGER}\footnote{$^2$}
{\eightrm  \raggedright
Department of Mathematics, Temple University,
Philadelphia, PA 19122, USA. 
\break
E-mail:{\eighttt zeilberg@math.temple.edu}.
WWW:{\eighttt http://www.math.temple.edu/\~{$\,$}zeilberg}.
Anon. ftp: {\eighttt ftp.math.temple.edu}, directory {\eighttt /pub/zeilberg}.
Supported in part by the NSF.
}
}
\bigskip
{\bf Checked by}\footnote{$^3$}
{\eightrm  
See the Exodion for affiliations, attribution, and
short bios.}: 
David Bressoud  and
 
{\eightrm  
Gert Almkvist,
Noga Alon,
George Andrews,
Anonymous,
Dror Bar-Natan,
Francois Bergeron,
Nantel Bergeron, 
Gaurav Bhatnagar,
Anders Bj\"orner,
Jonathan Borwein,
Mireille Bousquest-M\'elou,
Francesco Brenti,
E. Rodney Canfield,
William Chen,
Chu Wenchang,
Shaun Cooper,
Kequan Ding,
Charles Dunkl,
Richard Ehrenborg,
Leon Ehrenpreis,
Shalosh B. Ekhad,
Kimmo Eriksson,
Dominique Foata,
Omar Foda,
Aviezri Fraenkel,
Jane Friedman,
Frank Garvan,
George Gasper,
Ron Graham,
Andrew Granville,
Eric Grinberg,
Laurent Habsieger,
Jim Haglund,
Han Guo-Niu,
Roger Howe,
Warren Johnson,
Gil Kalai,
Viggo Kann,
Marvin Knopp,
Don Knuth,
Christian Krattenthaler,
Gilbert Labelle,
Jacques Labelle,
Jane Legrange,
Pierre Leroux,
Ethan Lewis,
Daniel Loeb,
John Majewicz,
Steve Milne,
John Noonan,
Kathy O'Hara,
Soichi Okada,
Craig Orr,
Sheldon Parnes,
Peter Paule,
Bob Proctor,
Arun Ram,
Marge Readdy,
Amitai Regev,
Jeff Remmel,
Christoph Reutenauer,
Bruce Reznick,
Dave Robbins,
Gian-Carlo Rota,
Cecil Rousseau,
Bruce Sagan,
Bruno Salvy,
Isabella Sheftel,
Rodica Simion,
R. Jamie Simpson,
Richard Stanley,
Dennis Stanton,
Volker Strehl,
Walt Stromquist,
Bob Sulanke,
X.Y. Sun,
Sheila Sundaram,
Rapha\"ele Supper,
Nobuki Takayama,
Xavier G. Viennot,
Michelle Wachs,
Michael Werman,
Herb Wilf,
Celia Zeilberger,
Hadas Zeilberger,
Tamar Zeilberger,
Li Zhang,
Paul Zimmermann .}
\bigskip
{\bf \qquad \qquad \qquad Dedicated to 
my Friend, Mentor, and Guru, Dominique Foata.}
\bigskip
\bigskip
\qquad\qquad\qquad 
{\it Two stones build two houses.
Three build six houses. Four build four and twenty houses. Five  build hundred
and twenty houses. Six build Seven hundreds and twenty houses. Seven build
five thousands and forty houses. From now on, [exit and] ponder
what the mouth cannot speak and the ear cannot hear.}
\smallskip
\qquad\qquad\qquad\qquad
\qquad\qquad\qquad\qquad\qquad\qquad \qquad\qquad\qquad\qquad\qquad\qquad 
(Sepher Yetsira IV,12)
 
\rm
{\bf Abstract:} The number of $n  \times n$ matrices whose
entries are either $-1$, $0$, or $1$, whose row- and column- sums
are all $1$, and such that in every row and every column the non-zero
entries alternate in sign, is proved
to be $[1!4! \dots (3n-2)!]/[n!(n+1)! \dots (2n-1)!]$, as conjectured
by Mills, Robbins, and Rumsey. 
\vfill
\eject
%begin macros
\baselineskip=14pt
\parskip=10pt
\def \inv{\mathop{\rm inv} \nolimits}
\def \sgn{\mathop{\rm sgn} \nolimits}
\def\epsilon{\varepsilon}
\def\B{{\cal B}}
\def\S{{\cal S}}
\def\N{{\cal N}}
\def\halmos{\hbox{\vrule height0.15cm width0.01cm\vbox{\hrule height
 0.01cm width0.2cm \vskip0.15cm \hrule height 0.01cm width0.2cm}\vrule
 height0.15cm width 0.01cm}}
\font\eightrm=cmr8  
\font\eighttt=cmtt8
\magnification=\magstephalf

\parindent=0pt
\headline={\rm  \ifodd\pageno  \RightHead  \else  \LeftHead  \fi}
\def\RightHead{\centerline{Proof  of the ASM Conjecture-Introduction}}
\def\LeftHead{ \centerline{Doron Zeilberger}}
 
\overfullrule=0in
 
%end macros
 
\centerline{\bf INTRODUCTION}
 
The number of permutations (``houses'') that can be made using
$n$ objects (``stones''), for $n \leq 7$, is given in 
{\it Sepher Yetsira} (Ch. IV, v. 12),
a Cabalistic text written more than $1700$ years ago. The general
formula, $n!$, was stated and proved about $1000$ years later by 
Rabbi Levi Ben Gerson (``Ralbag''). The Cabala, which is
a combinatorial Theory Of Everything (both physical and spiritual),
was interested in this problem because $n!$ is the number of
inverse images $\tau^{-1} (w)$ of a generic $n-$lettered word $w$ under the
canonical homomorphism $\tau$:
$$
\tau : \{ \aleph , bet, \dots , shin, tav \}^{*} \rightarrow
\{ \aleph , \dots , tav \}^{*} / 
\{ \aleph\, (bet) = (bet)\, \aleph , \dots , (shin) \, (tav) = (tav) \, 
(shin) \}
\qquad .
$$
The homomorphism $\tau$ is of considerable interest, since two words
$w_1$, and $w_2$ are {\it temura-equivalent} (anagrams) if 
$\tau (w_1) =\tau( w_2)$.
 
A coarser, but just as important, equivalence relation on Hebrew words
and sentences is the one induced by the homomorphism
$$
G: \{ \aleph \, , \, bet \, , \, \dots \, , \, shin \, , \, tav \}^{*}
\rightarrow \N \qquad ,
$$
defined by $G( \aleph ) =1$, $G(bet) =2 ,\, \dots \, , \, G(shin)=300, \,
G(tav)=400$, and
extended homomorphically.
Two words $w_1 , w_2$  are said to be {\it Gematrically equivalent}
if $G( w_1 ) = G( w_2 )$.
 
Modern-day numerologists, like Mills, Robbins, and Rumsey, use a
different kind of {\it Gematria}, one that searches for equality between 
{\it sequences} of cardinalities of
combinatorial families. This enabled them to conjecture a series
of remarkable and tantalizing enumeration identities[MRR1-3],[Stanl]. 
Most of these
conjectures concern a natural generalization of the notion of permutation,
called {\it alternating sign matrix}.
 
A permutation $\pi$ may be described in terms of its corresponding
{\it permutation matrix}, that is the $0-1$ matrix obtained 
by making the $i^{th}$ row have all zeros except
for a $1$ at the $\pi(i)^{th}$ column. What emerges is a $0-1$ matrix
whose row- and column- sums are all $1$.
 
Since $3$ is the cardinality of the set 
$\{ \aleph \,,\, mem \,,\, shin  \}$, the ``mother-letters'', dear
to the authors of Sepher Yetsira, they would 
have most likely enthusiastically approved of the generalization
of permutation matrices, {\it alternating sign matrices}, introduced 
by Robbins and Rumsey[RR] in their study of a determinant-evaluation rule
due to yet another wizard, the Rev. Charles Dodgson. Rather than use only 
the two symbols $0,1$ as entries, an alternating sign matrix is allowed the
use of the three symbols $\{ -1, 1 , 0\}$ (corresponding to guilt, innocence,
and the tongue of the law, respectively). The row- and column- sums
have still to be $1$, and in addition, in every row and every column,
the non-zero elements,
$1,-1$ (right and wrong), have to alternate.
 
Mills, Robbins, and Rumsey[MRR1] discovered that, like their predecessors the
permutations, that are enumerated by the beautiful formula $n!$,
these new mysterious objects seem to be enumerated by an almost equally
simple formula:
$$
A_n := \prod_{i=0}^{n-1} { { (3i+1)!} \over {(n+i)!} } \quad ,
$$
the now famous[R][Z3] sequence $1,2,7,42,429,7436 \ldots$,
first encountered by George Andrews[A1].
 
In this paper I present the {\it first} proof of this fact. 
I do hope that this is not the {\it last} proof, and that
a shorter, more direct, elegant, and combinatorial proof will be found one day.
Meanwhile, this paper should relieve at least some of Dave Robbins's
{\it mathschmerz}, expressed so eloquently in [R]:
 
{\eightrm 
These conjectures are of such compelling simplicity that it is
hard to understand how any mathematician can bear the pain of 
living without understanding why they are true.
}
 
Much more serious than not knowing whether a given fact
is true, is the agony of realizing that our cherished 
{\it tools of the trade} are inadequate to tackle a given problem.
The fact that a conjecture resists vigorous attacks by skilled practitioners
is an impetus for us either to sharpen our existing tools,
or else create new ones. The value of a proof of an outstanding conjecture
should be judged, not by its cleverness and elegance, and not
even by its ``explanatory power'', but by the extent in which it enlarges our
toolbox. By this standard, the present proof is adequate. Like most
new tools, the present method of proof
is a judicious assembly of existing tools, which
I will now describe.
 
The first ingredient consists of {\it partial recurrence equations}
(alias {\it partial difference equations})
and operators. The calculus of finite differences was introduced in the
last century by discrete mathematician George Boole,
but in this century was taken up, and almost monopolized by,
continuous number crunchers, who called them
{\it finite difference schemes}.
A notable exception was Dick Duffin (e.g. [D])
through whose writings I learnt about these objects, and 
immediately fell in love with them
(e.g. [Z1][Z2]). It was fun returning to my first 
love.\footnote{$^4$}{\eightrm ``Mathematicians, like Proust
and everyone else, are at their best when writing about their
first love''- Gian-Carlo Rota (in: `Discrete Thoughts', p.  3).}
 
Conspicuously missing from the present paper is my second love,
{\it bijective proofs} (e.g. [ZB]), that were 
taught to me by Dominique Foata, Xavier G. Viennot, Herb Wilf
and many others. However, doing bijections made me a better mathematician
and person, so their implicit impact is considerable.
 
The second ingredient is my third love, {\it constant term}
identities introduced to me by Dick Askey.
Dennis Stanton[Stant] and John Stembridge[Ste] showed me how to crack
them[Z4][Z5]. The {\it Stanton-Stembridge trick} (see below) 
was indeed crucial.
 
The third and last ingredient, which is not mentioned explicitly, but
without which this proof could never have come to be, is 
my current love: computer
algebra and Maple. Practically every lemma,
sublemma, subsublemma $\dots \,$, was first conjectured with the aid of
Maple,
and then tested by it. 
A Maple package, {\tt ROBBINS}, that empirically
(and in a few cases rigorously,) checks every non-trivial fact proved
in this paper, is given as a companion to this paper, and should be
used in conjunction with it. 
Almost every statement is followed by
a reference to the procedure in {\tt ROBBINS} that empirically
corroborates it. {\eightrm [These are enclosed in square
brackets; For example to see all alternating sign matrices of size 4,
type {\eighttt ASM(4):}.] }
 
Once you have downloaded ROBBINS to your favorite directory,
get into Maple by typing {\tt maple<CR>}. Once in Maple, type
{\tt  read ROBBINS:}, and follow the instructions given there.
(By the way {\it ezra} means `help' in Hebrew.) No knowledge of Maple
is required, except for the fact that every command must end with a 
colon or semi-colon.
For example,  typing {\tt  S15(4):} would verify sublemma $1.5$ for
$k=4$.
 
I wish to thank Shalosh B. Ekhad
for its diligent computations, 
and Russ de Flavia, our dedicated local UNIX guru, for his
constant technical support.
 
This paper would have been little more than a curiosity if not
for George Andrews's[A2] recent brilliant proof of another conjecture
of Mills, Robbins, and Rumsey[MRR3] (conj. 2 of [Stanl]), 
that the number of so-called Totally Symmetric, Self-
Complementary Plane Partitions (TSSCPP) is also given by $A_n$.
All that I show is, that the sequence enumerating ASMs is the same as the
one enumerating TSSCPPs, and then I take a free ride on Andrews's[A2]
result that the latter is indeed $\{ 1,2,7,42,429,  \dots \}$.
 
This paper only settles the first, and simplest, conjecture, concerning
the enumeration of alternating sign matrices. There are many 
variations and refinements
listed in [Stanl][R][MRR2,3].  I am sure that the method of this paper
should be capable of proving all of them. It is also possible that
the present method of proof,
combined with the multi-WZ method[WZ], could be used to prove a 
stronger conjecture of [MRR1,2](conj. 3 of [Stanl]), {\it directly},
in which case the present paper would also furnish an alternative proof
of Andrews's[A2] TSSCPP theorem.
 
\bigskip
\centerline{\bf A MORE GENERAL, AND HENCE EASIER, CONJECTURE}
 
The first step, already undertaken in [MRR2-3], is to find more congenial
``data structures'' for both alternating sign matrices, and for TSSCPPs.
Alternating sign matrices of size $n$ 
are in easy bijection ([MRR2][R]) with monotone
triangles. A {\it monotone triangle} is a triangular
array of  {\it positive}
integers $a_{i,j}$, $1 \leq i \leq n$, $1 \leq j \leq n-i+1$,
such that, $a_{i,j} < a_{i, j+1}$, and 
$a_{i,j} \leq a_{i+1,j} \leq a_{i,j+1}$. In addition we require that
the first row is $1,2, \cdots , n$, i.e. $a_{1,j}=j$, for $j=1 , \cdots , n$.
We will rename them {\it n-Gog} triangles, and the creatures
obtained by chopping off all entries with $j>k$, will be
called {\it $n \times k$-Gog} trapezoids. An $n$-Gog triangle
is an {\it $n \times n$-Gog} trapezoid.
For example the following is one of the $429$ 
$\quad 5-$Gog triangles (formerly called monotone triangles
of size $5$.)
$$
\matrix{ 1 & 2 & 3 & 4 & 5 \cr
1 & 3 & 4 & 5 \cr
2 & 4 & 5 \cr
3 & 5 \cr
4} \quad .
$$
{\eightrm [In order to view all of them type `GOG(5,5):' in ROBBINS.]}
Retaining only the first three columns of the above triangle, yields
one of the $387$ $ \quad 5 \times 3$-Gog trapezoids:
$$
\matrix{ 1 & 2 & 3 \cr
1 & 3 & 4 \cr
2 & 4 & 5 \cr
3 & 5 \cr
4} \quad .
$$
{\eightrm [In order to view all of them type `GOG(3,5):' in ROBBINS.]} 
 
On the TSSCPP side, it was shown in [MRR3] that TSSCPPs whose 3D Ferrers
graphs lie in the cube $[0,2n]^3$ are in trivial bijection with triangular
arrays $c_{i,j}$, $1 \leq i \leq n$, $1 \leq j \leq n-i+1$, 
of integers such that:
(i) $ 1 \leq c_{i,j} \leq j$, \quad (ii) $c_{i,j} \geq c_{i+1,j}$, 
and (iii) $c_{i,j} \leq c_{i,j+1}$. We will call such triangles
{\it n-Magog} triangles, and the corresponding chopped variety,
with exactly the same conditions as above, but $c_{i,j}$ is only
defined for $1 \leq i \leq k$ rather than for $1 \leq i \leq n$,
{\it $n \times k$-Magog trapezoids}. 
For example the following is one of the $429$ $\quad 5-$Magog triangles: 
$$
\matrix{ 1 & 2 & 3 & 3 & 5 \cr
1 & 2 & 2 & 3 \cr
1 & 2 & 2 \cr
1 & 2 \cr
1} \quad .
$$
{\eightrm [In order to view all of them type `MAGOG(5,5):' in ROBBINS.]}
Retaining only the first three rows of the above Magog-triangle, yields
one of the $387$ $\quad 5 \times 3$-Magog trapezoids:
$$
\matrix{ 1 & 2 & 3 & 3 & 5 \cr
1 & 2 & 2 & 3 \cr
1 & 2 & 2
} \quad .
$$
{\eightrm [ In order to view all of them type `MAGOG(3,5):' in ROBBINS.]} 
 
Our goal is to prove the following
statement, conjectured in [MRR3], and proved there for $k=2$.

{\bf Lemma 1:} For $n \geq k \geq 1$, the number of
$n \times k$-Gog trapezoids equals the number of $n \times k$-Magog
trapezoids.
 
{\eightrm [ The number of n by k Magog trapezoids, for specific n and k,
is obtained by typing {\eighttt b(k,n);} while
the number of n by k Gog trapezoids is given by
{\eighttt m(k,n);}. To verify lemma 1, type {\eighttt S1(k,n):}.]}
 
This would imply, by setting $n=k$, that,
 
{\bf Corollary 1':} For $n \geq 1$, the number of $n$-Gog triangles
equals the number of $n$-Magog triangles.
 
Since $n$-Gog triangles are equi-numerous with 
$n \times n$ alternating sign matrices, and $n$-Magog triangles
are equi-numerous with TSSCPPs bounded in $[0,2n]^3$, this
would imply,  together with Andrews's[A2] affirmative resolution of the
TSCCPP conjecture, the following result, that was conjectured in [MRR1].

{\bf The Alternating Sign Matrix Theorem}: The number of $n \times n$
alternating sign matrices, for $n \geq 1$, is:

$$
{
 {1!4! \dots (3n-2)!}
\over
 {n!(n+1)! \dots (2n-1)!}
}   \, = \,
\prod_{i=0}^{n-1} { { (3i+1)!} \over {(n+i)!} } \quad .
$$
 
The rest of this paper will consist of a proof of Lemma 1.
\bigskip 
\centerline{\bf NOMENCLATURE}
 
Throughout this paper, we will meet discrete functions 
$F(n\,;\, a_1 , \dots , a_k )$ of $k+1$ discrete variables, $n$,
$a_1 , \dots , a_k$, where $n \geq 0$, and $( a_1 , \dots , a_k )$
is confined to certain regions of discrete $k$-space $\N^k$, that
depend on $n$. We will make extensive use of the shift operators
$A_i , \, i=1, \, \dots ,\, k$, defined by:
 
$$
A_i F( a_1 , \dots , a_{i} , \dots , a_k ) =
F( a_1 , \dots , a_{i} +1 ,  \dots , a_k ) \qquad .
$$
 
For any (positive or negative) integer $r$, we have:
$$
A_i^r F( a_1 , \dots , a_{i} ,  \dots , a_k ) =
F( a_1 , \dots , a_{i} +r , \dots , a_k ) \qquad ,
$$
and more generally,
$$
A_1^{r_1} \dots A_i^{r_i} \dots A_k^{r_k} \,
F( a_1 , \dots , a_{i} , \dots , a_k ) =
F( a_1 + r_1 , \dots , a_{i} +r_i ,
\dots , a_k + r_k ) \qquad .
$$
 
We denote the identity operator by $I$.
A partial linear recurrence operator (with constant coefficients)
$P( A_1 , \dots , A_k )$ is
any Laurent polynomial in the fundamental shift operators
$A_1 , \dots , A_k$. For example 
$$
(I-A_1^{-1})(I-A_2^{-1})F( a_1 , a_2 ) = F( a_1 , a_2 ) - F( a_1 -1 , a_2 )
-F(a_1, a_2 -1 )+ F( a_1 -1 , a_2 -1 ) \quad .
$$
 
We will also meet polynomials and rational functions
$f( x_1 , \dots , x_k )$ that live on continuous $k$-space. The symmetric
group $\S_k$ of permutations acts on $f( x_1 , \dots , x_k )$  
as follows. For any permutation $\pi=[\pi(1), \ldots , \pi(k)]$, let
$\pi f( x_1 , \ldots , x_k )$ be the rational function obtained from
$f(x_1 , \ldots , x_k )$ by
replacing $x_i$ by $x_{\pi(i)}$, for $i=1 \ldots k$. In symbols:
 
$$
\pi f( x_1 , \dots , x_k )  =
f(x_1 , \dots , x_k)  
\vert_{ \{ x_1 \rightarrow x_{\pi (1)} , \dots , x_k \rightarrow 
x_{ \pi (k) } \} } \qquad ,
$$
 
for example, $[2,3,1] (x_1+2 x_2+3 x_3) = x_2+2 x_3+3 x_1 $.
 
Recall that the {\it number of inversions}, $\inv \pi$, of a permutation
$\pi \in S_k$ is the number of pairs $(i,j)$, with $1 \leq i < j \leq k$
such that $\pi(i) > \pi(j)$.  
Recall also that the {\it sign} of a permutation $\sgn (\pi)$ may be defined
by $(-1)^{\inv \pi }$.
 
A rational function 
$f( x_1 , \dots , x_k )$ is called {\it $\S_k$-antisymmetric},
or {\it antisymmetric} for short, if for $i=1 \, ,\, \dots \,,\, k-1$:
$$
f( x_1 , \dots , x_k ) \vert_{x_{i} \leftrightarrow x_{i+1}}  =
-f( x_1 , \dots , x_k ) \quad ,
$$
or equivalently, since the transpositions
$\{(i,i+1), 1 \leq i \leq k-1 \}$ generate
the symmetric group $\S_k$, the equality :$ \pi f( x) = sgn ( \pi ) f(x)$,
holds for every $\pi \in \S_k$.
 
The group of
{\it signed permutations}, that we will denote by $W(\B_k)$
(since it happens to be the Weyl group of the root system 
$\B_k$, but this is not (directly) relevant to our proof),
consists of pairs $( \pi , \epsilon )$, where $\pi \in \S_k$ 
and $\epsilon$ is a
{\it sign-assignment}:
$\epsilon = ( \epsilon_1 , \dots , \epsilon_k )$, where the
$\epsilon_i$ are either $+1$ or $-1$. A sign assignment
$\epsilon$ acts on $f( x_1 , \dots , x_k )$ by
$$
\epsilon f(x) = f( {\epsilon_1}(x_1) , \dots , {\epsilon_k}(x_k) ) \qquad ,
$$
where, for any variable $t$:
$$
\epsilon_i (t) := t \qquad \hbox{if} \qquad \epsilon_i =1 , \qquad \hbox {and}
\qquad \epsilon_i (t) := 1-t \quad \hbox{if} \quad \epsilon_i =-1 \qquad .
$$
For example $[-1,-1,+1]( x_1^2 x_2^3 x_3^4 )=(1-x_1)^2 (1-x_2)^3 x_3^4$.
 
A signed permutation $( \pi , \epsilon )$ acts on $f( x_1 , \dots , x_k )$
in the following way:
$$
(\pi , \epsilon ) f( x ) =  \epsilon \pi f(x)  \qquad.
$$
For example 
$([2,3,1] , [+1,-1,-1]) [ x_1^2+ 2 x_2+x_3^3] = (1-x_2)^2 +2(1-x_3)+x_1^3$.
 
Now it is time to define the {\it sign} of a signed permutation 
(no pun intended!), $g=(\pi , \epsilon )$.
 
{\bf Definition SIGN:} For any signed permutation $g=(\pi, \epsilon)$,
we define $\sgn (g) := \sgn (\pi)  \sgn( \epsilon)$, where
$\sgn (\epsilon)$ is $(-1)^{\hbox{ [number of $-1$'s]}}$.
 
Throughout this paper $\bar t := 1 -t$. Since this notation will be used
so frequently we will say it again, in bold face:
 
{\bf Crucial Notation :} $\bf { \bar t := 1 - t }$.
 
{\bf Warning:} $\bar t $ Is Not Complex Conjugation!
 
A rational function $f(x_1 , \dots , x_k)$ is $W( \B_k)$-antisymmetric if it is
$\S_k-$antisymmetric, and also 
 
$$
f( \bar x_1 , \dots , x_k ) = - f( x_1 , \dots , x_k ).
$$
 
Equivalently, $f( x_1 , \ldots , x_k)$ is $W(\B_k)$-antisymmetric if for
any signed-permutation $g \in W(\B_k)$, we have:
$g f(x) = \sgn(g) f(x)$. This follows from the well-known (and easy) fact
that the group $W(\B_k)$ is generated by the generators $(i,i+1)$ of
the symmetric group, along with the `sign change' $x_1 \rightarrow \bar x_1$.
 
A {\it Laurent} formal power series $f(x_1, \dots , x_k )$
is anything of the form
$$
f( x_1 , \dots , x_k ) = \sum_{\alpha_1 \geq L , \dots , \alpha_k \geq L }
a_{\alpha_1 , \dots , \alpha_k}
x_1^{\alpha_1} \dots x_k^{\alpha_k} \qquad ,
$$
where  $L$ is a (non-positive) integer. If $L=0$ then it is a run-of-the-mill
formal power series. If $a_{\alpha_1 , \dots , \alpha_k }$ is zero except for
finitely many $\alpha$'s then we have a {\it Laurent polynomial}.
The {\it constant term} of a Laurent formal power series 
$f(x)=f( x_1 , \dots , x_k )$, denoted
by $CT f(x)$, is the coefficient of $x_1^0 \dots x_k^0$, i.e.
$a_{0,0,\dots, 0,0}$.
 
\centerline{\bf CRUCIAL FACTS}
 
{\bf Crucial Fact ${\bf \aleph_0 }$}: A rational function
$f(x_1 , \dots , x_k)$ possesses a Laurent expansion if it is of the form
 
$$ 
f( x_1 , \dots , x_k )=
{{P( x_1 , \dots , x_k )} \over { x_1^{\gamma_1} \dots
x_k^{\gamma_k} Q(x_1 , \dots , x_k)}} \quad ,
$$
 
where $P$ and $Q$ are {\it polynomials} in $(x_1 , \dots , x_k )$,
and {\it most importantly}, $Q$ has a {\it non-zero} constant term,
i.e. $Q(0, 0 , \dots , 0 ) \neq 0$,
and the $\gamma_i$ are (not necessarily positive) {\it integers}.
 
{\bf Proof:} Let the constant term of $Q$, $Q(0,0, \dots , 0)$, be denoted
by $C$. Expanding, we get:
 
$$
{{1}\over {Q(x_1 , \dots , x_k )}  }= 
(C + \sum_{\alpha _1 + \dots + \alpha_k>0} 
q_{\alpha_1 , \dots , \alpha_k} x_1^{\alpha_1} \dots x_k^{\alpha_k})^{-1}=
\sum_{l=0}^{\infty}
(-1)^l {{ (\sum_{\alpha _1 + \dots + \alpha_k>0} 
q_{\alpha_1 , \dots , \alpha_k} x_1^{\alpha_1} \dots x_k^{\alpha_k })^l }
\over {C^{l+1}}} \quad .
$$
 
The right side ``converges'' in the ring of formal power series
(the coefficient of any fixed monomial gets contributions from only
finitely many terms in the above sum). It follows that, under the
condition on $Q$, $P/Q$ also possesses a formal power series, and
dividing by the monomial $x_1^{\gamma_1} \dots x_k^{\gamma_k}$ 
results in a {\it formal Laurent series}. \halmos
 
The condition $Q(0, \dots , 0) \neq 0$ is necessary, since,
for example, $1/(x+y)$ does {\it not} have a formal Laurent series.
(Try it!)
 
{\bf  Iterated  Constant-Terming}
 
The {\it constant term} $CT_x f(x)$  of a rational function 
$f(x)$ of a single variable $x$, is defined to be the
coefficient of $x^0$ in its Laurent expansion.
This is well-defined,
since such a rational function $f(x)$ of the {\it single} variable, $x$,
can always be written in the form $P(x)/(x^a Q(x))$, where
$P(x)$ is a polynomial and $Q(x)$ is a polynomial with
$Q(0) \neq 0$, and $a$ is a non-negative integer. By
crucial fact $\aleph_0$, it possesses a genuine Laurent
expansion, and $CT_x f(x)$ is always well-defined.
 
Consider now an arbitrary rational function $f(x_1 , \dots , x_k)$.
It can be viewed as a rational function in the variable $x_k$,
with coefficients that are rational functions of the other variables
$(x_1, \dots , x_{k-1}  )$. It follows that
$CT_{x_k} f(x_1 , \dots , x_k) $ is well-defined, and is
a certain rational function of $(x_1 , \dots , x_{k-1} )$,
Hence $CT_{x_{k-1}} CT_{x_k} f(x_1 , \dots , x_k)$ is
well-defined, and is a certain rational function of
$(x_1 , \dots , x_{k-2})$ and so on, until we get
that $CT_{x_1} CT_{x_2} \dots CT_{x_k} f(x_1 , \dots, x_k)$,
that we will abbreviate to $CT_{x_1 , \dots , x_k} f(x_1 , \dots , x_k)$,
is well-defined and is a certain number.
 
More formally we have the following recursive definition:
 
{\bf Definition ITERCT:} Let $f( x_1 , \dots , x_k)$ be any rational function.
For $r=1 , \ldots , k$, 
$CT_{x_1 , \dots , x_r} f(x_1 , \ldots , x_k )$ 
is equal to $CT_{x_1} f (x_1 , \ldots , x_k)$ if $r=1$,
and otherwise 
 
$$
CT_{x_1 , \dots , x_r} f( x_1 , \ldots , x_k ) :=
CT_{x_1 , \dots , x_{r-1}} \left [ CT_{x_r} f( x_1 , \ldots , x_k ) \right 
]  \quad .
$$
 
{\bf Crucial Fact ${\bf \aleph_1 }$}: 
Let $f(x_1, \dots , x_k)$ be a 
rational function that possesses a Laurent formal power series, and
that is antisymmetric w.r.t. 
two variables $x_i,x_j$ (i.e. $f \vert_{x_i \leftrightarrow x_j}= -f$.) 
Then $CT_{x_1 , \dots , x_k} f (x_1 , \dots , x_k)=0$. \halmos

{\bf Crucial Fact ${\bf \aleph_2}$}: Let $P(x)$ be a polynomial 
of degree $\leq 2A$, in 
(the single variable) $x$,  that is anti-symmetric w.r.t. to the operation
$x \rightarrow 1-x$ (i.e. $P(1-x) = -P(x)$). Then
 
$$
CT_{x} \left [  {{P(x)} \over { (1-x)^{A+1} x^A } } \right ] =0 .
$$
 
{\bf Proof:} Immediate from crucial fact $\aleph_5$ below. \halmos 
 
{\bf Crucial Fact ${\bf \aleph_3}$}: Let $f(x_1, \dots , x_k)$ be a 
rational function, and define a discrete function
$F( a_1 , \dots , a_k )$ by
 
$$
F( a_1 , \dots , a_k ) :=
CT_{x_1 , \dots , x_k}
 { { f( x_1 , \dots , x_k )} \over {x_1^{a_1} \dots x_k^{a_k} } }
\qquad .
$$
 
Then for any Laurent polynomial $P(x_1 , \dots , x_k )$, we have
 
$$
P( A_1^{-1} , \dots , A_k^{-1} )
F(a_1 , \dots , a_k ) =
CT_{x_1 , \dots , x_k} {P(x_1, \dots , x_k )
 { f( x_1 , \dots , x_k )} \over {x_1^{a_1} \dots x_k^{a_k} } } \qquad .
$$

{\bf Proof:} This is obvious if $P$ is a monomial. Since a
{\it poly}nomial is a linear combination of monomials, and both
sides are linear in $P$, it is also clear in general. \halmos

{\bf Crucial Fact ${\bf \aleph_4}$ (The Stanton-Stembridge trick)}: 
For any permutation $\pi \in \S_k$ and any rational function that
possesses a formal Laurent series,
$f(x)=f(x_1, \dots , x_k)$, we have
 
$$
CT_{x_1 , \dots , x_k} \left [ \pi f(x) \right ] = 
CT_{x_1 , \dots , x_k} f( x ) \quad .
$$
 
Equivalently:
 
$$
CT_{x_{\pi(1)} , \dots , x_{\pi(k)}} f( x ) = 
CT_{x_1 , \dots , x_k} f( x ) \quad .
$$

{\bf Proof:} Applying $\pi$ amounts to renaming the variables. However,
the constant term is obviously unaffected by this renaming. \halmos
 
{\bf Warning:} Crucial fact $\aleph_4$ is false if the rational function
$f$ does not possess a Laurent series, for example when 
$f(x_1 , x_2)=x_1/(x_1+x_2)$.
 
{\bf Crucial Fact ${\bf \aleph_5}$}: Let $A$ be
a non-negative integer, and let $P(x)$ be a Laurent
polynomial in $x$ of degree   $\leq 2A$,  then
 
$$
CT_{x} \left [ {{P(x)} \over { (1-x)^{A+1} x^A }}  \right ]
 = 
CT_{x} \left [ {{P(1-x)} \over { (1-x)^{A+1} x^A }}  \right ]
\quad.
$$
 
{\bf Proof}: $P(x)$ is a linear combination of powers $x^b$, for
$b \leq 2A$. If $P(x)=x^b$, with $A < b \leq 2A$,
both sides vanish, while for $P(x)=x^b$, with $b \leq A$, we have
 
$$
CT_{x} \left [ {{x^b} \over { (1-x)^{A+1} x^A }}  \right ]
=CT_{x} \left [ {{1} \over { (1-x)^{A+1} x^{A-b} }}  \right ]
={{2A-b} \choose {A-b}} 
$$
$$
= {{2A-b} \choose {A}} 
=CT_{x} \left [ {{1} \over { (1-x)^{A-b+1} x^{A} }}  \right ]
=CT_{x} \left [ {{(1-x)^b} \over { (1-x)^{A+1} x^{A} }}  \right ]
\quad . \quad \halmos
$$
 
We also need the celebrated
 
{\bf Vandermonde's Determinant Identity:} 
$$
\sum_{\pi \in \S_k} sgn( \pi ) \,\cdot \, \pi ( \prod_{i=1}^{k} x_i^{i-1} ) =
\prod_{1 \leq i < j \leq k} ( x_j - x_i ) \qquad .
$$
{\bf Proof:} View both sides as polynomials in $x_1$ of degree $k-1$.
Both sides agree (and in fact are zero) at $x_1=x_2 , \dots, x_1=x_k$,
and they also  agree at $x_1 =0$ by induction on $k$. 
Since they agree at $k$ distinct values, they agree
everywhere. \halmos 
 
A longer, but much nicer, proof 
(that turned out to be seminal[ZB]) was given by Ira Gessel[G].
 
{\bf Iterated Residues}
 
Later, we will find it necessary to convert 
our constant-term
expressions to ``residue'' expressions. People who, like myself, (and
John Riordan), are horrified by analysis, need not worry. 
The {\it residue}, $Res_x f(x)$, of a rational function of
a single variable  $x$, is defined as the coefficient of $x^{-1}$ 
in the (formal) Laurent expansion of $f(x)$. It is well-defined
for the same reason that $CT_x f(x)$ is: a rational function of a single
variable always possesses a Laurent series, thanks to crucial fact
$\aleph_0$.
The iterated residue $Res_{x_1 , \dots , x_k}f$ is defined
to be $Res_{x_1}[Res_{x_2} [ \dots [ Res_{x_k} f] \dots ] ]$.
{\eightrm [Try out IterRes in {\eighttt ROBBINS}]}. 
 
{\bf Trivial (Yet Crucial) Fact $\aleph_6$ :}
$$
CT_{x_1 , \dots , x_k} f( x_1 , \dots , x_k )=
Res_{x_1 , \dots , x_k} \left [ 
{ {f( x_1 , \dots , x_k )} \over  {x_1 \dots x_k} }
\right ] \quad . \quad \halmos
$$

For future reference, we need to recast crucial facts $\aleph_4$
and $\aleph_5$ in terms of residues.
 
{\bf Crucial Fact ${\bf \aleph_4'}$ (The Stanton-Stembridge trick)}: 
For any permutation $\pi \in \S_k$ and any rational function that
possesses a formal Laurent series,
$f(x)=f(x_1, \dots , x_k)$, we have
$$
Res_{x_1 , \dots , x_k} [\pi f(x) ]= 
Res_{x_1 , \dots , x_k} f( x ) \quad . 
$$
Equivalently:
$$
Res_{x_{\pi(1)} , \dots , x_{\pi(k)}} f( x ) = 
Res_{x_1 , \dots , x_k} f( x ) \quad . \quad \halmos
$$
 
{\bf Crucial Fact ${\bf \aleph_5'}$}: Let $A$ be
a non-negative integer, and let $P(x)$ be a Laurent
polynomial in $x$ of degree   $\leq 2A$,  then
 
$$
Res_{x} \left [ {{P(x)} \over { (1-x)^{A+1} x^{A+1} }}  \right ]
 = 
Res_{x} \left [ {{P(1-x)} \over { (1-x)^{A+1} x^{A+1} }}  \right ]
\quad. \quad \halmos 
$$
 
{\bf Antisymmetrizers}
 
{\bf Crucial Fact ${\bf \aleph_7}$}: For any rational function
$f(x_1 , \ldots , x_k)$, its {\it antisymmetrizer} w.r.t. $\S_k$:
$$
\sum_{\pi \in \S_k} \sgn(\pi)\cdot  \pi f(x_1 , \ldots , x_k ) \quad ,
$$
is an $\S_k$- antisymmetric function.  \halmos 
 
{\eightrm [ Try  {\eighttt antisymmetrizerS\_k} in ROBBINS]}.
 
{\bf Crucial Fact ${\bf \aleph_7'}$}: For any rational function
$f(x_1 , \ldots , x_k)$, its {\it antisymmetrizer} w.r.t. $W(\B_k)$:
 
$$
\sum_{g \in W(\B_k)} \sgn(g) \cdot g f(x_1 , \ldots , x_k ) \quad ,
$$
 
is a $W(\B_k)$-anti-symmetric function.  \halmos
 
{\eightrm [ Try  {\eighttt antisymmetrizerWB\_k} in ROBBINS]}.
 
{\bf Crucial Fact ${\bf \aleph_8}$}: Any 
$S_k-$ antisymmetric polynomial is divisible by
 
$$
\prod_{1 \leq i < j \leq k} ( x_j - x_i) \quad .
$$
 
{\bf Proof:} Let's view the polynomial as a polynomial in $x_1$.
By anti-symmetry, it vanishes at $x_1=x_i$, for $i=2,3, \ldots , k$,
hence is divisible by $(x_1-x_2) \cdots (x_1-x_n)$. Similarly
(or by (anti-)symmetry) it is also divisible by all $(x_i-x_j)$,
$2 \leq i < j \leq k$, and hence by their product.
 
{\bf Crucial Fact ${\bf \aleph_8'}$}: Any 
$W(\B_k)-$ antisymmetric polynomial is divisible by
 
$$
\Delta_k ( x_1 , \dots , x_k) :=
\prod_{i=1}^{k} (1-2 x_i )
\prod_{1 \leq i < j \leq k} ( x_j - x_i) (x_j+x_i-1) \quad .
\eqno(Delta)
$$
 
{\bf Proof:} Let's view the polynomial as a polynomial in $x_1$.
By $W(\B_k)$-antisymmetry, it vanishes at $x_1=x_i$, for $i=2,3, \ldots , k$,
as well as at $x_1=\bar x_i$, and at $x_1=1/2$.
Hence it is divisible by 
 
$$
(1-2 x_1) \prod_{j=2}^{k} (x_1-x_j) (x_1+x_j-1) \quad.
$$
 
Similarly, (or by (anti-)symmetry) it is divisible by all the other
factors of $\Delta_k$, and hence by $\Delta_k$ itself. \halmos
\vfill
\eject
%%%The Alternating Sign Matrix Conj.- Act I; Version of April 1995
%%%File asmI.tex
%begin macros
\baselineskip=14pt
\parskip=10pt
\def \inv{\mathop{\rm inv} \nolimits}
\def \sgn{\mathop{\rm sgn} \nolimits}

\def\epsilon{\varepsilon}
\def\B{{\cal B}}
\def\S{{\cal S}}
\def\halmos{\hbox{\vrule height0.15cm width0.01cm\vbox{\hrule height
 0.01cm width0.2cm \vskip0.15cm \hrule height 0.01cm width0.2cm}\vrule
 height0.15cm width 0.01cm}}
\font\eightrm=cmr8  
\font\eighttt=cmtt8
\magnification=\magstephalf

\parindent=0pt
\overfullrule=0in
%\pageno=12
\headline={\rm  \ifodd\pageno  \RightHead  \else  \LeftHead  \fi}
\def\RightHead{\centerline{Proof  of the ASM Conjecture-Act I}}
\def\LeftHead{ \centerline{Doron Zeilberger}}
 
%end macros
 
\centerline{\bf Act I. COUNTING MAGOG}
 
\qquad\qquad\qquad {\it  And I will send a fire on Magog $\dots$
 and they shall know that I am the Lord.}
\smallskip
\qquad\qquad\qquad\qquad
\qquad\qquad\qquad\qquad\qquad\qquad \qquad\qquad\qquad\qquad\qquad\qquad 
(Ezekiel XXXIX,6)
 
Recall that an $n \times k$-Magog trapezoid,
where $n \geq k \geq 1$, is a trapezoidal array of integers
$( c_{i,j} )$, $1 \leq i \leq k$, $1 \leq j \leq n-i+1$, such that:
$$
(i)\quad 1 \leq c_{i,j} \leq j \quad , \quad (ii)\quad c_{i,j} 
\geq c_{i+1,j} \quad,  \quad (iii) \quad c_{i,j} \leq c_{i,j+1} \quad .
$$
{\bf Sublemma 1.1:} The total number of $n \times k -$ Magog trapezoids, let's
call it $b_k (n)$, is given by the following constant term expression:
$$
b_k (n) = CT_{x_1, \dots , x_k}  \,
\left \{ 
{
{\Delta_k (x_1 , \dots , x_k )} \over
{\prod_{i=1}^k x_i^{n+k-i-1} ( \bar x_i)^{n+k+1} \,  
 \prod_{1 \leq i < j \leq k} (1- x_i x_j ) }
} \right \} ,
\eqno(MagogTotal)
$$
where,
$$
\Delta_k (x_1 , ... , x_k ) :=
\prod_{i=1}^{k}( 1- 2 x_i )
\prod_{1 \leq i < j \leq k}
(x_j - x_i )(x_j+x_i-1) \qquad .
\eqno(Delta)
$$
 
{\eightrm  [ Type `{\eighttt S11(k,n):}' in ROBBINS, for specific k,n.]}
 
{\bf Proof:}  
Let $B_k ( n ; a_1 , ... , a_k )$ be the number of $n \times k$
Magog trapezoids $( c_{i,j} )$, such that for $i = 1 , \dots , k$,
$c_{i,n-i+1} = a_i$. In other words, the rightmost border is prescribed
by the $a_i$. By the definition of Magog trapezoids, the natural domain of
the $(k+1)-variable$ discrete function $B_k( - ; -, \dots , -)$ is
$$
\displaylines{
Land\_Of\_Magog_k :=
\cr
\{  
(n \, ; \, a_1 , \dots , a_k) \, | \, n \geq k \, , \quad
n \geq a_1 \geq a_2 \geq \dots \geq a_k \geq 1 \, , \quad \hbox{and} \quad
a_i \leq n-i+1  \quad \hbox{for} \quad i= 1\,,\, \dots \, , \, k \}.
\cr}
$$
 
It is convenient  to extend this region to the larger region
$$
\{  
(n \, ; \, a_1 , \dots , a_k) \, | \, n \geq k \, , \quad
n \geq a_1 \geq a_2 \geq \dots \geq a_k \geq 1 \, \},
$$
and to define $B_k(n; a_1 , \dots , a_k)$ to be $0$ whenever
$a_i > n-i+1$ for one or more $i$ ($1 \leq i \leq k$). This
makes perfect combinatorial sense, since the number of
Magog-trapezoids with such $a_i$'s, that break the rules, is $0$.
 
The following subsublemma gives a constant term expression for
$B_k$.
 
{\bf Subsublemma 1.1.1:}  Let
$$
C_k(n; a_1 , \dots , a_k ) :=
CT_{x_1 , \dots , x_k} \left \{ 
{ { \Delta_k ( x_1 , ... , x_k ) } 
\over { \prod_{i=1}^k x_i^{a_i+k-i-1} (\bar x_i)^{k+n} } }  \right \} \quad ,
$$
where $\Delta_k$ is defined in Eq. $(Delta)$ in
the statement of sublemma $1.1$ above. For $n \geq k \geq 1$ and 
for all $(n ; a_1 , \dots , a_k )$  for which 
$n \geq a_1 \geq \dots \geq a_k \geq 1$, we have:
$$
B_k (n\,;\,a_1 , ... , a_k ) =
C_k (n\,;\,a_1 , ... , a_k ) \quad .
$$
{\eightrm  [ Type `{\eighttt S111(k,n):}' in ROBBINS, for specific k,n.]}
 
{\bf Proof:} It is convenient to extend $B_k$ 
even further to the following larger domain:
$$
\overline {Land\_Of\_Magog_k} :=
$$
$$
\{  
(n \, ; \, a_1 , \dots , a_k) \, | \, n \geq k \geq 1 , \quad
n - a_1 \geq -1 \, , \,  a_1 - a_2 \geq -1 , \, \dots , \,
 a_{k-1} -  a_k  \geq -1 \, , \, a_k \geq 0  \quad \} \quad ,
$$
and to define $B_k$ to be zero in the ``no-man's-land''
points that are in 
$\overline {Land\_Of\_Magog_k}$ but not in 
$Land\_Of\_Magog_k$. 
 
We will prove that subsublemma $1.1.1$  
holds in this larger domain $\overline {Land\_Of\_Magog_k}$.
This would follow from the following three subsubsublemmas.
 
{\bf Subsubsublemma 1.1.1.1:} Let $k \geq 1$. The following partial
recurrence holds in the subset of $\overline{Land\_Of\_Magog_k}$ 
of points $(n; a_1 , \dots , a_k)$
for which $n>k$, and $n \geq a_1 \geq a_2 \geq \dots \geq a_k \geq 1$:
$$
\left \{ \prod_{i=1}^{k}
(I- A_i^{-1}) \right \} \, B_k (n \, ; \, a_1 , \dots , a_k )\,
=\, B_k ( n-1 \, ; \, a_1 , \dots , a_k ) \quad .
\eqno(P \Delta E (B_k)) 
$$
 
$B_k (n; a_1 , \dots , a_k )$ also satisfies the
following initial/boundary conditions in \break
$ (n \,; \, a_1 , \dots , a_k ) \in \overline {Land\_Of\_Magog_k}$:
 
For $i=1,2, \dots, k-1$,
$$
B_k (n \, ; \, a_1 , \dots , a_k ) = 0 \quad, \qquad \hbox{on} \qquad 
a_i - a_{i+1}=-1 \quad .
\eqno(B1_i)
$$
Also:
$$
B_k (n \, ; \, a_1 , \dots , a_k ) = 0 \quad, \qquad \hbox{on} \qquad 
a_k =0 \quad .
\eqno(B1_k)
$$
Furthermore:
$$
B_k (n \, ; \, n+1 \,,\, a_2 \,,\, \dots \,,\, a_k ) = 0  \quad ,
\hbox{and}
\eqno(B2)
$$
$$
B_k ( k \, ; \, a_1 , \dots , a_{k-1} , a_k ) = 
B_{k-1} ( k \, ; \, a_1 , \dots , a_{k-1} ) \delta ( a_k , 1 )  \quad ,
\eqno(B3)
$$
where, as usual, $\delta (i,j)$ is $1$ or $0$ according
to whether $i=j$ is true or false respectively.
Finally,  $B_k(n; -)$ satisfies the ``Eve'' condition:
$$
B_1(1; a_1 ) = \delta ( a_1 , 1 ) , \qquad 
 (1 \,; \, a_1 ) \in \overline {Land\_Of\_Magog_1} \quad .
\eqno(B4)
$$
{\eightrm  [ Type `{\eighttt S1111all(k,n):}' in ROBBINS, for specific k,n.]}
 
{\bf Proof:} Let's try to compute $B_k (n \,;\, -)$ from 
$B_k ( n-1 \,;\, - )$.
By looking at what vectors $( b_1 , \dots , b_k)$ 
can occupy the diagonal immediately to the
left of the rightmost one, $c_{i,n-i}=b_i \,$, say, we get
a recurrence, valid for $n > k$ and
all $(n; a_1 , \dots , a_k ) \in Land\_Of\_Magog_k$:
$$
B_k (n\,;\,a_1 , \dots , a_k ) =
\sum_{( b_1 , \dots , b_k )} B_k ( n-1 \,;\, b_1 , \dots , b_k ) \quad .
\eqno(Ekhad)
$$
 
The summation extends over all the $b=(b_1 , \dots , b_k)$ such that
$a_{i+1} \leq b_i \leq \min ( a_i , n-i )$ , $\, i= 1 \dots k-1 \,$, 
$\, b_k \leq \min ( a_k , n-k )$. Note  that these conditions imply that
$n-1 \geq b_1 \geq b_2 \dots \geq b_k$.
 
$(Ekhad)$ can be used to compile a table of $B_k$, together with the
initial condition
$$
B_k(k; a_1 , \dots , a_{k-1}, 1)=
B_{k-1} (k; a_1 , \dots , a_{k-1})  \quad ,
$$
that holds because when $n=k$, we have a Magog triangle, and for it
necessarily $c_{k,1}=1$ and there is a one-one correspondence
between $k \times k$ Magog triangles and $ k \times (k-1)$ Magog-trapezoids
obtained by deleting that $c_{k,1}$ entry (or putting it back, if
one wishes to go the other way.) This enables us to
go down in the $k-$ladder until we reach $k=1$ for which
$B_1(1;1)=1$.
 
We can extend the set of $(n-1; b_1 , \dots , b_k )$ over which
the summation in $(Ekhad)$ takes place, so that $(Ekhad)$ becomes $(Ekhad')$:
$$
B_k ( n\,;\, a_1 , \dots , a_k )  =
\sum_
{a_{i+1} \leq b_i \leq  a_i  \,\, ( 1 \leq i \leq k )}
B_k ( n-1 \,;\, b_1 , b_2 ,  \dots , b_k ) \quad,
\eqno(Ekhad')
$$
(where we put $a_{k+1}=0$,)
since, whenever $b_i > n-i$, the contribution is zero, by our extended
definition of $B_k$. The next thing to observe is that $(Ekhad')$ holds
not only for $(n; a_1 , \dots , a_k) \in Land\_Of\_Magog_k$ 
but for all points for which $n \geq a_1 \geq \dots \geq a_k \geq 1$,
since for these extra points, 
the left side is $0$ by the extended definition of $B_k$,
and all the terms on the right side are also $0$ for the same reason.
Now,
$$
\displaylines{
\{ (n-1 \, ; \, b_1 , \dots , b_k ) \, \vert \, a_2 \leq b_1 \leq a_1 \, , \,\,
a_3 \leq b_2 \leq a_2 \, ,\, \dots \,,\, a_{k} \leq b_{k-1} \leq a_{k-1} \, ,
\, b_k \leq a_k \,\,  \} \, \backslash 
\cr
\{ (n-1 \, ; \, b_1 , \dots , b_k ) \, \vert \, a_2 \leq b_1 \leq (a_1 -1), 
\,\,
a_3 \leq b_2 \leq a_2 \, , \, \dots \, , \, a_{k} \leq b_{k-1} \leq a_{k-1}
 \, , b_k \leq a_k \, \}
\cr
=
\{ (n-1 \, ; \, a_1 , b_2 , \dots , b_k ) \,\vert \,
a_3 \leq b_2 \leq a_2 \, , \, \dots \, , \, a_{k} \leq b_{k-1} \leq a_{k-1} 
\, , \,
b_k \leq a_k  \, \} \quad .
\cr}
$$
 
Using $(Ekhad')$, we get that for all 
$n \geq a_1 \geq \dots \geq a_k \geq 1$:
$$
(I- A_1^{-1} ) B_k ( n\,;\, a_1 , \dots , a_k )  =
\sum_
{a_{i+1} \leq b_i \leq  a_i  \,\,(2 \leq i \leq k) }
B_k ( n-1 \,;\, a_1 , b_2 ,  \dots , b_k ) \quad,
$$
(with the convention that $a_{k+1}=0$,)
for $(n; a_1 , \dots , a_k )$ for which $n>k$ and 
$n \geq a_1 \geq a_2 \geq \dots \geq a_k \geq 1$.
 
Iterating the process, we get, in turn, that
$$
(I- A_1^{-1} ) \dots (I- A_l^{-1}) B_k ( n\,;\, a_1 , \dots , a_k )  =
\sum_
{a_{i+1} \leq b_i \leq  a_i \,\, (l+1 \leq i \leq k) }
B_k ( n-1 \,;\, a_1 , \dots , a_l ,b_{l+1} ,  \dots , b_k ) \quad,
$$
for $l=2, \dots ,k$. The final $l=k$  is  $(P \Delta E(B_k))$.
 
The boundary/initial conditions (B1)-(B2) follow from the fact
that $B_k$ is defined to be $0$ in the `boundary'
$\overline{Land\_Of\_Magog_k}\,\,\backslash \,\,Land\_Of\_Magog_k$, while
$(B3)$ and $(B4)$ follow straight from the (combinatorial) definition
of $B_k$. This completes the proof of  $(sub)^3$lemma $1.1.1.1$. \halmos
 
{\bf Subsubsublemma 1.1.1.2:}  Let $C_k (n: a_1 , \dots , a_k)$ be
the constant term expression defined in the statement of
subsublemma $1.1.1$.
The following partial
recurrence holds in the subset of $\overline{Land\_Of\_Magog_k}$ 
of points $(n; a_1 , \dots , a_k)$
for which $n>k$, and $n \geq a_1 \geq a_2 \geq \dots \geq a_k \geq 1$:
$$
\left \{ \prod_{i=1}^{k}
(I- A_i^{-1}) \right \} \, C_k (n \, ; \, a_1 , \dots , a_k )\,
=\, C_k ( n-1 \, ; \, a_1 , \dots , a_k ) \quad .
\eqno(P \Delta E (C_k)) 
$$
$C_k (n; a_1 , \dots , a_k )$ also satisfies the
following initial/boundary conditions for \break
$ (n; a_1 , \dots , a_k ) \in \overline {Land\_Of\_Magog_k}$.
 
For $i=1,2, \dots, k-1$,
$$
C_k (n \, ; \, a_1 , \dots , a_k ) = 0 \quad, \qquad \hbox{on} \qquad 
a_i - a_{i+1}=-1 \quad .
\eqno(C1_i)
$$
Also:
$$
C_k (n \, ; \, a_1 , \dots , a_k ) = 0 \quad, \qquad \hbox{on} \qquad 
a_k =0 \quad .
\eqno(C1_k)
$$
Furthermore:
$$
C_k (n \, ; \, n+1 \,,\, a_2 \,,\, \dots \,,\, a_k ) = 0  \quad ,\hbox{and}
\eqno(C2)
$$
$$
C_k ( k \, ; \, a_1 , \dots , a_{k-1} , a_k ) = 
C_{k-1} ( k \, ; \, a_1 , \dots , a_{k-1} ) \delta ( a_k , 1 )  \quad ,
\eqno(C3)
$$
and
$$
C_1(1; a_1 ) = \delta ( 1, a_1 ) \quad , \quad  (1; a_1 ) \in 
\overline{Land\_Of\_Magog_1} \qquad .
\eqno(C4)
$$
{\eightrm  [ Type `{\eighttt S1112(k):}' in ROBBINS, for specific k.]}
 
{\bf Proof:} 
$(P \Delta E (C_k))$ follows from crucial fact $\aleph_3$.
(In fact it holds for all $(n; a_1 , \dots , a_k)$, without
the indicated restriction, but we only need it there.)
 
$(C1_i)$, ($1 \leq i <k$), is satisfied by crucial fact
$\aleph_1$. $(C1_k)$ is satisfied, since when $a_k=0$, the
constant-termand of the constant term expression defining
$C_k$ is a multiple of $x_k$, and hence its constant term
w.r.t. to $x_k$, and hence w.r.t. to all variables, is $0$.
$(C2)$ is satisfied  thanks to
crucial fact $\aleph_2$, since the degree of the numerator in $x_1$,
$2k-1$ is $\leq 2(n+k-1)$ (because we have $n \geq k \geq 1$.)
 
The proof of  $(C3)$ is not quite so easy. 
First let's prove $(C3)$ when $a_k=1$. Since
$$
\Delta_k ( x_1 , \dots , x_k ) =
\Delta_{k-1} ( x_1 , \dots , x_{k-1} ) \cdot
(1- 2 x_k )
\prod_{i=1}^{k-1} ( x_k - x_i ) ( x_k - \bar x_i )
\quad ,
$$
we have
$$
\displaylines{
C_k ( k ; a_1 , \dots , a_{k-1} , 1)=
\cr
CT_{x_1 , \dots , x_{k-1}} \left \{
{
{ \Delta_{k-1} ( x_1 , \dots , x_{k-1} ) }
\over
{ \prod_{i=1}^{k-1} (x_i^{a_i+k-i-1} \bar x_i^{2k}) }
} \, \cdot \,
CT_{x_k} \left [
{
{ (1- 2 x_k )  \prod_{i=1}^{k-1} (x_i - x_k ) ( \bar x_i - x_k )  }
\over
{{x_k}^0 (\bar x_k)^{2k} }
} \right ]
\right \}
\cr
= CT_{x_1 , \dots , x_{k-1}}
\left \{
{
{ \Delta_{k-1} ( x_1 , \dots , x_{k-1} ) }
\over
{ \prod_{i=1}^{k-1} (x_i^{a_i+k-i-1} \bar x_i^{2k}) }
} \, \cdot \,
\left [ \prod_{i=1}^{k-1} (x_i \bar x_i) \right ]
\right \}
\cr
= CT_{x_1 , \dots , x_{k-1}}
\left \{
{
{ \Delta_{k-1} ( x_1 , \dots , x_{k-1} ) }
\over
{ \prod_{i=1}^{k-1} (x_i^{a_i+(k-1)-i-1} \bar x_i^{(k-1)+ k}) }
}
\right \}
\, =
C_{k-1} (k; a_1 , \dots , a_{k-1} ) \quad.
\cr}
$$
 
We will now show that $(C3)$ holds when $a_k \geq 2$.
When $(k;a_1 , \dots , a_{k-1},a_k) \in \overline {Land\_Of\_Magog_k}$
does not satisfy $k \geq a_1 \geq a_2 \geq \dots \geq a_k $,
then the left side of $(C3)$ is already known to be $0$ in virtue
of $(C1)$ or $(C2)$. We are left with the task of showing that 
when $k \geq a_1 \geq \dots \geq a_k \geq 2$,
the following expression vanishes:
$$
C_k(k; a_1 , \dots , a_k ) =
CT_{x_1 , \dots , x_k}  \left \{
{ { \Delta_k ( x_1 , ... , x_k ) } 
\over { \prod_{i=1}^k x_i^{a_i+k-i-1} (\bar x_i)^{2k} } } \right \}
=
CT_{x_1 , \dots , x_k}  \left \{
{ { \left [ \prod_{i=1}^k x_i^{k-a_i+i} \right ] \Delta_k ( x_1 , ... , x_k ) } 
\over { \prod_{i=1}^k x_i^{2k-1} (\bar x_i)^{2k} } } \right \}
$$
$$
=
Res_{x_1 , \dots , x_k} \left \{
{ { \left [ \prod_{i=1}^k x_i^{k-a_i+i} \right ] \Delta_k ( x_1 , ... , x_k ) } 
\over { \prod_{i=1}^k x_i^{2k} (\bar x_i)^{2k} } } \right \} \quad .
\eqno(Efes)
$$
 
In order to show that the right side of $(Efes)$ indeed vanishes,
we will first prove a $(sub)^4$lemma that asserts
that, under the stated conditions on the $a_i$, the right side of $(Efes)$
remains invariant whenever its residuand is replaced by
any of its images under the action of the group of signed permutations
$W(\B_k)$.
 
{\bf Subsubsubsublemma 1.1.1.2.1:} Let 
$F_{k;a_1, \dots , a_k} (x_1 , \dots , x_k )$ be the residuand of
the right side of $(Efes)$, i.e.
$$
F_{k;a_1, \dots , a_k} (x_1 , \dots , x_k ) :=
{ { \left [ \prod_{i=1}^k x_i^{k-a_i+i} \right ] \Delta_k ( x_1 , ... , x_k ) } 
\over { \prod_{i=1}^k x_i^{2k} (\bar x_i)^{2k} } } \quad ,
$$
where $k \geq a_1 \geq a_2 \dots \geq a_k \geq 2$, and $\Delta_k$ is
defined in Eq. $(Delta)$ in the statement of sublemma $1.1$. Then
for any signed permutation $g=( \pi , \epsilon) \in W(\B_k)$, we 
have
$$
Res_{x_1 , \dots , x_k} \left \{
g \, [ \, F_{k;a_1, \dots , a_k}(x) \, ]  \,\right \} =
Res_{x_1 , \dots , x_k} \left \{ F_{k;a_1, \dots , a_k} (x)  \right \} \quad .
$$
{\eightrm  [ Type `{\eighttt S11121(k,a):}' in ROBBINS, for specific k and a.]}
 
{\bf Proof:} $\Delta_k$, defined in  $(Delta)$, is a polynomial
of degree $2k-1$ in any one of its variables.
By Crucial Fact $\aleph_5'$ we can repeatedly change
all the $-1$'s in $\epsilon$ to $+1$'s since the numerator is
a polynomial of degree $\leq (2k-1) +(k+k-2)=4k-3 \leq 2(2k-1)$.
Finally we can change the $\pi$ in $g$ into the identity element, by
Crucial Fact $\aleph_4'$. This completes the proof of
$(sub)^4$lemma $1.1.1.2.1$. \halmos
 
We can now replace the right side of $(Efes)$ by the iterated residue of the
average of all its images under $W(\B_k)$:
$$
Res_{x_1 , \dots , x_k} \left \{ F_{k;a_1, \dots , a_k} (x) \right \} =
{{1} \over {2^k k!}} \sum_{g \in W(\B_k)}
Res_{x_1 , \dots , x_k} \left \{ g[ F_{k;a_1, \dots , a_k} (x)] \right \}
$$
$$
=
{{1} \over {2^k k!}} Res_{x_1 , \dots , x_k} \left \{
 \sum_{g \in W(\B_k)}  g \, [ \,  F_{k;a_1, \dots , a_k} (x) \, ] \,
 \right \} \quad .
\eqno(Efes')
$$
 
We now will prove  that not only is the
{\it iterated residue} on the right side of $(Efes')$ identically zero,
but so is the whole {\it residuand}.
 
{\bf Subsubsubsublemma 1.1.1.2.2:} Let 
$F_{k;a_1, \dots , a_k} (x_1 , \dots , x_k )$ be the residuand of
the right side of $(Efes)$, i.e.
$$
F_{k;a_1, \dots , a_k} (x_1 , \dots , x_k ) :=
{ { \left [\prod_{i=1}^k x_i^{k-a_i+i} \right ] \Delta_k ( x_1 , ... , x_k ) } 
\over { \prod_{i=1}^k x_i^{2k} (\bar x_i)^{2k} } } \quad ,
$$
where $k \geq a_1 \geq a_2 \dots \geq a_k \geq 2$, and $\Delta_k$ is
defined in Eq. $(Delta)$ in the statement of sublemma $1.1$. Then
$$
\sum_{g \in W(\B_k)}  g \left [ F_{k;a_1, \dots , a_k} (x) \right ] =0 \quad .
\eqno(Efes'')
$$
{\eightrm  [ Type `{\eighttt S11122(k,a):}' in ROBBINS, for specific k and a.]}
 
{\bf Proof:}  By the symmetry of the denominator and the anti-symmetry
of $\Delta_k$ w.r.t. the group of signed permutations $W(\B_k)$,
we have:
$$
 \sum_{g \in W(\B_k)} g \left [ F_{k;a_1, \dots , a_k} (x) \right ] =
{   {\left \{ 
 \sum_{g \in W(\B_k)} \sgn(g) \cdot 
g \left [ \prod_{i=1}^k x_i^{k-a_i+i} \right ]
    \right \}
\Delta_k ( x_1 , ... , x_k ) 
} 
\over { 
       \prod_{i=1}^k x_i^{2k} (\bar x_i)^{2k}  
      } } \quad .
\eqno(Efes''')
$$
 
By crucial fact $\aleph'_7$,
the expression inside the braces of $(Efes''')$ is a 
$W(\B_k)$- {\it anti-symmetric} polynomial.  By crucial fact $\aleph'_8$,
it is divisible by  $\Delta_k( x_1 , \dots , x_k )$.
The degree in $x_1$ of the polynomial inside the braces of $(Efes''')$ 
is $\leq k-2+k=2k-2$, while the degree of $\Delta_k(x)$, in $x_1$, is
$2k-1$. Hence it must be the zero polynomial. This completes
the proof of subsubsubsublemma $1.1.1.2.2$ \halmos.
 
This completes the proof that
$C_{k}$ satisfies $(C3)$. Finally $(C4)$ is completely routine.
This completes the proof of subsubsublemma $1.1.1.2$. \halmos
 
{\bf Subsubsublemma 1.1.1.3:}  There is a unique sequence of discrete
functions $X_k(n; a_1 , \dots , a_k)$, defined for $k \geq 1$ and
$(n; a_1 , \dots , a_k ) \in \overline{Land\_Of\_Magog_k}$, satisfying
the following partial difference equation:
$$
\left \{ \prod_{i=1}^{k} (I- A_i^{-1}) \right \} 
X_k (n \, ; \, a_1 , \dots , a_k )\,
=\, X_k ( n-1 \, ; \, a_1 , \dots , a_k ) \quad ,
\eqno(P \Delta E (X_k)) 
$$
for $n>k$ and $n \geq a_1 \geq a_2 \geq \dots \geq a_k \geq 1$,
and the following boundary/initial conditions
for $ (n; a_1 , \dots , a_k ) \in \overline {Land\_Of\_Magog_k}$.
 
For $i=1,2, \dots, k-1$,
$$
X_k (n \, ; \, a_1 , \dots , a_k ) = 0 \quad, \qquad \hbox{on} \qquad 
a_i - a_{i+1}=-1 \quad .
\eqno(X1_i)
$$
Also:
$$
X_k (n \, ; \, a_1 , \dots , a_k ) = 0 \quad, \qquad \hbox{on} \qquad 
a_k =0 \quad .
\eqno(X1_k)
$$
Furthermore:
$$
X_k (n \, ; \, n+1 \,,\, a_2 \,,\, \dots \,,\, a_k ) = 0  \quad \hbox{and}
\eqno(X2)
$$
$$
X_k ( k \, ; \, a_1 , \dots , a_{k-1} , a_k ) = 
X_{k-1} ( k \, ; \, a_1 , \dots , a_{k-1} ) \delta ( a_k , 1 )  \quad ,
\eqno(X3)
$$
and
$$
X_1(1; a_1 ) = \delta ( 1, a_1 ) \quad , \quad  (1; a_1 ) \in 
\overline{Land\_Of\_Magog_1} \qquad .
\eqno(X4)
$$
{\eightrm  [ Type `{\eighttt S1113(k,n):}' in ROBBINS, for specific k and n.]}
 
{\bf Proof:} \halmos.
 
Subsublemma $1.1.1$ now follows from subsubsublemmas $1.1.1.1$, 
$1.1.1.2$, and $1.1.1.3$. 
This completes the proof of subsublemma $1.1.1$. \halmos 
 
Subsublemma 1.1.1 gave us a constant term expression 
for the number of $n \times k -$
Magog trapezoids with a prescribed rightmost border. 
In order to complete the proof of sublemma $1.1$, we would have to
sum them all up. In order to establish the constant term expression
$(MagogTotal)$ for $b_k(n)$, we would need to go through an intermediary,
not-so-nice expression, given by subsublemma $1.1.2$ below.
 
{\bf Subsublemma 1.1.2:} The number of $n \times k$ Magog-trapezoids,
$b_k(n)$, is given by the following constant term
expression:
$$
b_k(n)
= CT_{x_1 , \dots , x_k}
\left \{
{
{ x_1 x_2^2 \dots x_k^k \, \cdot \, \Delta_k ( x_1 , \dots , x_k )}
\over
{\prod_{i=1}^{k} (\bar x_i)^{k+n} x_i^{n+k-1} }
} \cdot
{
{1}
\over
{(1- x_k ) (1- x_{k-1}x_k) \dots (1- x_k x_{k-1} \dots x_1 )}
} \right \} \qquad .
\eqno(George)
$$
{\eightrm  [ Type `{\eighttt S112(k,n):}' in ROBBINS, for specific k and n.]}
 
{\bf Proof:} We will establish $(George)$ by starting from the right side,
and proving that it equals the left side. By the well-known symmetry property
of the `$=$' relation, it would follow that the left side is equal to its right 
side.
 
Since $(1-y)^{-1}= \sum_{\alpha =0}^{\infty} y^{\alpha}$,  we have
$$
\displaylines{
(1- x_k)^{-1} (1- x_{k-1}x_k )^{-1} \dots (1- x_k x_{k-1} \dots x_1 )^{-1}
\,=\,
[ \sum_{\alpha_k=0}^{\infty} x_k^{\alpha_k}]
[ \sum_{\alpha_{k-1}=0}^{\infty} (x_{k-1} x_k)^{\alpha_{k-1}}]
\dots
[ \sum_{\alpha_{1}=0}^{\infty} (x_1 \dots x_{k-1} x_k)^
{\alpha_1}]
\cr
=\sum_{ 0 \leq \alpha_1 , \dots , \alpha_k < \infty} 
x_1^{\alpha_1} x_2^{\alpha_1 + \alpha_2}
\dots x_k^{\alpha_1 + \alpha_2 + \dots + \alpha_k} 
=\sum_{0 \leq \lambda_1 \leq \dots \leq \lambda_k < \infty}
 x_1^{\lambda_1} \dots x_k^{\lambda_k} \qquad .
\cr}
$$
 
The right side of $(George)$ is hence equal to:
 
$$
CT_{x_1 , \dots , x_k}
\left \{
{
{ x_1 x_2^2 \dots x_k^k \, \cdot \, \Delta_k ( x_1 , \dots , x_k )}
\over
{\prod_{i=1}^{k} (\bar x_i)^{k+n} x_i^{n+k-1} }
} \cdot
\sum_{0 \leq \lambda_1 \leq \dots \leq \lambda_k < \infty}
 x_1^{\lambda_1} \dots x_k^{\lambda_k}
 \right \} 
$$
$$
= \sum_{0 \leq \lambda_1 \leq \dots \leq \lambda_k < \infty}
CT_{x_1 , \dots , x_k}
\left \{
{
{ x_1 x_2^2 \dots x_k^k \, \cdot \, \Delta_k ( x_1 , \dots , x_k )}
\over
{\prod_{i=1}^{k} (\bar x_i)^{k+n} x_i^{(n- \lambda_i)+k-1} }
} 
 \right \} 
$$
$$
=\sum_{n \geq  (n- \lambda_1) \geq \dots \geq \ (n- \lambda_k ) > - \infty}
CT_{x_1 , \dots , x_k}
\left \{
{
{ x_1 x_2^2 \dots x_k^k \, \cdot \, \Delta_k ( x_1 , \dots , x_k )}
\over
{\prod_{i=1}^{k} (\bar x_i)^{k+n} x_i^{(n- \lambda_i)+k-1} }
} 
 \right \} 
$$
$$
=\sum_{n \geq a_1 \geq \dots \geq \ a_k > - \infty}
CT_{x_1 , \dots , x_k}
\left \{
{
{ x_1 x_2^2 \dots x_k^k \, \cdot \, \Delta_k ( x_1 , \dots , x_k )}
\over
{\prod_{i=1}^{k} (\bar x_i)^{k+n} x_i^{a_i+k-1} }
} 
 \right \}  \quad .
\eqno(George')
$$
 
The sum on the right side of $(George')$ can be replaced by
a sum over the restricted region $n \geq a_1 \geq \dots \geq a_k \geq 1$,
since whenever $a_k <1$, the constant-termand there is a multiple of
$x_k$, and hence its constant term is $0$. So the right side of
$(George')$ (and hence the right side of $(George)$) equals:
$$
\sum_{n \geq a_1 \geq \dots \geq \ a_k \geq 1}
CT_{x_1 , \dots , x_k}
\left \{
{
{ x_1 x_2^2 \dots x_k^k \, \cdot \, \Delta_k ( x_1 , \dots , x_k )}
\over
{\prod_{i=1}^{k} (\bar x_i)^{k+n} x_i^{ a_i +k-1} }
} 
 \right \} 
$$
$$
=\sum_{n \geq a_1 \geq \dots \geq \ a_k \geq 1}
C_k(n; a_1 , \dots , a_k)
=\sum_{n \geq a_1 \geq \dots \geq \ a_k \geq 1}
B_k(n; a_1 , \dots , a_k)= b_k(n)
 \quad ,
\eqno(George'')
$$
where the last equality in $(George'')$ is obvious, and the second-to-last
equality follows from subsublemma $1.1.1$. This completes the proof
of subsublemma $1.1.2$. \halmos
 
By crucial fact $\aleph_4$, we can 
replace the right side of $(George)$
by the constant term of the average of all the $\S_k$-images of its
constant-termand.
Since $\Delta_k$ is antisymmetric in $(x_1 , \dots , x_k )$, we see that
$b_k(n)$ equals
$$
{ {1} \over {k!} }
CT_{x_1 , \dots , x_k} 
\left \{
{
{ \Delta_k ( x_1 , \dots , x_k )}
\over
{\prod_{i=1}^{k} (\bar x_i)^{k+n} x_i^{n+k-1}}
}
\left ( \sum_{\pi \in \S_k} \sgn (\pi) \cdot
\pi \left [ 
{
{x_1 x_2^2 \dots x_k^k}
\over
{(1- x_k)(1- x_k x_{k-1} ) \dots (1 - x_k x_{k-1} \dots x_1 ) }
}
\right ]  \right )  \right  \} \, .
\eqno(Stanley)
$$
 
We now need:
 
{\bf Subsublemma 1.1.3:}
$$
 \sum_{\pi \in \S_k} \sgn (\pi) \cdot
\pi \left [ 
{
{x_1 x_2^2 \dots x_k^k}
\over
{(1- x_k)(1- x_k x_{k-1} ) \cdots (1 - x_k x_{k-1} \dots x_1 ) }
}
\right ]  =
{
{x_1 \cdots x_k \prod_{1 \leq i < j \leq k} ( x_j - x_i )}
\over
{ \prod_{i=1}^k (1- x_i ) \prod_{1 \leq i < j \leq k} (1- x_i x_j ) }
} \quad .
\eqno(Issai)
$$
{\eightrm  [ Type `{\eighttt S113(k);}' in ROBBINS, for specific k.]}
 
{\bf Proof :} See [PS], problem VII.47. Alternatively,
(Issai) is easily seen to be equivalent to Schur's identity
that sums all the Schur functions ([Ma], ex I.5.4, p. 45). 
This takes care of subsublemma $1.1.3$. \halmos
 
Inserting $(Issai)$ into $(Stanley)$,  
expanding ${\prod_{1 \leq i < j \leq k} ( x_j - x_i )}$ by
Vandermonde's expansion, 
$$
\sum_{ \pi \in \S_k} \sgn( \pi) \cdot 
\pi ( x_1^0 x_2^1 \cdots x_k^{k-1} ) \quad ,
$$
using the antisymmetry of $\Delta_k$ once again, and employing crucial 
fact $\aleph_4$, we get the following string of equalities:
 
$$
b_k(n)=
{ {1} \over {k!} }
CT_{x_1 , \dots , x_k}
\left \{
{
{ \Delta_k ( x_1 , \dots , x_k )}
\over
{\prod_{i=1}^{k} (\bar x_i)^{k+n} x_i^{n+k-1}}
}
\left ( 
{
{x_1 \cdots x_k \prod_{1 \leq i < j \leq k} ( x_j - x_i )}
\over
{ \prod_{i=1}^k (1- x_i ) \prod_{1 \leq i < j \leq k} (1- x_i x_j ) }
}
 \right )  \right  \}
$$
$$
=
{ {1} \over {k!} }
CT_{x_1 , \dots , x_k}
\left \{
{
{ \Delta_k ( x_1 , \dots , x_k )}
\over
{\prod_{i=1}^{k} (\bar x_i)^{k+n+1} x_i^{n+k-2} 
\prod_{1 \leq i < j \leq k} (1- x_i x_j )}
}
\left ( 
\sum_{ \pi \in \S_k}  \sgn( \pi) \cdot \pi ( x_1^0 x_2^1 \dots x_k^{k-1} )
 \right )  \right  \}
$$
$$
={ {1} \over {k!} } \sum_{ \pi \in \S_k}
CT_{x_1 , \dots , x_k}
 \left  \{ \pi \left [
{
{ \Delta_k ( x_1 , \dots , x_k )}
\over
{\prod_{i=1}^{k} (\bar x_i)^{k+n+1} x_i^{n+k-2} 
\prod_{1 \leq i < j \leq k} (1- x_i x_j )}
}
\left ( 
\prod_{i=1}^{k} x_i^{i-1}
 \right )  \right ] \right  \}
$$
$$
={ {1} \over {k!} } \sum_{ \pi \in \S_k}
CT_{x_1 , \dots , x_k}
 \left  \{ \pi \left [
{
{ \Delta_k ( x_1 , \dots , x_k )}
\over
{\prod_{i=1}^{k} (\bar x_i)^{k+n+1} x_i^{n+k-i-1} 
\prod_{1 \leq i < j \leq k} (1- x_i x_j )}
} \right ]
  \right  \}
$$
$$
={ {1} \over {k!} } \sum_{ \pi \in \S_k}
CT_{x_1 , \dots , x_k}
\left \{
{
{ \Delta_k ( x_1 , \dots , x_k )}
\over
{\prod_{i=1}^{k} (\bar x_i)^{k+n+1} x_i^{n+k-i-1} 
\prod_{1 \leq i < j \leq k} (1- x_i x_j )}
}
  \right  \}
$$
$$
={ {1} \over {k!} } \left ( \sum_{ \pi \in \S_k} 1 \right )
CT_{x_1 , \dots , x_k}
\left \{
{
{ \Delta_k ( x_1 , \dots , x_k )}
\over
{\prod_{i=1}^{k} (\bar x_i)^{k+n+1} x_i^{n+k-i-1} 
\prod_{1 \leq i < j \leq k} (1- x_i x_j )}
}
  \right  \}
$$
$$
=CT_{x_1 , \dots , x_k}
\left \{
{
{ \Delta_k ( x_1 , \dots , x_k )}
\over
{\prod_{i=1}^{k} (\bar x_i)^{k+n+1} x_i^{n+k-i-1} 
\prod_{1 \leq i < j \leq k} (1- x_i x_j )}
}
  \right  \}
 \, , 
\eqno(George''')
$$
 
where in the last equality we have used Levi Ben Gerson's celebrated
result that the number of elements in $\S_k$ 
(the symmetric group on $k$ elements,)
equals $k!$. The extreme right of $(George''')$ is exactly the
right side of $(MagogTotal)$. 
This completes the proof of sublemma $1.1$. \halmos
 
\vfill
\eject
%%The Alternating Sign Matrix Conjecture Act II
%%File asmII.tex
%begin macros
\baselineskip=14pt
\parskip=10pt
\def \inv{\mathop{\rm inv} \nolimits}
\def \sgn{\mathop{\rm sgn} \nolimits} 
\def\hat{\widehat}
\def\tilde{\widetilde}
\def\epsilon{\varepsilon}
\def\B{{\cal B}}
\def\S{{\cal S}}
\def\halmos{\hbox{\vrule height0.15cm width0.01cm\vbox{\hrule height
 0.01cm width0.2cm \vskip0.15cm \hrule height 0.01cm width0.2cm}\vrule
 height0.15cm width 0.01cm}}
\font\eightrm=cmr8  
\font\eighttt=cmtt8
\magnification=\magstephalf

\parindent=0pt
\overfullrule=0in
\headline={\rm  \ifodd\pageno  \RightHead  \else  \LeftHead  \fi}
\def\RightHead{\centerline{Proof  of the ASM Conjecture-Act II}}
\def\LeftHead{ \centerline{Doron Zeilberger}}
%\pageno=22
%end macros
 
\centerline{\bf Act II. COUNTING GOG}
 
\qquad\qquad\qquad 
{\it And say, Thus saith the Lord God; Behold, I am against thee, O 
Gog $\dots$}
\smallskip
\qquad\qquad\qquad\qquad
\qquad\qquad\qquad\qquad\qquad\qquad \qquad\qquad\qquad\qquad\qquad\qquad 
(Ezekiel XXXIIX,3)
 
An $n \times k -$ Gog trapezoid (where $n \geq k \geq 1$) is
a trapezoidal array of integers
$$
 ( d_{i,j} ) \quad, \quad  \, 1 \leq i \leq n \, 
\quad , \quad  1 \leq j \leq \min (k, n+1-i) \quad ,
$$
such that:
$$
(i) \quad d_{i,j} < d_{i,j+1} \qquad
(ii) \quad d_{i,j} \leq d_{i+1,j} \qquad
(iii) \quad d_{i,j} \geq d_{i+1,j-1} \qquad
$$
$$
(iv) \quad d_{1,j} = j  
\qquad (v) \quad d_{i,k} \leq i+k-1 .
$$ 
(Note that these five conditions imply that  $1 \leq d_{i,j} \leq n$.)
 
{\eightrm  [ To view all n by k Gog trapezoids, for any given k and n,
type `{\eighttt GOG(k,n):}' in ROBBINS. For example,
GOG(3,4) would display all 4 by 3 Gog trapezoids.]}
 
{\bf Sublemma 1.2:} The total number of $n \times k -$ Gog trapezoids, let's
call it $m_k (n)$, is given by the following constant term expression:
$$
m_k(n) \, = \,
CT_{x_1 , \dots , x_k} \left [
 { { \Phi_k ( x_1 , \dots , x_k ) } \over
   { \left \{ \prod_{i=1}^k x_i^n ( \bar x_i)^{n+i+1} \right \} 
\left \{
 \prod_{ 1 \leq i < j \leq k}
(1-x_i x_j ) (1- \bar x_i x_j ) \right \}
}
}
 \right  ] \qquad ,
\eqno(GogTotal)
$$
where the polynomial $\Phi_k(x_1 , \dots , x_k)$ is defined by:
$$
\Phi_k ( x_1 , \dots , x_k ) = 
(-1)^k \,
\sum_{ g \in W( \B_k ) } sgn(g) \cdot
g \left [ \prod_{i=1}^{k} \bar x_i^{k-i} x_i^{k} 
\prod_{1 \leq i < j \leq k} (1- x_i \bar x_j )(1- \bar x_i \bar x_j ) \right ]
\qquad.
\eqno(Gog_1)
$$
{\eightrm  [ Type `{\eighttt S12(k,n):}' in ROBBINS, for specific k and n.]}
 
{\bf Proof:}  
Let $M_k ( n \, ; \, a_1 , a_2 , \dots , a_k )$ be the number of $n \times k$-
Gog trapezoids 
$$
(d_{i,j}), \quad 1 \leq i \leq n \quad ,\quad 
 1 \leq j \leq \min(k,n+1-i) \quad ,
$$
such that $d_{n-k+i,k-i+1} = a_i$, for $i=1, \dots , k$. In other words:
$$
d_{n-k+1,k} = a_1 \quad ,\quad d_{n-k+2, k-1} = a_2 \quad , \quad \dots 
\quad , \quad
d_{n-k+i,k-i+1}=a_i \quad , \quad \dots \quad,\quad d_{n,1}=a_k \quad .
$$
 
The  conceivable $(n; a_1 , \dots , a_k )$ range over the set
$$
\displaylines{
Land\_Of\_Gog_k \, :=
\{ \quad (n \, ; \, a_1 , \dots , a_k )
\quad  | \cr \quad n \geq k \, , 
\quad n \geq a_1 \geq a_2 \geq \dots \geq a_k \geq 1 \quad , \quad
a_1 \geq k \quad , \quad a_2 \geq k-1 \quad , \quad \dots \quad , 
\quad a_k \geq 1 \quad \} \quad . \cr}
$$
{\eightrm  [ To view this set, for specific k and n, 
type `{\eighttt LOGOG(k,n):}' in ROBBINS.]}
 
For any $(n; a_1 , \dots , a_k ) \in Land\_Of\_Gog_k$, define the
set:
$$
\tilde T_k (n \, ; \, a_1 , \dots , a_k ) :=
\{ (b_1 , \dots , b_k ) \quad | \quad 
n-1 \geq b_1 \geq b_2 \geq \dots \geq b_k \geq 1 \quad , \quad
k-i+1 \leq b_i \leq  a_i  \} \qquad .
$$
 
It turns out to be convenient (and most likely necessary) to introduce
a related discrete function that we will call $\tilde M_k$. For $n>k$, 
and $(n \, ; \, a_1 \, , \, \dots \, , \, a_k \,) \in Land\_Of\_Gog_k$,
define:
$$
\tilde M_k ( n \, ; \, a_1 , a_2 , \dots , a_k ) :=
\sum_{b \in \tilde T_k (n \, ; \, a_1 , \dots , a_k )}
M_k ( n-1 \, ; \, b_1 , \dots , b_k ) \qquad ,
\eqno(\tilde M_k)
$$
and when $n=k>1$ {\it define} $\tilde M_k ( k ; k, - , \dots , -)$ by
$$
\tilde M_k ( k \, ; \, k , a_2 , \dots , a_k ) =
\tilde M_{k-1} ( k \, ; \, a_2 , \dots , a_k ) \qquad , ( k > 1 ) \quad ,
\eqno(\tilde M_k(n=k))
$$
while for $n=k=1$ define 
$$
\tilde M_1 (1;1):=1 \quad .
\eqno(\tilde M_k (n=k=1))
$$
 
The following subsublemma gives a constant term expression for
$\tilde M_k$.
 
{\bf Subsublemma 1.2.1:}  Let
$$
H_k(n; a_1 , \dots , a_k ) :=
CT_{x_1 , \dots , x_k} \left \{ 
{ { \Phi_k ( x_1 , ... , x_k ) } 
\over { \left \{ \prod_{i=1}^k x_i^{a_i-1} (\bar x_i)^{n+i} \right \} 
 \prod_{ 1 \leq i < j \leq k}
(1-x_i x_j ) (1- \bar x_i x_j )
} }  \right \} \quad ,
\eqno(H_k)
$$
where $\Phi_k$ is defined in the statement of sublemma $1.2$ above. 
For all 
$(n \, ; \, a_1 \, , \, \dots \, , \, a_k \,) \in Land\_Of\_Gog_k$, we have
$$
\tilde M_k (n\,;\,a_1 , ... , a_k ) =
H_k (n\,;\,a_1 , ... , a_k ) \quad .
\eqno(Gog)
$$
{\eightrm  [ Type `{\eighttt S121(k,n):}' in ROBBINS, for specific k and n.]}
 
{\bf Proof:} It is convenient to extend $\tilde M_k$ 
to the following larger domain (except that we slightly `chop' it
when $n=k$):
$$
\displaylines{
\overline {Land\_Of\_Gog_k} \, :=
\{ \quad (n \, ; \, a_1 , \dots , a_k )
\quad  | n \geq k \, , \,  n- a_1 \geq -1 \, ,
\cr
 a_1 \geq a_2 \geq \dots \geq a_k \quad , \quad
a_1 \geq k-1 \quad , \quad a_2 \geq k-2 \quad , \quad \dots \quad , 
\quad a_k \geq 0  
\cr
\hbox{and} \quad a_1=n+1 \Rightarrow a_2 \leq n \quad \hbox{and} \quad
n=k=a_1 \Rightarrow  a_2<k )
\} \quad, 
\cr}
$$
{\eightrm  [ To view this set, for specific k and n, 
type `{\eighttt ELOGOG(k,n):}' in ROBBINS]},
and to {\it define} $\tilde M_k$ to be $0$ at all the boundary points
that are in $\overline{Land\_Of\_Gog_k} \,\,\,\backslash\, Land\_Of\_Gog_k$.
Note that $(\tilde M_k)$ extends to those boundary points
for which $a_i=k-i$ for one or more $i$'s, since then the
summation-set, $\tilde T_k$, is the empty set. However,
$(\tilde M_k)$ does {\it not} hold for those boundary
points for which $a_1=n+1$.
 
We will prove that subsublemma $1.2.1$  
even holds in this larger domain $\overline {Land\_Of\_Gog_k}$.
This would follow from the following three subsubsublemmas.
 
{\bf Subsubsublemma 1.2.1.1:} Let $k \geq 1$ and for any non-increasing vector
of non-negative integers $a=(a_1 ,\dots , a_k)$, define the partial difference
operator  $P^{(a)}$ by:
$$
P^{(a)} ( A_1 , \dots , A_k ) :=
\prod_{
\{i| a_i - a_{i+1} > 0 \}
}
(I- A_i^{-1}) \qquad,
$$
where we declare that $a_{k+1}=0$. $\tilde M_k$
satisfies the following partial recurrence: 
$$
P^{(a)} (A_1 , \dots , A_k ) \tilde M_k ( n \, ; \, a_1 , \dots , a_k ) 
$$
$$
=
\tilde M_k (n-1 \, ; \, a_1 \,\, ,\,\, \min ( a_1 -1 , a_2 ) \,\, , \,\,
\dots \,\, , \,\, \min (a_{k-1} -1,a_k ) ) \qquad , 
\eqno(P\Delta E(\tilde M_k)) 
$$
valid whenever $n > k$ and $(n; a_1 , \dots , a_k) \in Land\_Of\_Gog_k$.
$\tilde M_k (n; a_1 , \dots , a_k )$ also satisfies the
following initial/boundary conditions for 
$(n \, ; \, a_1 \, , \, \dots \, , \, a_k \,) \in \overline{Land\_Of\_Gog_k}$:
$$
\tilde M_k (n \, ; \, a_1 , \dots , a_k ) = 0  \quad \hbox{on} \quad
a_i = k-i \qquad  ,
\eqno(M1_i)
$$
$$
\tilde M_k (n \, ; \, n+1 , a_2 , \dots , a_k ) = 0 \quad ,
\eqno(M2) 
$$
$$
\tilde M_k ( k \, ; \, k , a_2 , \dots , a_k ) =
\tilde M_{k-1} ( k \, ; \, a_2 , \dots , a_k ) \qquad ,
\eqno(M3) 
$$
and, finally the ``Adam'' condition:
$$
\tilde M_1 ( 1 \, ; \, 1 ) = 1  \qquad .
\eqno(M4)
$$
{\eightrm  [ Type `{\eighttt S1211all(k,n):}' in ROBBINS, 
for specific k and n.]}
 
{\bf Proof:} Consider a typical $n \times k$ Gog trapezoid $(d_{i,j})$
that is counted by $M_k(n; a_1 , \dots , a_k)$, where $n>k$ and
$(n \, ; \, a_1 \, , \, \dots \, , \, a_k \,) \in Land\_Of\_Gog_k$.
This means that 
$$
d_{n-k+1,k} = a_1 \quad ,\quad d_{n-k+2, k-1} = a_2 \quad , \quad \dots 
\quad , \quad
d_{n-k+i,k-i+1}=a_i \quad , \quad \dots \quad,\quad d_{n,1}=a_k \quad .
$$
Now let 
$$
d_{n-k,k} = b_1 \quad , \quad d_{n-k+1,k-1}=b_2 \quad , \quad
\dots \quad , \quad d_{n-k+i-1,k-i+1}
=b_i \quad , \quad \dots \quad , \quad d_{n-1,1}=b_k \qquad. 
$$
 
By the conditions $(i)-(v)$ defining 
Gogs, given right at the beginning of this section,
it follows that the allowed $b= (b_1 , \dots , b_k )$ range over the set
$$
T_k (n \, ; \, a_1 , \dots , a_k ) :=
$$
$$ 
\{ (b_1 , \dots , b_k ) \quad | \quad
n-1 \geq b_1 \geq b_2 \geq \dots \geq b_k \geq 1 \quad \hbox{and} \quad
k-i+1 \leq b_i \leq \min ( a_i , a_{i-1}-1 ) \quad (i=1, \dots , k) \,
 \} \quad ,
$$
where $a_0 := \infty$. By deleting $d_{n,1} = a_k \quad,\quad
 d_{n-1,2}=a_{k-1} \quad ,\quad \dots \quad ,\quad
d_{n-k+1,k}=a_1 \,$, we obtain a certain $(n-1) \times k -$ Gog trapezoid.
Looking at all the conceivable possibilities for $b= (b_1 , \dots , b_k )$,
we get the following (non-local) recurrence, valid for $n > k$, 
and $(n \, ; \, a_1 \, , \, \dots \, , \, a_k \,) \in Land\_Of\_Gog_k \,$,
$$
M_k (n \, ; \, a_1 , \dots , \, , a_k ) =
\sum_{b \in T_k (n \, ; \, a_1 , \dots , a_k ) }
M_k ( n-1 \, ; \, b_1 , \dots , b_k ) \qquad .
\eqno(Howard)
$$
 
(Howard) Can be used to compile a table of $M_k$ 
together with the obvious initial condition
$$
M_k ( k \, ; \, k , a_2 , \dots , a_k ) =
M_{k-1} ( k \, ; \, a_2 , \dots , a_k ) \qquad , ( k > 1 ) \quad ,
\eqno(Howard')
$$
(note that, when $n=k$, $a_1$ {\it must} be $k$,) and the ``Adam'' condition
$M_1 (1 \, ; \, 1)=1$.
 
Using the definition $(\tilde M_k)$, 
Eq. (Howard) can be rewritten as:
$$
M_k ( n \, ; \, a_1 , \dots , a_k ) =
\tilde M_k (n \, ; \, a_1 \, ,\, \min (a_1 -1 , a_2 )\,  , \, \dots \, , \,
\min ( a_{k-1}-1 , a_k ) ) \qquad ,
\eqno(Bill)
$$
valid for $n>k$ and $(n; a_1 , \dots , a_k) \in Land\_Of\_Gog_k$. 
We are now ready to prove that $\tilde M_k$ satisfies 
($P\Delta E(\tilde M_k)$) in the statement of the current
subsubsublemma $1.2.1.1$.
Note that since $a_k >0$, the factor $(I- A_k^{-1})$ is
always present in $P^{(a)}$.
 
Let $(a_1 , \dots , a_k)$ be arranged in maximal blocks of equal
components as follows:
$$
a_1 = a_2 = \dots = a_{r_1} > a_{r_1 +1}= \dots = a_{r_2}
> \dots >  a_{r_{l-1} +1} = \dots =
a_{r_l} \geq 1 \qquad ,
$$
where $r_l=k$, and $r_1$, $r_2 - r_1 , \dots , r_l - r_{l-1}$, are
the lengths
of (maximal) blocks of consecutive equal components of $a$.
Let's put 
$$
c_1:=a_{r_1} \, , \, \dots \, , \, c_{l-1}:=a_{r_{l-1}}
\, , \, c_l:=a_{r_l} (=a_k) \quad.
$$ 
We have:
$$
P^{(a)} = (I- A_{r_1}^{-1}) \dots (I- A_{r_l}^{-1}) \quad .
$$
 
The recurrence that we have to show,
$(P \Delta E( \tilde M_k))$, spells out to:
$$
\displaylines{
(I-A_{r_1}^{-1}) (I- A_{r_2}^{-1}) \dots (I-A_{r_l}^{-1})
\tilde M_k (n \, ; \, c_1 , \dots , c_1 ,\, c_2, \dots , c_2, \,
\dots \,,\,\dots , c_l ,\dots , c_l ) 
\cr
=
\tilde M_k (n-1 ; c_1 , c_1-1 , \dots , c_1 -1 ,
c_2, c_2-1, \dots , c_2-1, \dots , \dots , c_l , c_l -1, \dots , c_l -1 )
\quad , \cr}
$$
where the length of the `$c_i$' block is $r_{i}-r_{i-1}$, for
$i=1 \dots l$ (where we agree that $r_0=0$.)
 
Suppose first that $c_1=a_1<n$, so that 
$(n-1 ; \, a_1 \, , \, \dots \, , \, a_k \,) \in Land\_Of\_Gog_k$.
 
Applying $(I-A_k^{-1})$ (note that $k=r_l$) to the sum $(\tilde M_k)$,
yields a sum over the subset:
$$
\displaylines{
\tilde T_k (n; a_1 , \dots , a_{k-1} , a_k ) \,\, \backslash \,\,
\tilde T_k (n; a_1 , \dots , a_{k-1} , a_k-1 )=
\cr
\tilde T_k (n; c_1 , \dots , c_1 , c_2 , \dots , c_2 ,
\dots , \dots , c_{l-1} , \dots , c_{l-1}, c_l , \dots , c_l,c_l ) 
\,\, \backslash \,\,
\cr
\tilde T_k (n; c_1 , \dots , c_1 , c_2 , \dots , c_2 ,
\dots , \dots , c_{l-1} , \dots , c_{l-1}, c_l , \dots , c_l,c_l-1 )
\cr
=
\{ (b_1 , \dots , b_k ) \quad | \quad 
n-1 \geq b_1 \geq b_2 \geq \dots \geq b_k \geq 1 \quad , \quad
\cr
k-i+1 \leq b_i \leq  a_i  \quad \hbox{for} \quad i=1, \dots , r_{l-1} \, , \,
\hbox{and} \quad
b_i=c_l \, , \, \hbox{for} \quad i=r_{l-1}+1, \dots ,  k \} \qquad .\cr}
$$
 
The last equality follows from the fact that 
$$
b_{r_{l-1}+1} \geq \dots \geq b_{r_l}=c_l \quad \hbox{and} \quad
b_{r_{l-1}+1}  \leq c_l \quad ,
$$
implies that $b_{r_{l-1}+1}= \dots = b_{r_l}=c_l$.
Applying next $(I-A_{r_{l-1}}^{-1})$ to this smaller sum, yields a sum
with the same summand, but over the smaller subset:
$$
\displaylines{
\{ (b_1 , \dots , b_k ) \quad | \quad 
n-1 \geq b_1 \geq b_2 \geq \dots \geq b_k \geq 1 \quad , \quad
k-i+1 \leq b_i \leq  a_i  \, , \quad \hbox{for} \quad i=1, \dots , r_{l-2}
\quad,
\cr
b_i=c_{l-1} \, , \, \hbox{for} \quad i=r_{l-2}+1 \dots r_{l-1} \, , \,
\hbox{and} \quad
b_i=c_l \, , \, \hbox{for} \quad i=r_{l-1}+1 , \dots ,  k \} \qquad . 
\cr}
$$
 
Continuing to apply 
$(I-A_{r_i}^{-1})$ with $i=l-2  \, , \, l-3 \,,\,\dots \, ,\, 1$,  we keep
getting a sum over continuously shrinking subsets of $\tilde T_k$, until
at the end we get the sum of $M_k (n-1; b_1  , \dots , b_k )$ over
the singleton set 
$$
\{ (n-1; b_1 , \dots , b_k )= (n-1; a_1 , \dots , a_k ) \} \quad.
$$
 
We have thus shown that, whenever $n >k$ and 
$(n; a_1 , \dots , a_k) \in Land\_Of\_Gog_k$, we have
$$
P^{(a)} (A_1 , \dots , A_k ) \tilde M_k ( n \, ; \, a_1 , \dots , a_k )
=
 M_k (n-1 \, ; \, a_1  ,  \dots , a_k )
\qquad .
$$
 
Now use $(Bill)$, with $n$ replaced by $n-1$.
($(Bill)$ holds with $n$ replaced by $n-1$ if $n-1>k$ and 
$(n-1;a_1,...,a_k) \in Land\_Of\_Gog_k$.
For this we need the assumed hypothesis $a_1<n$.)
 
This leaves those elements $(n;a_1,...,a_k)$ in $Land\_Of\_Gog_k$, $n>k$, 
for which $(n-1;a_1,...,a_k)$ does not lie in $Land\_Of\_Gog_k$, namely
the cases $n-1=k$ or $a_1=n$. These are done in the next two paragraphs.
 
When $n=k+1$, the sole surviving term, after applying $P^{(a)}$
to $(\tilde M_k)$, is $M_k(k;k,a_2, \dots ,a_k)$,
which by $(Howard')$, equals $M_{k-1}(k;a_2 , \dots , a_k)$,
which by $(Bill)$ equals
$\tilde M_{k-1} (k; a_2 , \min(a_2-1,a_3), \dots , \min(a_{k-1}-1,a_k) )$,
which, in turn, by the {\it definition} $(\tilde M_k (n=k))$ equals
$\tilde M_{k} (k;k, a_2 , \min(a_2-1,a_3), \dots , \min(a_{k-1}-1,a_k) )$, 
which disposes of the case $n=k+1$.
 
Now suppose that $a_1=n$, so that 
$a_1 = \dots = a_{r_1}=n$. As before, applying
$(I- A_{r_2}^{-1}) \dots (I- A_{r_l}^{-1})$ to the sum on the
right of Eq. $(\tilde M_k)$, leaves us a much reduced sum over
the subset of $\tilde T_k (n; a_1 , \dots , a_k)$ for which
$b_{r_1+1}=a_{r_1+1} , \dots , b_{k}=a_{k}$. Since we {\it must}
have $b_1  \leq n-1$ this is the same set as the one for 
which $a_1= \dots =a_{r_1-1}=n$ and  $a_{r_1}=n-1$, so that applying
$(1-A_{r_1}^{-1})$ to that sum gives $0$. But the right side 
of $(P \Delta E (\tilde M_k))$ is also $0$ then,
by the extended  definition of $\tilde M_k$, since $a_1=n=$`$n-1$'$+1$,
and $n=$`$n-1$'.
 
$(M1_i)$  and $(M2)$ hold  because $\tilde M_k$ was {\it defined} to be
$0$ there,
while $(M3)$ is the definition at $n=k$, given by
$(\tilde M_k (n=k))$,
and $(M4)$ follows from the definition ($\tilde M_1$).
This completes the proof of $(sub)^3$lemma $1.2.1.1$.
\halmos
 
{\bf Subsubsublemma 1.2.1.2:} Let $k \geq 1$, and for any non-increasing vector
of non-negative integers $a=(a_1 ,\dots , a_k)$, define the partial difference
operator  $P^{(a)}$ by:
$$
P^{(a)} ( A_1 , \dots , A_k ) :=
\prod_{
\{i| a_i - a_{i+1} > 0 \}
}
(I- A_i^{-1}) \qquad,
$$
where we declare that $a_{k+1}=0$. The sequence of discrete functions 
$H_k(n \, ; \,  a_1 , \dots , a_k)$,  ($n \geq k \geq 1$), defined
in the statement of the subsublemma $1.2.1$,
satisfies the following partial recurrence: 
$$
P^{(a)} (A_1 , \dots , A_k ) H_k ( n \, ; \, a_1 , \dots , a_k )
$$
$$
=
H_k (n-1 \, ; \, a_1 \,\, ,\,\, \min ( a_1 -1 , a_2 ) \,\, , \,\,
\dots \,\, , \,\, \min (a_{k-1} -1,a_k ) ) \qquad ,
\eqno(P\Delta E(H_k))
$$
valid whenever $n > k$ and $(n; a_1 , \dots , a_k) \in Land\_Of\_Gog_k$.
$H_k (n; a_1 , \dots , a_k )$ also satisfies the
following initial/boundary conditions for 
$(n \, ; \, a_1 \, , \, \dots \, , \, a_k \,) \in \overline{Land\_Of\_Gog_k}$:
$$
H_k (n \, ; \, a_1 , \dots , a_k ) = 0  \quad \hbox{on} \quad
a_i = k-i \qquad  ,
\eqno(H1_i)
$$
$$
H_k (n \, ; \, n+1 , a_2 , \dots , a_k ) = 0 \quad ,
\eqno(H2) 
$$
$$
H_k ( k \, ; \, k , a_2 , \dots , a_k ) =
H_{k-1} ( k \, ; \, a_2 , \dots , a_k ) \qquad ,
\eqno(H3) 
$$
and, finally the ``Adam'' condition
$$
H_1 ( 1 \, ; \, 1 ) = 1  \qquad .
\eqno(H4)
$$
 
{\bf Proof:} The task of proving this $(sub)^3$lemma will be divided
between the following $(sub)^4$lemmas: $1.2.1.2.1$, $1.2.1.2.2$,
$1.2.1.2.3$, and $1.2.1.2.4$, and $1.2.1.2.5$, which will prove
$(P \Delta E(H_k))$, $(H1_i)$, $(H2)$, $(H3)$, and $(H4)$ respectively.
 
{\bf Subsubsubsublemma 1.2.1.2.1:} Let $H_k(n; a_1 , \dots , a_k)$ be
the discrete function defined in the statement of subsublemma $1.2.1$,
i.e. :
$$
H_k(n; a_1 , \dots , a_k ) :=
CT_{x_1 , \dots , x_k} \left \{ 
{ { \Phi_k ( x_1 , ... , x_k ) } 
\over { \left \{ \prod_{i=1}^k x_i^{a_i-1} (\bar x_i)^{n+i} \right \} 
 \prod_{ 1 \leq i < j \leq k}
(1-x_i x_j ) (1- \bar x_i x_j )
} }  \right \} \quad ,
$$
where $\Phi_k$ is the {\it anti-symmetric} polynomial 
defined in the statement of sublemma $1.2$. (In fact all we need is
that $\Phi_k$ is {\it some} anti-symmetric polynomial.)
For all $n>k$, and
$(n \, ; \, a_1 \, , \, \dots \, , \, a_k \,) \in Land\_Of\_Gog_k$
(in fact whenever  $a_1 \geq  \dots \geq a_k$,) we have
that the partial-recurrence equation $(P \Delta E (H_k))$,
given in the statement of the parent $(sub)^3$lemma $1.2.1.2$, holds.
 
{\eightrm  [ Type `{\eighttt S12121all(k,n):}' in ROBBINS, 
for specific k and n.]}
 
{\bf Proof:} Suppose that we have:
$$
a_1 = a_2 = \dots = a_{r_1} > a_{r_1 +1}= \dots = a_{r_2}
> \dots >  a_{r_{l-1} +1} = \dots =
a_{r_l} \geq 1 \qquad ,
$$
where $r_l=k$, and $r_1$, $r_2 - r_1 , \dots , r_l -r_{l-1}$, are
the lengths of maximal blocks of equal components of $a=(a_1 ,\dots,a_k)$.
Putting $c_i:=a_{r_i}$, for $i=1 , \dots , l$,
we have to show that
$$
(I-A_{r_1}^{-1}) (I- A_{r_2}^{-1}) \dots (I-A_{r_l}^{-1})
H_k (n \, ; \, c_1 , \dots , c_1 ,\, c_2, \dots , c_2, \,
\dots \,,\,\dots , c_l ,\dots , c_l ) 
$$
$$
=
H_k (n-1 ; c_1 , c_1-1 , \dots , c_1 -1 ,
c_2, c_2-1, \dots , c_2-1, \dots , \dots , c_l , c_l -1, \dots , c_l -1 )
\quad , 
\eqno(Rodica)
$$
where, 
in the left, the block of $c_1$ is $r_1$-component-long,
the block of $c_2$ is $r_2$- component-long, $\dots$, the block
of $c_l$ is $r_l$-component long. In the right side, the $c_i-1$
block is $r_i-1$ -component long, for $i=1, \dots , l$.
 
By crucial fact $\aleph_3$, 
the difference between the left and right sides of $(Rodica)$ equals:
$$
CT_{x_1 , \dots , x_k}
 \left [ 
 { 
  { \{ \bar x_{r_1} \bar x_{r_2} \dots \bar x_{r_l} \quad - \quad
[ \bar x_1 \bar x_2 \dots \bar x_k ]
( x_2  \dots x_{r_1}) (x_{r_1+2} \dots x_{r_2})  \dots
( x_{r_{l-1}+2} \dots x_{r_l} )
 \} 
\Phi_k ( x_1 , \dots , x_k ) } 
     \over
   {
\left \{ \prod_{i=1}^k x_i^{a_i-1} ( \bar x_i)^{n+i} \right \}
\prod_{ 1 \leq i < j \leq k}
(1-x_i x_j ) (1- \bar x_i x_j )
   } 
  }
 \right ] \qquad 
\eqno(Dave)
$$
In order to prove $(sub)^4$lemma $1.2.1.2.1$, we need to show that
the constant-term expression in $(Dave)$ is identically zero.
To this end we need the following $(sub)^5$ lemma:
 
{\bf Subsubsubsubsublemma 1.2.1.2.1.1:} Let $Jamie(x_1 ,\dots , x_k)$
be the polynomial in $(x_1 , \dots ,x_k)$ inside the braces of
$(Dave)$, i.e. the polynomial:
$$
Jamie(x_1 , \dots , x_k):=
\prod_{j=1}^l \bar x_{r_j} -
\left \{ \prod_{i=1}^{k} \bar x_i \right \} 
\left \{ \prod_{j=1}^{l} \prod_{i=r_{j-1}+2}^{r_{j}} x_i \right \}
\quad.
$$
 
$Jamie(x_1 , \dots , x_k)$ can be written as follows
(recall that $r_l=k$ and we put $r_0=0$):
$$
Jamie(x_1 , \dots , x_k)=
\sum_{ 
{{p=2} \atop { p \not \in \{ r_1 +1 ,  r_2+ 1 , \dots , r_{l-1}+1 \} } }
     }^{k}
POL( \hat x_{p-1} , \hat x_p ) \, \cdot \, \bar x_p \, (1- \bar x_{p-1} x_p ) ,
\eqno(Herb)
$$
where $POL( \hat x_{p-1} , \hat x_p )$ is shorthand for:
``some polynomial of $(x_1 , \dots , x_k)$ that does {\it not} depend
on the variables $x_{p-1}$ and $x_p$''. (In other words,
$POL( \hat x_{p-1} , \hat x_p )$ is some polynomial in the $k-2$ variables
$(x_1 , \dots , x_{p-2} , x_{p+1} , \dots , x_k)$.)
 
{\bf Proof:}  Since (let $\bar x_0 :=1$)
$$
\prod_{i=1}^{k} \bar x_i =
\left \{ \prod_{j=1}^{l} \prod_{i=r_{j-1}+1}^{r_{j}} \bar x_{i-1} \right \}
\cdot \bar x_k
=
\left \{ \prod_{j=1}^{l} \prod_{i=r_{j-1}+2}^{r_{j}} \bar x_{i-1} \right \}
\cdot 
\left \{ \prod_{j=1}^{l} \prod_{i=r_{j-1}+1}^{r_{j-1}+1} \bar x_{i-1} \right \}
\cdot \bar x_k
$$
$$
=
\left \{ \prod_{j=1}^{l} \prod_{i=r_{j-1}+2}^{r_{j}} \bar x_{i-1} \right \}
\cdot 
\left \{ \prod_{j=1}^{l} \bar x_{r_j} \right \} \quad ,
$$
we have that
$$
Jamie(x_1 , \dots , x_k)=
\prod_{j=1}^l \bar x_{r_j} -
\left \{ \prod_{i=1}^{k} \bar x_i \right \} 
\left \{ \prod_{j=1}^{l} \prod_{i=r_{j-1}+2}^{r_{j}} x_i \right \}
$$
$$
=
\prod_{j=1}^l \bar x_{r_j} -
\left \{ \prod_{j=1}^{l} \prod_{i=r_{j-1}+2}^{r_{j}} ( \bar x_{i-1} x_i )
\right \}
\cdot 
\prod_{j=1}^{l} \bar x_{r_j} =
\left \{ \prod_{j=1}^{l} \bar x_{r_j} \right \} \cdot
\left [ 1-  \prod_{j=1}^{l} \,\,\prod_{i=r_{j-1}+2}^{r_{j}} ( \bar x_{i-1} x_i )
\right ] \quad .
\eqno(Marvin)
$$
 
We now need the following $(sub)^6$ lemma:
 
{\bf Subsubsubsubsubsublemma 1.2.1.2.1.1.1:} Let $U_j$, $j=1, \dots , l$, be
quantities in an associative algebra, then:
$$
1- \prod_{j=1}^{l} U_j =
\sum_{j=1}^l \left \{ \prod_{h=1}^{j-1} U_h \right \} (1- U_j) \quad .
$$
 
{\bf Proof:} The series on the right  telescopes to the expression
on the left. Alternatively, use increasing induction on $l$, starting
with the tautologous ground case $l=0$. \halmos
 
Using $(sub)^6$lemma $1.2.1.2.1.1.1$ with
$$
U_j = \prod_{i=r_{j-1}+2}^{r_j} ( \bar x_{i-1} x_i ) \quad ,
$$
we get that  $(Marvin)$ implies:
$$
Jamie(x_1 , \dots , x_k)=
\left \{ \prod_{m=1}^{l} \bar x_{r_m} \right \} \cdot
\sum_{j=1}^l \left \{ \prod_{h=1}^{j-1}
\prod_{i=r_{h-1}+2}^{r_h} ( \bar x_{i-1} x_i ) 
\right \} \cdot
\left ( \,\, 
1-  \prod_{i=r_{j-1}+2}^{r_j} ( \bar x_{i-1} x_i ) \,\, \right ) \quad.
\eqno(Marvin')
$$
 
We can split $(Marvin')$ yet further apart, with the aid of the
following $(sub)^6$lemma:
 
{\bf Subsubsubsubsubsublemma 1.2.1.2.1.1.2:} Let $U_j$, ($j=K, \dots , L$), be
quantities in an associative algebra, then:
$$
1- \prod_{i=K}^{L} U_i =
\sum_{p=K}^L (1- U_p) \left \{ \prod_{h=p+1}^L U_h \right \}  \quad .
$$
 
{\bf Proof:} The sum on the right telescopes to the expression on
the left. (Note that it is in the opposite direction to the way
in which it happened in $1.2.1.2.1.1.1$.) Alternatively,
the identity is tautologous when $K=L+1$, and follows by decreasing
induction on $K$.  This completes the proof of $(sub)^6$ lemma
$1.2.1.2.1.1.2$. \halmos .
 
Going back to $(Marvin')$, we use the last
$(sub)^6$lemma ($1.2.1.2.1.1.2$), with $K=r_{j-1}+2$, $L=r_j$,
and $U_i:= ( \bar x_{i-1} x_i )$, to rewrite:
$$
Jamie(x_1 , \dots , x_k)=
\left \{ \prod_{j=1}^{l} \bar x_{r_j} \right \} \cdot
\sum_{j=1}^l 
\,\, \left \{ \prod_{h=1}^{j-1}
\prod_{i=r_{h-1}+2}^{r_h} ( \bar x_{i-1} x_i ) \right \} \cdot
\sum_{p=r_{j-1}+2}^{r_j} (1- \bar x_{p-1} x_p )
 \prod_{i=p+1}^{r_j} ( \bar x_{i-1} x_i )   
$$
$$
=
\sum_{j=1}^l  \sum_{p=r_{j-1}+2}^{r_j}  \,\,
\left \{ \prod_{m=1}^{l} \bar x_{r_m} \right \}
\left \{ \prod_{h=1}^{j-1}
\prod_{i=r_{h-1}+2}^{r_h} ( \bar x_{i-1} x_i ) 
 \right \} \cdot
 (1- \bar x_{p-1} x_p )
\left \{  \prod_{i=p+1}^{r_j} ( \bar x_{i-1} x_i )   \right \}
\quad .
\eqno(Marvin'')
$$
 
The right side of $(Marvin'')$ is a sum over all $p$ in the
range $1 \leq p \leq k$
except that $p=r_j+1$, for $j=0, \dots , l-1$, are  omitted (where
$r_0:=0$.) For each such participating $p$, the summand is
$(1- \bar x_{p-1} x_p)$ times a product of $x_i$'s and $\bar x_i$'s,
neither of which is $x_p, \bar x_p , x_{p-1}, \bar x_{p-1}$ {\it except}
for a single $\bar x_p$, which comes out of the first term of
the product $ \prod_{i=p+1}^{r_j} ( \bar x_{i-1} x_i )$ in
$(Marvin'')$, when $p< r_j$ for some $j$, or from the $j$'th term
of the product $\prod_{j=1}^{l} \bar x_{r_j}$ in case $p=r_j$ for
one of the $j$'s. Thus the polynomial $Jamie(x_1 , \dots , x_k)$
can indeed be expressed as claimed in $(Herb)$. This completes the
proof of $(sub)^5$ lemma $1.2.1.2.1.1$. \halmos
 
Going back to the proof of $(sub)^4$lemma $1.2.1.2.1$, we insert
Eq. $(Herb)$ into Eq. $(Dave)$, to get that the expression in $(Dave)$
equals:
$$
\sum_{ 
{{p=2} \atop { p \not \in \{ r_1 +1 ,  r_2+ 1 , \dots , r_{l-1}+1 \} } }
     }^{k}
CT_{x_1 , \dots , x_k} \left \{ 
{ { 
POL( \hat x_{p-1} , \hat x_p ) \, \cdot \, \bar x_p \, (1- \bar x_{p-1} x_p )
\Phi_k ( x_1 , ... , x_k ) } 
\over { \left \{ \prod_{i=1}^k x_i^{a_i-1} (\bar x_i)^{n+i} \right \} 
 \prod_{ 1 \leq i < j \leq k}
(1-x_i x_j ) (1- \bar x_i x_j )
} }  \right \} .
\eqno(Dave')
$$
 
Let $RAT( \hat x_{p-1} , \hat x_p )$ denote a rational function
(that possesses a Laurent series) of the $k-2$ variables
$(x_1 , \dots , x_{p-2},x_{p+1}, \dots , x_k)$. We can express the
constant-termand of the summand of $(Dave')$ as follows:
$$
{ {POL( \hat x_{p-1} , \hat x_p ) } 
\over { \left \{ \prod_
{{{i=1} \atop {i \neq p-1,p}}}^k x_i^{a_i-1} (\bar x_i)^{n+i} \right \} 
 \prod_{ {{1 \leq i < j \leq k} \atop {i,j \neq p-1,p}} }
(1-x_i x_j ) (1- \bar x_i x_j )
} } \cdot
$$
$$
{{ \bar x_p \, (1- \bar x_{p-1} x_p )
\Phi_k ( x_1 , ... , x_k )} \over
{ x_{p-1}^{a_{p-1}-1} x_p^{a_p-1} \bar x_{p-1}^{n+p-1} \bar x_p^{n+p}
 (1-x_{p-1} x_p ) (1- \bar x_{p-1} x_p )
\prod_{i=1}^{p-2} (1- x_i x_{p-1})(1- \bar x_i x_{p-1})
(1- x_i x_{p})(1- \bar x_i x_{p})
}} \cdot
$$
$$
{{1} \over {
\prod_{j=p+1}^{k} (1-  x_{p-1} x_j )(1- \bar x_{p-1} x_j)
(1-  x_{p} x_j )(1- \bar x_{p} x_j)
}}
$$
$$
=
{{RAT( \hat x_{p-1} , \hat x_p ) \Phi_k ( x_1 , ... , x_k )} \over
{ x_{p-1}^{a_{p-1}-1} x_p^{a_p-1} \bar x_{p-1}^{n+p-1} \bar x_p^{n+p-1}
 (1-x_{p-1} x_p )
\prod_{i=1}^{p-2} (1- x_i x_{p-1})(1- \bar x_i x_{p-1})
(1- x_i x_{p})(1- \bar x_i x_{p})}} \cdot
$$
$$
{{1} \over {
\prod_{j=p+1}^{k} (1-  x_{p-1} x_j )(1- \bar x_{p-1} x_j)
(1-  x_{p} x_j )(1- \bar x_{p} x_j) 
}} \quad .
$$
 
But this constant-termand is manifestly anti-symmetric w.r.t.
$x_{p-1} \leftrightarrow x_p$, since $\Phi_k$ is, and the rest is
symmetric, since $a_{p-1}=a_p$ (recall that $p \neq r_j +1$,
$j=0, \dots , l-1$.) It follows from crucial fact $\aleph_1$
that its constant-term, $CT_{x_1, \dots , x_k}$ vanishes.
Since every single term in $(Dave')$ vanishes, it follows
that $(Dave')$, and hence $(Dave)$ must also vanish. This
completes the proof of $(sub)^4$lemma $1.2.1.2.1$. \halmos
 
{\bf Subsubsubsublemma 1.2.1.2.2:} For $i=1, \dots , k$, we have, for
$(n; a_1 , \dots , a_k ) \in \overline {Land\_Of\_Gog_k}$,
$$
H_k (n \, ; \, a_1 , \dots , a_k ) = 0 \quad on \quad
a_i=k-i \quad .
\eqno(H1)_i
$$
{\eightrm  [ Type `{\eighttt S12122(k,n):}' in ROBBINS, for specific 
k and n.]}
 
{\bf Proof:} By the definition of $H_k$,
$$
\displaylines{
H_k(n; a_1 , \dots , a_{i-1} , k-i , a_{i+1} , \dots , a_k )=
\cr
CT_{x_1, \dots , x_k} \left [
 { { \Phi_k ( x_1 , \dots , x_k ) } \over
   { x_1^{a_1 -1} \cdots x_{i-1}^{a_{i-1} -1} x_i^{k-i-1} \cdots x_k^{a_k -1} }
} \cdot
\prod_{i=1}^k ( \bar x_i)^{-n-i} \prod_{ 1 \leq i < j \leq k}
(1-x_i x_j )^{-1} (1- \bar x_i x_j )^{-1}   \right ] \qquad .
\cr}
$$
 
For each term of the polynomial $\Phi_k$ (see $(Gog_1)$ in the statement of
sublemma $1.2$ for its definition), there
are at least $i$ distinct $x$'s that appear as factors with multiplicity
$\geq k-i$. But the denominator of the constant termand has
at least $k-i+1$ $x$'s in the denominator with 
exponent $\leq k-i-1$. 
By the pigeonhole principle, these two sets must have a non-empty 
intersection. Say $x_l$ appears in this intersection. Then, regarded as
a Laurent series in $x_l$, our constant termand is a positive power of $x_l$
times a power series in $x_l$. Thus its constant term with respect to
$x_l$ is zero. This completes the proof of $(sub)^4$lemma
$1.2.1.2.2$.  \halmos
 
{\bf Subsubsubsublemma 1.2.1.2.3:} For
$(n; n+1 ,  a_2 , \dots , a_k ) \in \overline {Land\_Of\_Gog_k \,\,}$,
(note that $a_2 \leq n$,)
$$
H_k( n; n+1 , a_2 , \dots , a_k ) =0 \quad .
\eqno(H2)
$$
{\eightrm  [ Type `{\eighttt S12123(k,n):}' in ROBBINS, for specific k and n.]}
 
{\bf Proof:}  By definition,
$$
H_k( n; n+1 , a_2 , \dots , a_k ) =
$$
$$
CT_{x_1 , \dots , x_k} \left [
 { { \Phi_k ( x_1 , \dots , x_k ) } \over
   { x_1^{n} x_2^{a_2 -1} \dots x_k^{a_k -1} }
} \cdot
 \bar x_1^{-n-1} \bar x_2^{-n-2} \dots \bar x_k^{-n-k}
 \prod_{ 1 \leq i < j \leq k}
(1-x_i x_j )^{-1} (1- \bar x_i x_j )^{-1}  \right ] \, .
\eqno(H2a)
$$
 
$(H2)$ is obvious for $a_2=a_3= \dots = a_k =0$.
It is reasonable to try
induction on the size of the $a_i$, or more precisely, on
their sum. As is often the case when trying to apply induction,
we would have to consider points that we do not care about. So let's
invite all the vectors $( a_2 , \dots , a_k)$ for which
the components are all $\leq n$, to be included in the statement
of subsubsubsublemma 1.2.1.2.3. Never mind that $\tilde M_k$ doesn't make
sense unless $a_2 \geq a_3 \geq \dots \geq a_k$, right now we are
talking about $H_k$, for which it is meaningful to consider
all such $(a_2 , \dots , a_k)$, and we are thus lead to prove the more general
statement that $(H2a)$ vanishes for all $a_2 , \dots , a_k \leq n$.
By multiplying the numerator by
$$
\bar x_2^{n+2} \dots \bar x_k^{n+k} \prod_{2 \leq i < j \leq k}
(1- x_i x_j ) (1- \bar x_i x_j )= 1 + O(x_2) + \dots +O(x_k) \quad ,
$$
and using the inductive hypothesis on $a_2 + \dots + a_k$, it
is clear that  the assertion that
$(H2a)$ vanishes for all $a_2 , \dots , a_k \leq n$ 
is equivalent to the following simpler statement:
 
{\bf Subsubsubsubsublemma 1.2.1.2.3.1}: For all $n \geq k$, and
$a_2 , \dots , a_k \leq n$, the following identity holds.
(Recall that for any variable $t$, its {\it bar}, $\bar t$,
stands for $1-t$.)
$$
F_k( n;  a_2 , \dots , a_k ):=
CT_{x_1 , \dots , x_k} \left [
 { { \Phi_k ( x_1 , \dots , x_k ) } \over
   { x_1^{n}   \bar x_1^{n+1} x_2^{a_2 -1} \dots x_k^{a_k -1} }
} \cdot
 \prod_{i=2}^{k}
(1-x_1 x_i )^{-1} (1- \bar x_1 x_i )^{-1}  \right ] \,= \,0\, , 
\eqno(H2b)
$$
where $\Phi_k ( x_1 , \dots , x_k )$ is the polynomial defined
in $(Gog_1)$ (in the statement of sublemma $1.2$.)
 
{\bf Proof:}  We can write $F_k (n; a_2 , \dots , a_k )$ as
$$
F_k( n;  a_2 , \dots , a_k ):=
CT_{x_1} \left [
 { { P( x_1 ) } \over
   { x_1^{n}   \bar x_1^{n+1} }
} \right ] \quad ,
$$
where, $P( x_1 )$ is the polynomial in $x_1$ given by:
$$
P( x_1 ) :=
CT_{ x_2 , \dots , x_ k} \left [
 { { \Phi_k ( x_1 , \dots , x_k ) } \over
   {  x_2^{a_2 -1} \dots x_k^{a_k -1}  \prod_{i=2}^{k}
(1-x_1 x_i )(1- \bar x_1 x_i )}
} 
  \right ] \quad .
\eqno(H2c)
$$
 
Since $P( x_1)$ is obviously anti-symmetric w.r.t. 
to $x_1 \rightarrow 1 - x_1$, it follows
by Crucial Fact $\aleph_2$, that we would be done once we can prove that
the degree of $P( x_1)$ is $\leq 2n$.
 
Let's perform a partial-fraction decomposition of the constant-termand
of $(H2c)$, w.r.t. to $x_1$. Since the degree of $\Phi_k ( x_1 , \dots , x_k)$
in $x_1$ is $\leq 4k-3$, we can set:
$$
 { { \Phi_k ( x_1 , \dots , x_k ) } \over
   {  \prod_{i=2}^{k}
(1-x_1 x_i ) (1- \bar x_1 x_i ) }
} 
 =
Q_k ( x_1 ; x_2 ,\dots , x_k ) + 
\sum_{i=2}^{k} {{B_i} \over {1- x_1 x_i } } +
\sum_{i=2}^{k} {{C_i} \over {1- \bar x_1 x_i } }  \quad ,
\eqno(H2d)
$$
where $B_i$ and $C_i$ do not depend on $x_1$, 
and $Q_k$ is a polynomial of degree $\leq (4k-3)-(2k-2) = 2k-1$,
in $x_1$, with coefficients that are rational functions of
$(x_2 , \dots , x_k )$. Inserting  this into $P(x_1)$ above shows
that the contribution to it, due to $Q_k$ has degree $\leq 2k-1 < 2n$.
We now have to determine the $B_i$'s. By the anti-symmetry of 
$\Phi_k$, we have that $C_i = -B_i$. $B_i$ is given by the following
$(sub)^6$lemma:
 
{\bf Subsubsubsubsubsublemma 1.2.1.2.3.1.1:} $B_i$ in $(H2d)$ is given by:
$$
\displaylines{
B_i \,\, = \,\,
\pm  {{(1- 2 x_i ) (1- x_i^2) (1- x_i \bar x_i ) 
 \Phi_{k-2}( x_2 , x_3 , \dots , \hat x_i , \dots , x_k ) }
\over { x_i^{k+1} } }\cdot
\cr
\prod_{ {{j=2} \atop { j \neq i }}}^k
\{ (1- x_i x_j ) ( 1 - \bar x_i x_j )( 1- x_i \bar x_j )
(x_j + x_i -1 ) \} \quad. \cr}
$$
 
{\bf Proof:} This is a special case, $R=1$, and $z_i=x_i$, for all
$2 \leq i \leq k$, of subsubsublemma $1.4.1.1$ that is proved later, in act IV.
This completes the proof of  $(sub)^6$ lemma $1.2.1.2.3.1.1$. \halmos
 
Since $C_i = -B_i$, one can write:
$$
P(x_1)
 =
CT_{x_2 , \dots , x_k } {{Q_k ( x_1 ; \dots , x_k )} \over
{x_2^{a_2-1} \dots x_k^{a_k-1}}} \,\, +  
CT_{x_2 , \dots x_k }
\sum_{i=2}^{k} B_i ( {{1} \over { x_2^{a_2-1} \dots x_k^{a_k -1}
(1- x_1 x_i) } } -
{{1} \over {x_2^{a_2-1} \dots x_k^{a_k -1} (1- \bar x_1 x_i )} } ) \quad .
\eqno(H2e)
$$
 
We have already observed that the degree (in $x_1$) of the first piece
is $\leq 2n$.  Since $B_i/(1- x_1 x_i )$ 
possesses a Laurent series, we can use crucial fact $\aleph_4$
to change the order of the constant-terming, for
the $i^{th}$ term, and get that its contribution is
$$
CT_{x_2 , \dots , \hat x_i , \dots , x_k }
 \{ {
{
\pm  \Phi_{k-2}( x_2 , x_3 , \dots , \hat x_i , \dots , x_k ) 
}
 \over { x_2^{a_2-1} \dots x_{i-1}^{a_{i-1}-1}
x_{i+1}^{a_{i+1}-1} \dots x_k^{a_k -1} } }
\cdot
$$
$$
CT_{x_i} \{ 
{ 
{(1- 2 x_i ) (1- x_i^2) (1- x_i \bar x_i ) 
\prod_{ {{j=2} \atop { j \neq i }}}^k
[ (1- x_i x_j ) ( 1 - \bar x_i x_j )( 1- x_i \bar x_j )
(x_j + x_i -1 )] 
}
\over 
{ x_i^{a_i + k} }
} 
 \cdot {{1} \over {(1- x_i x_1)}} \}\} \quad .
$$
 
Using the Maclaurin expansion  (in $x_i$)
$$
(1- x_1 x_i )^{-1} = \sum_{r=0}^{\infty} x_1^r x_i^r \quad ,
$$
and expanding the numerator of the last constant-termand 
(w.r.t.  $x_i$) in powers of $x_i$:
$$
\sum_{r=0}^{4k-3} L_r (x_2 , \dots , \hat x_i , \dots , x_k ) x_i^r \quad ,
$$
multiplying out, and then extracting the coefficient of
$x_i^{a_i +k}$, shows that the inner constant-term (w.r.t. $x_i$,)
can be written as:
$$
\sum_{r=0}^{\min(a_i +k ,4k-3)} 
L_r (x_2 , \dots , \hat x_i , \dots , x_k ) x_1^{a_i + k-r}
 \quad ,
$$
which is a polynomial of degree
$\leq k+ a_i \leq n+n=2n$ in $x_1$, with coefficients that are
polynomials in the remaining variables 
$(x_2 , \dots , \hat x_i , \dots , x_k )$. Taking the outer constant-term
does not change this fact, except that now we arrive at a polynomial
in $x_1$ with {\it numerical} coefficients. The same argument
(with $\bar x_1$ replaced by $x_1$) holds for the term
$-B_i/(1-\bar x_1 x_i )$.
 
Since each and every contribution to $P(x_1)$ has degree $\leq 2n$,
the same must be true of $P(x_1)$ itself, and hence completes the
proof of subsubsubsubsublemma $1.2.1.2.3.1$ . \halmos
 
This completes the proof of subsubsubsublemma 1.2.1.2.3 \halmos.   
 
{\bf Subsubsubsublemma 1.2.1.2.4:} For
$(k \, ; \, k ,  a_2 , \dots , a_k ) \in \overline{Land\_Of\_Gog_k}$,
$$
H_k (k \, ; \, k ,  a_2 , \dots , a_k ) \,=\,
H_{k-1}(k; a_2 , \dots , a_k ) \quad .
\eqno(H3)
$$
{\eightrm  [ Type `{\eighttt S12124(k):}' in ROBBINS, for specific k.]}
 
{\bf Proof:} From the definition of $H_k$,
$$
H_k ( k ; k , a_2 , \dots , a_k )=
CT_{x_1 , \dots , x_k} 
\left \{
{
{\Phi_k( x_1 , \dots , x_k )} \over
{x_1^{k-1} x_2^{a_2-1} \dots x_k^{a_k-1} \bar x_1^{k+1} \dots
\bar x_k^{2k} }
}
 \prod_{ 1 \leq i < j \leq k}
(1- x_i x_j )^{-1} (1- \bar x_i x_j )^{-1}
\right \} \quad .
\eqno(1.2.1.2.4a)
$$
 
Recall, from the definition of $\overline{Land\_Of\_Gog_k}$
(given at the beginning of the proof of $1.2.1$), that $k > a_2$.
Let's look at the contributions of the various terms of $\Phi_k$
(defined in ($Gog_1$)), to the constant term (1.2.1.2.4a).
Because of the low exponents ($<k$) of the $x_i$'s that appear
in the denominator of (1.2.1.2.4a),
it is clear that only those $g=( \pi , \epsilon )$ for which
all the $\epsilon_i$'s equal $-1$, have any hope of
contributing a non-zero term. It follows that (note that
we ``lose'' a factor of $(-1)^k$ due to fixing
all the $k$ $\epsilon_i$'s to be $-1$):
$$
H_k ( k ; k , a_2 , \dots , a_k )=
CT_{x_1 , \dots , x_k} 
\left \{
{
{ \Psi_k ( x_1 , \dots , x_k )} \over
{x_1^{k-1} x_2^{a_2-1} \dots x_k^{a_k-1} \bar x_1^{k+1} \dots
\bar x_k^{2k} }
}
 \prod_{ 1 \leq i < j \leq k}
(1- \bar x_i x_j )^{-1} (1-  x_i x_j )^{-1}
\right \} \quad ,
\eqno(1.2.1.2.4b)
$$
where
$$
\Psi_k ( x_1 , \dots , x_k ) = 
\,
\sum_{ \pi \in \S_k  } \sgn(\pi) \, \cdot \,
\pi \left [ \prod_{i=1}^{k}  x_i^{k-i} \bar x_i^{k} 
\prod_{1 \leq i < j \leq k} (1- \bar x_i x_j )(1-  x_i x_j ) \right ]
\qquad.
\eqno(Gog_1')
$$
 
Even now, there are lots of terms of $\Psi_k$ that do not contribute
to the constant term of (1.2.1.2.4b). All those terms corresponding to
$\pi$ for which $\pi(1) \neq 1$, will have one of the $x_j$, for
$j>1$ raised to a power $k-1$, which will introduce a factor of
$x_j$ and shatter any hope for a non-zero 
contribution to the constant term. So let's write:
$$
\displaylines{
\Psi_k ( x_1 , \dots , x_k ) =
 x_1^{k-1} \bar x_1^k  (1- \bar x_1 x_2 ) \dots (1- \bar x_1 x_k)
(1-  x_1 x_2 ) \dots (1-  x_1 x_k) \Psi_{k-1} (x_2 , \dots , x_k )
\, +  \, GARBAGE
\cr
= x_1^{k-1} \bar x_1^k ( \bar x_2 \dots \bar x_k ) (1) \dots (1)
\Psi_{k-1} (x_2 , \dots , x_k ) + O(x_1^k) + \, GARBAGE \quad ,
\cr}
$$
where $GARBAGE$ denotes terms that clearly contribute $0$ to
the constant term, and in the second equality we have used that
$(1- \bar x_1 x_2 ) \dots (1- \bar x_1 x_k )=$
$( \bar x_2 + x_1 x_2 ) \dots ( \bar x_k + x_1 x_k ) =$
$\bar x_2 \dots \bar x_k +O( x_1)$.
 
Thus,
$$
\displaylines{
H_k ( k \, ; \, k , a_2 , \dots , a_k )=
 CT_{x_2 , \dots , x_k} CT_{x_1}
 \{
{
{ x_1^{k-1} \bar x_1^k \bar x_2 \dots \bar x_k
\Psi_{k-1} ( x_2 , \dots , x_k )} \over
{x_1^{k-1} x_2^{a_2-1} \dots x_k^{a_k-1} \bar x_1^{k+1} \dots
\bar x_k^{2k} }
} \cdot
\cr
\prod_{1 < j \leq k} (1- x_1 x_j )^{-1} (1- \bar x_1 x_j )^{-1}
 \prod_{ 2 \leq i < j \leq k}
(1- x_i x_j )^{-1} (1- \bar x_i x_j )^{-1}
 \} 
\cr
=
 CT_{x_2 , \dots , x_k}
\{
{
   { \bar x_2 \dots \bar x_k \Psi_{k-1} ( x_2 , \dots , x_k )} 
  \over
{ x_2^{a_2-1} \dots x_k^{a_k-1} \bar x_2^{k+2} \dots
\bar x_k^{2k} }
} 
\cdot
 \prod_{ 2 \leq i < j \leq k}
(1- x_i x_j )^{-1} (1- \bar x_i x_j )^{-1} \cdot
\cr
CT_{x_1} \{ (1/ \bar x_1 )  
\prod_{1 < j \leq k} (1- x_1 x_j )^{-1} (1- \bar x_1 x_j )^{-1} \} \}
\cr
=  CT_{x_2 , \dots , x_k}
\{
{
{
\Psi_{k-1} ( x_2 , \dots , x_k )} \over
{ x_2^{a_2-1} \dots x_k^{a_k-1} \bar x_2^{k+1} \dots
\bar x_k^{2k-1} }
} \cdot
 \prod_{ 2 \leq i < j \leq k}
(1- x_i x_j )^{-1} (1- \bar x_i x_j )^{-1}
 \, \cdot 1 \ \} \quad .
\cr}
$$
$$
\eqno(H_k (k))
$$
 
I claim that this equals $H_{k-1} (k; a_2 , \dots , a_k)$.
By the definition $(H_k)$ given in the statement of
subsublemma $1.2.1$, we have:
$$
H_{k-1}(k; a_2 , \dots , a_k ) :=
CT_{x_2 , \dots , x_k} \left \{ 
{ { \Phi_{k-1} ( x_2 , ... , x_k ) } 
\over { \left \{ \prod_{i=2}^k x_i^{a_i-1} (\bar x_i)^{(k-1)+i} \right \} 
 \prod_{ 2 \leq i < j \leq k}
(1-x_i x_j ) (1- \bar x_i x_j )
} }  \right \} \quad ,
\eqno(H_{k-1}(k))
$$
 
Since $k> a_2 \geq a_3 \dots \geq a_k$, a parallel
argument can be given to show that the right side of $(H_{k-1}(k))$ equals
the right side of $(H_{k}(k))$. This completes the proof of
subsubsubsublemma 1.2.1.2.4. \halmos 
 
{\bf Subsubsubsublemma 1.2.1.2.5:} $(H4)$ in the statement
of subsubsublemma $1.2.1.2$ holds. \halmos 
 
This completes the proof of subsubsublemma $1.2.1.2$. \halmos
 
{\bf Subsubsublemma 1.2.1.3:} Let $k \geq 1$ and for any vector
of integers $a=(a_1 ,\dots , a_k)$, define the partial difference
operator  $P^{(a)}$ by:
$$
P^{(a)} ( A_1 , \dots , A_k ) :=
\prod_{
\{i| a_i - a_{i+1} > 0 \}
}
(I- A_i^{-1}) \qquad,
$$
where we declare that $a_{k+1}=0$. There is a unique
sequence of discrete functions $Y_k(n; a_1 , \dots , a_k )$,
defined on $\overline{Land\_Of\_Gog_k}$ that 
satisfies the following partial recurrence 
$$
P^{(a)} (A_1 , \dots , A_k ) Y_k ( n \, ; \, a_1 , \dots , a_k )
$$
$$
=
Y_k (n-1 \, ; \, a_1 \,\, ,\,\, \min ( a_1 -1 , a_2 ) \,\, , \,\,
\dots \,\, , \,\, \min (a_{k-1} -1,a_k ) ) \qquad ,
\eqno(P \Delta E(Y_k))
$$
valid for $n>k$ and for $(n; a_1 , \dots , a_k) \in Land\_Of\_Gog_k$,
and the following initial/boundary conditions for 
$(n \, ; \, a_1 \, , \, \dots \, , \, a_k \,) \in \overline{Land\_Of\_Gog_k}$:
$$
Y_k (n \, ; \, a_1 , \dots , a_k ) = 0  \quad \hbox{on} \quad
a_i = k-i \qquad  ,
\eqno(Y1_i)
$$
$$
Y_k (n \, ; \, n+1 , a_2 , \dots , a_k ) = 0 \quad ,
\eqno(Y2) 
$$
$$
Y_k ( k \, ; \, k , a_2 , \dots , a_k ) =
Y_{k-1} ( k \, ; \, a_2 , \dots , a_k ) \qquad ,
\eqno(Y3) 
$$
and, finally the ``Adam'' condition
$$
Y_1 ( 1 \, ; \, 1 ) = 1  \qquad .
\eqno(Y4)
$$
{\eightrm  [ Type `{\eighttt S1213(k,n):}' in ROBBINS, for specific k and n.]}
 
{\bf Proof:} Eq. $P\Delta E(Y_k)$ can be rewritten such that
$Y_k(n;a_1 , \dots ,a_k)$ remains on the left side, and all the rest
goes to the right side. The right side involves either points
$(n; b_1 , \dots , b_k)$ that belong to $\overline{Land\_Of\_Gog}$,
with $b_1+ \dots + b_k< a_1 + a_2 + \dots + a_k$, so they were already
computed before, or points that belong to the boundary, for which the
boundary conditions are used. Note that
the reason that it was possible
to impose the extra restriction that $n=k=a_1 \Rightarrow a_2<k$,
in the definition of $\overline {Land\_Of\_Gog_k}$ ,
was that when we use the local recurrence ($P \Delta E(\tilde M_k)$) 
to compute $\tilde M_k$,
starting with $n=k+1$, the right hand side of ($P \Delta E(\tilde M_k)$)
only requires knowledge
of $\tilde M_k$, when $n=k=a_1$, for those $a$'s for which $a_2<k$, since
even if $a_2=k$ on the left side, the evaluation at the right is
at $\tilde M_k (k;k, \min(k-1,k), \dots )=\tilde M_k (k;k,k-1, \dots , )$.
This completes the proof of subsubsublemma $1.2.1.3$.
\halmos.
 
Combining subsubsublemmas $1.2.1.1-3$ yields subsublemma $1.2.1$. \halmos.
 
Going back to the proof of sublemma $1.2$, note that
$m_k (n)$, the total number of $n \times k$-Gog trapezoids, equals
$ \tilde M_k ( n+1 ; n+1 \,,\, n+1 \,,\, \dots \,,\, n+1 )$, as
follows directly from the definition $(\tilde M_k)$ of $\tilde M_k$.
By subsublemma $1.2.1$, this equals 
$  H_k ( n+1 ; n+1 \,,\, n+1 \,,\, \dots \,,\, n+1 )$ which coincides with
the right side of
$(GogTotal)$. This completes the proof of
sublemma $1.2$. \halmos.
\bigskip
\centerline{\bf INTERACT}  
 
We have now established constant term expressions both
for $b_k(n)$, the number of $n \times k$ Magog trapezoids,
given by $(MagogTotal)$ in sublemma  $1.1$ , and for 
$m_k(n)$, the number of $n \times k$ Gog trapezoids,
given by $(GogTotal)$ in sublemma  $1.2$. It would have been
nice if the {\it constant-termands} were the same.
Unfortunately, they are not. In order to prove that
the {\it constant terms} are,  we have to resort
to the omnipotent Stanton-Stembridge trick once  again,
only now w.r.t. to the action of the group of
signed permutations, $W(\B_k)$. Justifying its use, however,
will require many more pages.
It would be necessary first  to convert the constant-term
expressions to ``residue'' expressions. 
 
{\bf Sublemma 1.1':} The total number of $n \times k -$ Magog trapezoids,
$b_k (n)$, is given by the following iterated-residue expression:
$$
b_k (n) = Res_{x_1, \dots , x_k}  \,
\left \{ 
{
{\Delta_k (x_1 , \dots , x_k )} \over
{\prod_{i=1}^k x_i^{n+k-i} ( \bar x_i)^{n+k+1} \,  
 \prod_{1 \leq i < j \leq k} (1- x_i x_j ) }
} \right \} ,
\eqno(MagogTotal')
$$
where,
$$
\Delta_k (x_1 , ... , x_k ) :=
\prod_{i=1}^{k}( 1- 2 x_i )
\prod_{1 \leq i < j \leq k}
(x_j - x_i )(x_j+x_i-1) \qquad .
\eqno(Delta)
$$
 
{\bf Sublemma 1.2':} The total number of $n \times k -$ Gog trapezoids, 
$m_k (n)$, is given by the following iterated-residue expression:
$$
m_k(n) \, = \,
Res_{x_1 , \dots , x_k} \left [
 { { \Phi_k ( x_1 , \dots , x_k ) } \over
   { \left \{ \prod_{i=1}^k x_i^{n+1} ( \bar x_i)^{n+i+1} \right \} 
 \prod_{ 1 \leq i < j \leq k}
(1-x_i x_j ) (1- \bar x_i x_j )
}
}
 \right  ] \qquad ,
\eqno(GogTotal')
$$
where the polynomial $\Phi_k(x_1 , \dots , x_k)$ is defined by:
$$
\Phi_k ( x_1 , \dots , x_k ) = 
(-1)^k \,
\sum_{ g \in W( \B_k ) } \sgn(g) \, \cdot \,
g \left [ \prod_{i=1}^{k} \bar x_i^{k-i} x_i^{k} 
\prod_{1 \leq i < j \leq k} (1- x_i \bar x_j )(1- \bar x_i \bar x_j ) \right ]
\qquad.
\eqno(Gog_1)
$$
 
{\bf Proof of Sublemmas 1.1' and 1.2':}
Use Crucial fact $\aleph_6$. \halmos
 
\vfill
\eject
 %%%%Proof of the Alternating Sign Matrix Conj., Act III
 
%begin macros
\baselineskip=14pt
\parskip=10pt
\def\B{{\cal B}}
\def\S{{\cal S}}
\def \inv{\mathop{\rm inv} \nolimits}
\def \sgn{\mathop{\rm sgn} \nolimits} 
\def\hat{\widehat}
\def\tilde{\widetilde}
\def\epsilon{\varepsilon}
\def\halmos{\hbox{\vrule height0.15cm width0.01cm\vbox{\hrule height
 0.01cm width0.2cm \vskip0.15cm \hrule height 0.01cm width0.2cm}\vrule
 height0.15cm width 0.01cm}}
\font\eightrm=cmr8  
\font\eighttt=cmtt8
\magnification=\magstephalf

\parindent=0pt
\overfullrule=0in
\headline={\rm  \ifodd\pageno  \RightHead  \else  \LeftHead  \fi}
\def\RightHead{\centerline{Proof  of the ASM Conjecture-Act III}}
\def\LeftHead{ \centerline{Doron Zeilberger}}
%\pageno=38
%end macros
 
\centerline{\bf Act III. (MagogTotal') IS UNFAZED BY THE SIGNED PERMUTATIONS'
ACTION} 
 
\qquad\qquad\qquad {\it  
Talk no more so exceeding proudly; let not arrogancy come out of
your mouth: for the Lord is a God of knowledge, and by him actions are
weighed. }\qquad\qquad\qquad ---(Samuel II,3)
\bigskip 
{\bf Sublemma 1.3 :} Let $n \geq k \geq 1$ and let
$f_{n,k} (x)= f_{n,k}( x_1 , \dots , x_k )$
be the residuand of $(MagogTotal')$:
$$
f_{n,k} ( x_1 , \dots , x_k ) :=
 { {\Delta_k ( x_1 , \dots , x_k ) } \over
   { \prod_{i=1}^{k} ( \bar x_i )^{n+k+1} x_i^{n+k-i}
   \prod_{ 1 \leq i < j \leq k } ( 1 - x_i x_j ) }
} \quad ,
$$
where $\Delta_k ( x_1 , \dots , x_k)$ is the polynomial defined in ($Delta$)
(in the statement of sublemma $1.2$).
Let $g=(\pi,\epsilon)$ be any signed permutation in $W(\B_k)$; then
$$
Res_{x_1 , \dots , x_k}[g f_{n,k}(x) ] = Res_{x_1 , \dots , x_k}
 [  f_{n,k} (x ) ]  \quad (=b_k (n)) \quad .
$$
{\eightrm  [ Type `{\eighttt S13(k,n):}' in ROBBINS, for specific k and n.]}
 
{\bf Proof}: Recall that a signed permutation $g=(\pi,\epsilon)$ acts
on any rational function $f(x)$ of $x=(x_1, \dots , x_k)$, by
$g f(x)= f(\epsilon(\pi(x)))$.
By renaming the variables, and replacing $\pi^{-1}$ by $\pi$, one finds that
the statement of sublemma 1.3 is equivalent to:
 
{\bf Sublemma 1.3' :} Let $f_{n,k} (x)= f_{n,k}( x_1 , \dots , x_k )$
be the residuand of $(MagogTotal')$, reproduced just above.
Let $\pi$ be any permutation of
$\{ 1 , \dots , k \}$, and let 
$\epsilon = (\epsilon_1 , \dots , \epsilon_k )$ be any sign-assignment,
then
$$
Res_{x_{\pi(1)} , \dots , x_{\pi(k)}}
[ f_{n,k} ( \epsilon_1(x_1) , \dots , \epsilon_k(x_k) ) ] = 
Res_{x_1 , \dots , x_k} [  f_{n,k} (x_1 , \dots , x_k ) ]  \quad  .
\eqno(1.3')
$$
 
The proof will be by induction on the number of $\epsilon_i$'s
that are equal to $-1$. If there are no $-1$'s in $\epsilon$, then
$1.3'$ follows immediately from crucial fact $\aleph_4'$,
since then the residuand $f_{n,k}( x_1 , \dots , x_k )$ has a Laurent
series. Otherwise, let $u$ be the smallest integer $i$, ($1 \leq i \leq k$)
such that $\epsilon_{\pi(i)} = -1$.  
The inductive step is provided by
the following subsublemma, that asserts that it is permissible to change
$\epsilon_{\pi(u)}$ from $-1$ to $+1$ 
inside the residuand of the iterated residue on the left side 
of $(1.3')$, without affecting its value.
 
{ \bf Subsublemma 1.3.1:} Let $f_{n,k} ( x_1 , \dots , x_k )$ (where, 
as always, $n \geq k \geq 1$,) be
the residuand of $(MagogTotal')$, let $\pi$ be any permutation of
$\{ 1 , \dots , k \}$, and let 
$\epsilon = ( \epsilon_1 , \dots , \epsilon_k )$ be any sign-assignment
with at least one $\epsilon_i$ equal to $-1$.
Let $u$ be the smallest $i$ ($1 \leq i \leq k$) such that
$\epsilon_{\pi(i)} = -1$. (i.e., we have
$\epsilon_{\pi(1)}= \dots  = \epsilon_{\pi(u-1)}= +1$,
and $\epsilon_{\pi(u)}= -1$.) Then:
$$
Res_{x_{\pi(1)} , \dots , x_{\pi(k)}}
[ f_{n,k} ( \epsilon_1(x_1) , \dots ,  x_{\pi(u)} , \dots ,
  \epsilon_k(x_k) ) ] = 
Res_{x_{\pi(1)} , \dots , x_{\pi(k)}}
[ f_{n,k} ( \epsilon_1(x_1) , \dots ,  \bar x_{\pi(u)} , 
\dots , \epsilon_k(x_k) ) ] 
\quad  .
$$
{\eightrm  [ Type `{\eighttt S131(k,n):}' in ROBBINS, for specific k and n.]}
 
{\bf Proof:} Let's abbreviate $R:= \pi (u)$, and
$z_i := \epsilon_i (x_i)$ ( $1 \leq i \leq k , i \neq R $.)
Recall that this means that $z_i=x_i$ if $\epsilon_i=+1$,
and $z_i = \bar x_i = 1 - x_i$ if $\epsilon_i = -1$.
By the $W(\B_k)$-antisymmetry of $\Delta_k(x)$, we have
$$
f_{n,k} ( \epsilon_1(x_1) , \dots ,  x_R , \dots , \epsilon_k(x_k) )
=
{ {1} \over { x_R^{n+k-R} \bar x_R^{n+k+1} } }
\prod_{ {{j=1} \atop {j \neq R } }}^{k}
{{1} \over {z_j^{n+k-j} \bar z_j^{n+k+1} }}
\prod_{ {{1 \leq r < s \leq k } \atop {r,s \neq R} }}
{{1} \over {(1- z_r z_s ) }} \cdot
\left \{
{ {\pm \Delta_k ( x_1 , \dots , x_k ) } \over
   { \prod_{{{i=1} \atop { i \neq R}} }^{k} (1- z_i x_R ) }} \right \} \quad .
\eqno(Jane)
$$
 
The next step is to perform a partial-fraction decomposition,
with respect to the variable $x_R$, of the fraction inside
$(Jane)$'s braces.
 
{\bf Subsubsublemma 1.3.1.1:}
$$
{ {\pm \Delta_k ( x_1 , \dots , x_k ) } \over
   { \prod_{{{i=1} \atop { i \neq R}} }^{k} 
(1- z_i x_R ) 
} }
=  P_R +
\sum_{ {{i=1} \atop {i \neq R } }}^{k} 
{  {B_i} \over {1- z_i x_R} } 
\quad ,
\eqno(Tamar)
$$
where $P_R$ is a polynomial of degree $k$ in $x_R$,
with coefficients that are rational functions of the other
variables $(x_1 , \dots , x_{R-1} , x_{R+1} , \dots , x_k )$, and
$$
B_i =  {{\pm ( z_i -2 ) (1- 2 z_i )(1- z_i^2 ) (1- z_i \bar z_i ) 
\Delta_{k-2} (x_1 , \dots ,
\hat {x_i} , \dots , \hat x_R , \dots , x_k ) }
\over {z_i^{k+1}} } \cdot
$$
$$
\left \{ 
\prod_{ {{j=1} \atop {j \neq i,R}} }^{k}
(1- z_i z_j ) (1- z_i \bar z_j ) (z_j + z_i -1)
\right \}
\quad .
\eqno(Bi)
$$
{\eightrm  [ Type `{\eighttt S1311(k,R):}' in ROBBINS, for specific k and R.]}
 
{\bf Proof:} By Freshman Calculus (or rather the algebraic part of it),
we know that the left side of $(Tamar)$, when viewed as a rational function
of $x_R$, can be written as on its right side, for
{\it some } polynomial $P_R$, and {\it constants} (in the present
context, i.e. quantities  not involving $x_R$) $B_i$.
The degree of $\Delta_k$, in $x_R$
(or for that matter, in any variable),
from its definition in $(Delta)$, is
$ 2k-1$, and the degree of the denominator of the left side of
$(Tamar)$ is $(k-1)$. Hence the degree of $P_R$, in $x_R$ is
$(2k-1)-(k-1)=k$. 
The coefficients of $P_R$, as a polynomial of $x_R$, are certain
rational functions of the remaining variables whose exact form
does not interest us. The exact form of the $B_i$'s
{\it should} interest us,  since it is used heavily in the proof.
So let's brace ourselves and embark on partial-fractioning.
 
Multiplying both sides of $(Tamar)$ by the denominator of
its left side yields the equation:
$$
\pm \Delta_k ( x_1 , \dots , x_k ) 
=
 \left \{   \prod_{{{j=1} \atop {j \neq R} }}^{k} 
(1- z_j x_R ) \right \} P_R 
+ \sum_{ {{i=1} \atop {i \neq R } }}^{k} 
 \left \{   \prod_{{{j=1} \atop {j \neq R,i} }}^{k} 
(1- z_j x_R ) \right \} B_i \quad .
\eqno(Tamar')
$$
 
To decipher $B_i$ ($1 \leq i \leq k$, $i \neq R$), we plug in
$x_R=1/z_i$ on both sides of $(Tamar')$ to get:
$$
B_i=
{  {\pm \Delta_k ( x_1 , \dots , x_k ) \vert_{x_R=1/z_i} } \over
{ \left \{   \prod_{{{j=1} \atop {j \neq R,i} }}^{k} 
(1- z_j/z_i ) \right \} }}
= {{ \pm  z_i^{k-2} \Delta_k ( x_1 , \dots , x_k ) \vert_{x_R=1/z_i}} \over
{ \left \{   \prod_{{{j=1} \atop {j \neq R,i} }}^{k} 
( z_i - z_j ) \right \} }} \quad .
\eqno(PP1)
$$
 
{\bf Subsubsubsublemma 1.3.1.1.1:}
$$
\Delta_k ( x_1 , x_2 , \dots , x_k ) \vert_{x_R=1/z_i} =
\pm {{(z_i -2 ) (1- 2 z_i ) (1- z_i^2) (1- z_i \bar z_i ) 
\cdot \Delta_{k-2}
( x_1 , x_2 , x_3 , \dots , \hat x_i , \dots , \hat x_R , \dots ,
x_k ) }
\over { z_i^{2k-1} }} \cdot
$$
$$
\prod_{ {{j=1} \atop { j \neq i,R }}}^{k}
\{ (1- z_i z_j ) ( 1 -  z_i \bar z_j )(z_j - z_i ) (z_j + z_i -1) \} \quad .
\eqno(1.3.1.1.1)
$$
{\eightrm  [ Type `{\eighttt S13111(k,R,i):}' in ROBBINS, for specific k,R 
and i.]}
 
{\bf Proof:} This follows routinely from the definition $(Delta)$ of
$\Delta_k$. \halmos .
 
Plugging $(1.3.1.1.1)$ into $(PP1)$, and cancelling
out $\prod_{ {{j=1 } \atop {j \neq i,R} }}^{k} ( z_j - z_i ) $
from top and bottom,
finishes up the proof of $(sub)^3$lemma $1.3.1.1$. \halmos
 
For the 
rest of the proof of $(sub)^2$lemma $1.3.1$, we need some notation.
Let $POL$ stand for ``some polynomial of $(x_1 , \dots , x_k )$'',
and let $RAT$ stand for ``some rational function of $(x_1 , \dots , x_k )$''.
$POL( \hat x_R )$ stands for such a polynomial that does not involve $x_R$.
$POLRAT( x_R; x_1 , \dots , \hat x_R , \dots , x_k)$, or $POLRAT(x_R \,\, ; \, \, \dots )$
for short, stands for a polynomial in $x_R$ with coefficients that
are rational functions of the rest of the variables.
$POLRAT(x_i ; \hat x_R)$ stands for a polynomial in $x_i$
with coefficients that are rational function of the remaining
variables {\it and} that does not involve $x_R$.

Inserting $(Tamar)$(in the statement of $1.3.1.1$) 
into $(Jane)$ is going to yield the following $(sub)^3$lemma:
 
{\bf Subsubsublemma 1.3.1.2:}
$$
f_{n,k} ( \epsilon_1(x_1) , \dots ,  x_R , \dots , \epsilon_k(x_k) )
=  {{\tilde P_R} \over {x_R^{n+k-R} \bar x_R^{n+k+1} } } +
\sum_{ {{i=1} \atop {i \neq R } }}^{k} 
{  {\tilde B_i} \over {1- z_i x_R} }  \quad ,
\eqno(Celia)
$$
where $\tilde P_R$ is a polynomial of degree $k$ in $x_R$,
with coefficients that are rational functions of the rest of the
variables, and $\tilde B_i$ can be written as follows.
$$
\tilde B_i = {
{POL(\hat x_R) } \over { x_R^{n+k-R} \bar x_R^{n+k+1}
                        z_i^{n+2k-i+1} \bar z_i^{n+k+1} }} 
\cdot \prod_{ {{j=1} \atop {j \neq i,R} }}
^k {{1} \over {z_j^{n+k-j} \bar z_j^{n+k+1} }} \cdot
\prod_{ {{1 \leq r < s \leq k } \atop {r,s \neq R,i} }}
{{1} \over {(1- z_r z_s ) }}
$$
$$
={{POLRAT(x_i ;\hat x_R) } \over
   {x_R^{n+k-R} \bar x_R^{n+k+1} z_i^{n+2k-i+1} \bar z_i^{n+k+1} } }
\quad .
\eqno(\tilde Bi)
$$
{\eightrm  [ Type `{\eighttt S1312(k,n,eps,R):}' in ROBBINS, for specific k,n,
eps and R.]}
 
{\bf Proof:} Insert $(Tamar)$ into $(Jane)$, bringing the quantity
in front of the sum inside, and for each $B_i$ split
the double product:
$$
\prod_{ {{1 \leq r < s \leq k } \atop {r,s \neq R} }}
{{1} \over {(1- z_r z_s ) }}=
\prod_{ {{1 \leq r < s \leq k } \atop {r,s \neq R,i} }}
{{1} \over {(1- z_r z_s ) }} \cdot
\prod_{ {{j=1 } \atop {j \neq i,R} }}^{k}
{{1} \over {(1- z_j z_i ) }} \quad ,
\eqno(Split)
$$
and cancel out 
$\prod_{ {{j=1 } \atop {j \neq i,R} }}^{k} (1- z_j z_i ) $
on the top and bottom.
This completes  the proof of $(sub)^3$lemma
$1.3.1.2$. \halmos
 
We will say that a rational function $h=h( x_1 , \dots , x_k )$ has
property $(Hadas)$ if:
$$
Res_{x_{\pi(1)} , \dots , x_{\pi(k)}} [h] =
Res_{x_{\pi(1)} , \dots , x_{\pi(k)}} [h( x_R \rightarrow \bar x_R )] \quad ,
\eqno(Hadas)
$$
where  $h( x_R \rightarrow \bar x_R )$ is the rational function obtained from $h$
by replacing the variable $x_R$ by $\bar x_R$ ($=1- x_R $.)
 
The statement of subsublemma $1.3.1$ is equivalent to saying that
$$
f_{n,k} ( \epsilon_1(x_1) , \dots ,  x_R , \dots ,
\epsilon_k(x_k) ) 
$$
has property $(Hadas)$. Property $(Hadas)$ is obviously additive
(since the action of taking residue, and performing $x_R \rightarrow \bar x_R$ are),
so in order to complete the proof of $1.3.1$ (and hence of $1.3$), it
would suffice to show that each term in $(Celia)$ has property
$(Hadas)$. This will be accomplished by the following two
subsubsublemmas ($1.3.1.3-4$.)
 
{\bf Subsubsublemma 1.3.1.3:}
Let $\tilde P_R$ be as in the statement of $1.3.1.2$, i.e.
a polynomial of degree $k$ in $x_R$ with coefficients
that are rational functions of the remaining variables.
Then ${{\tilde P_R} \over {x_R^{n+k-R} \bar x_R^{n+k+1} } }$ 
has property $(Hadas)$:
$$
Res_{x_{\pi(1)} , \dots , x_{\pi(k)}} 
\left [{{\tilde P_R} \over {x_R^{n+k-R} \bar x_R^{n+k+1} } } \right ] =
Res_{x_{\pi(1)} , \dots , x_{\pi(k)}} 
\left [ {{\tilde P_R} \over {x_R^{n+k-R} \bar x_R^{n+k+1} } } 
( x_R \rightarrow \bar x_R ) \right ] 
\quad .
\eqno(Hadas_1)
$$
{\eightrm  [ Type `{\eighttt S1313(k,n):}' in ROBBINS, for specific k and n.]}
 
{\bf Proof:} Recall that $R=\pi(u)$. Since $\tilde P_R$ is
a $POLRAT( x_R ; x_{\pi(1)} , \dots , \hat x_R , \dots , x_k )$, we have that
$$
Res_{x_{\pi(u+1)} , \dots , x_{\pi(k)}} 
\left [{{\tilde P_R} \over {x_R^{n+k-R} \bar x_R^{n+k+1} } } \right ]
=
Res_{x_{\pi(u+1)} , \dots , x_{\pi(k)}} 
\left [ {{POLRAT(x_R; x_{\pi(1)} , \dots , \hat x_R , x_{\pi(u+1)},
\dots , x_{\pi(k)})} \over {x_R^{n+k-R} \bar x_R^{n+k+1} } }  \right ]
$$
$$
={{POLRAT(x_R; x_{\pi(1)} , \dots , x_{\pi(u-1)})} 
\over {x_R^{n+k-R} \bar x_R^{n+k+1} } } 
={{x_R^{R+1} \cdot POLRAT(x_R; x_{\pi(1)} , \dots , x_{\pi(u-1)})} 
\over {x_R^{n+k+1} \bar x_R^{n+k+1} } } \quad .
$$
 
Since  the degree of
$x_R^{R+1} POLRAT(x_R \,\, ; \, \, \dots)$, in $x_R$, is 
$\leq (R+1)+ k \leq 2(n+k)$, it follows from crucial fact $\aleph_5'$ that 
$$
Res_{x_R} \left [ {{POLRAT ( x_R \,\, ; \, \, \dots) } \over 
{x_R^{n+k-R} \bar x_R^{n+k+1} } } \right ] =
Res_{x_R} \left [ {{POLRAT ( x_R \,\, ; \, \, \dots)} \over 
{x_R^{n+k-R} \bar x_R^{n+k+1} } }
( x_R \rightarrow \bar x_R ) \right ] \quad .
$$
 
Applying $Res_{x_{\pi(1)} , \dots , x_{\pi(u-1)}}$ to both sides of
the above equation establishes that 
${{\tilde P_R} \over {x_R^{n+k-R} \bar x_R^{n+k+1} } } $ indeed
has property $(Hadas)$, so that $(Hadas_1)$ holds. This
completes the proof of subsubsublemma $1.3.1.3$. \halmos
 
{\bf Subsubsublemma 1.3.1.4:}
Let $\tilde B_i$ be as in the statement of $1.3.1.2$, i.e.
given by $(\tilde Bi)$. Let $\pi$ be a permutation, $R$ an integer
in $[1,k] \backslash \{i\}$, and $u:=\pi^{-1}(R)$.
Then:
$$
Res_{x_{\pi(1)} , \dots , x_{\pi(k)}} 
\left [ {{\tilde B_i} \over {1- z_i x_R}} \right ] =
Res_{x_{\pi(1)} , \dots , x_{\pi(k)}} 
\left [ {{\tilde B_i} \over {1- z_i x_R}} ( x_R \rightarrow \bar x_R ) \right ]
\quad.
\eqno(Hadas_2)
$$
{\eightrm  [ Type `{\eighttt S1314(k,n):}' in ROBBINS, for specific k and n.]}
 
{\bf Proof:}
 
{\bf  Case I:} $i \in \{ \pi (u+1) , \dots , \pi (k) \}$ .
 
Let $v := \pi^{-1} (i)$, so that $i= \pi(v)$, and $v > u$.
We have that:
$$
{{\tilde B_i} \over {1- z_i x_R}} =
{
 {
  POLRAT( x_i ; x_{\pi(1)} , \dots , \hat x_R ,
    \dots , x_{\pi(v-1)} , 
  x_{\pi(v+1)} , \dots , x_{\pi(k)} )
  } \over
   {x_R^{n+k-R} \bar x_R^{n+k+1} z_i^{n+2k-i+1} 
\bar z_i^{n+k+1} (1- z_i x_R) }}
\quad .
$$
Applying $Res_{x_{\pi(v+1)} , \dots , x_{\pi(k)}}$ to it, gets
rid of the dependence on $(x_{\pi(v+1)} , \dots , x_{\pi(k)})$:
$$
Res_{x_{\pi(v+1)} , \dots , x_{\pi(k)}} \left [
{{\tilde B_i} \over {1- z_i x_R}} \right ]
={{
POLRAT(x_i ;x_{\pi(1)} , \dots , \hat x_R,
 \dots , x_{\pi(v-1)}) } \over
   {x_R^{n+k-R} \bar x_R^{n+k+1} z_i^{n+2k-i+1} 
\bar z_i^{n+k+1} (1- z_i x_R) }}
  \quad .
$$
Applying $Res_{x_i}$ to the above gives:
$$
Res_{x_i , x_{\pi(v+1)} , \dots , x_{\pi(k)}}
\left [ {{\tilde B_i} \over {1- z_i x_R}} \right ]
=Res_{x_i} \left [ 
{{
POLRAT(x_i ; x_{\pi(1)} , \dots , \hat x_R,
 \dots , x_{\pi(v-1)}) } \over
   {x_R^{n+k-R} \bar x_R^{n+k+1} z_i^{n+2k-i+1} 
\bar z_i^{n+k+1} (1- z_i x_R) }}
 \right ]
$$
$$
={{1} \over {x_R^{n+k-R} \bar x_R^{n+k+1}}} \cdot
Res_{x_i} \left [ 
{{
POLRAT(x_i; x_{\pi(1)} , \dots , \hat x_R,
 \dots , x_{\pi(v-1)} )} \over
   { z_i^{n+2k-i+1} \bar z_i^{n+k+1} (1- z_i x_R) }}
 \right ] \quad .
\eqno(Doron)
$$
 
We now distinguish two subcases:
 
{\bf Subcase Ia:} $\epsilon_i = +1$, that is: $z_i = x_i$.
 
In this subcase, the right side of $(Doron)$ equals:
$$
 {{1 } \over
{x_R^{n+k-R} \bar x_R^{n+k+1}}}
\cdot
Coeff_{x_i^{n+2k-i}} \left [ {{
 POLRAT(x_i ;x_{\pi(1)} , \dots ,
 \hat x_R  , \dots , x_{\pi(v-1)} )
 } \over
   { \bar x_i^{n+k+1} (1- x_i x_R) }}
 \right ]
  \quad .
\eqno(Gil)
$$
Expanding  everything involving $x_i$ as a power series in $x_i$:
$$
POLRAT(x_i;x_{\pi(1)} , \dots , \hat x_R , \dots , x_{\pi(v-1)}) =
\sum_{t=0}^{degree}
RAT_t (x_{\pi(1)} , \dots , \hat x_R , \dots , x_{\pi(v-1)} ) x_i^t \quad ,
$$
$$
{{1} \over {\bar x_i^{n+k+1} }} = \sum_{t=0}^{\infty} NUMBER_t \cdot x_i^t
\quad ,
$$
and, most importantly:
$$
(1- x_i x_R)^{-1} = \sum_{t=0}^{\infty} x_R^t x_i^t
\quad .
$$
Multiplying out, and collecting the coefficient of $x_i^{n+2k-i}$, gives that
$Coeff_{x_i^{n+2k-i}}$ in $(Gil)$ is a certain {\it polynomial }
in $x_R$, of degree $\leq n+2k-i$, with coefficients that are rational
functions of $(x_{\pi(1)} , \dots , \hat x_R , \dots , x_{\pi(v-1)})$.
In other words it is a $POLRAT(x_R; x_{\pi(1)} , \dots , x_{\pi(v-1)})$
with $deg_{x_R} \leq n+2k-i$. Hence,
$$
Res_{x_R, x_{\pi(u+1)} , \dots ,  x_{\pi(k)}}
\left [ {{\tilde B_i} \over {1- x_i x_R}} \right ]
=
Res_{x_R,x_{\pi(u+1)} , \dots ,  x_{\pi(v-1)}}
\left [ 
{{POLRAT(x_R ; x_{\pi(1)} , \dots , \hat x_R, \dots , x_{\pi(v-1)}) 
} \over
{x_R^{n+k-R} \bar x_R^{n+k+1}}} \right ] \quad ,
\eqno(Gil')
$$
where the $POLRAT$  in $(Gil')$
has degree $\leq n+2k-i$ in $x_R$.
Applying $Res_{x_{\pi(u+1)} , \dots , x_{\pi(v-1)}}$ to
the residuand of $(Gil')$ gets rid of the dependence on
$(x_{\pi(u+1)} , \dots , x_{\pi(v-1)})$, so we have:
$$
Res_{x_R , x_{\pi(u+1)} , \dots ,  x_{\pi(k)}}
\left [ {{\tilde B_i} \over {1- x_i x_R}} \right ]
=
Res_{x_R}
\left [ {{POLRAT(x_R ; x_{\pi(1)} , \dots , x_{\pi(u-1)}) 
} \over
{x_R^{n+k-R} \bar x_R^{n+k+1}}} \right ] \quad ,
\eqno(Gil'')
$$
where $POLRAT(x_R; x_{\pi(1)}, \dots , x_{\pi(u-1)} )$ 
has degree $\leq n+2k-i$ in $x_R$.
But, the right of $(Gil'')$ can be rewritten as:
$$
Res_{x_R}
\left [ {{
POLRAT( x_R ; x_{\pi(1)} , \dots , x_{\pi(u-1)}) 
} \over
{x_R^{n+k-R} \bar x_R^{n+k+1}}} \right ] =
Res_{x_R}
\left [ {{
x_R^{R+1} \cdot POLRAT(x_R; x_{\pi(1)} , \dots , x_{\pi(u-1)}) 
} \over
{x_R^{n+k+1} \bar x_R^{n+k+1}}} \right ] \quad .
\eqno(Gil''')
$$
 
By crucial fact $\aleph_5'$, since $(R+1)+(n+2k-i) \leq 2(n+k)$
(recall that always $n \geq k \geq 1$, $R \leq k$, and $i \geq 1$), we can 
perform $x_R \rightarrow \bar x_R$ in the residuand of $(Gil''')$, so
$$
Res_{x_R , x_{\pi(u+1)}, \dots ,  x_{\pi(k)}}
\left [ {{\tilde B_i} \over {1- x_i x_R}} \right ]
=
Res_{x_R , x_{\pi(u+1)}, \dots ,  x_{\pi(k)}}
\left [ {{\tilde B_i} \over {1- x_i x_R}} 
( x_R \rightarrow \bar x_R ) \right ] \quad.
$$
Applying $Res_{x_{\pi(1)} , \dots , x_{\pi(u-1)}}$ to both sides
of the above equation finishes up subcase Ia of the proof
of $(sub)^3$lemma $1.3.1.4$.
 
{\bf Subcase Ib:} $\epsilon_i = -1$, that is: $z_i = \bar x_i$.
 
In this subcase, the right side of $(Doron)$ equals:
$$
 {{1 } \over
{x_R^{n+k-R} \bar x_R^{n+k+1}}}
\cdot
Coeff_{x_i^{n+k}} \left [ {{
 POLRAT(x_i ; x_{\pi(1)} , \dots ,
 \hat x_R  , \dots , x_{\pi(v-1)})
 } \over
   { \bar x_i^{n+2k-i+1} (1- \bar x_i x_R) }}
 \right ]
  \quad .
\eqno(Gil)
$$
Expanding  everything involving $x_i$ as a power series in $x_i$:
$$
POLRAT(x_i ; x_{\pi(1)} , \dots , \hat x_R , \dots , x_{\pi(v-1)}) =
\sum_{t=0}^{degree}
RAT_t (x_{\pi(1)} , \dots , \hat x_R , \dots , x_{\pi(v-1)} ) x_i^t \quad ,
$$
$$
{{1} \over {\bar x_i^{n+2k-i+1} }} = \sum_{t=0}^{\infty} NUMBER_t \cdot x_i^t
\quad ,
$$
and, most importantly:
$$
(1- \bar x_i x_R)^{-1} = (1- x_R + x_i x_R )^{-1} =
( \bar x_R + x_i x_R )^{-1} =
\sum_{t=0}^{\infty} (-1)^t {{x_R^t} \over {\bar x_R^{t+1}}} x_i^t \quad ,
$$
multiplying out, and collecting the coefficient of $x_i^{n+k}$, gives that
$Coeff_{x_i^{n+k}}$ in $(Gil)$ is a certain sum of terms of the form
$$
RAT_t (x_{\pi(1)} , \dots , \hat x_R, \dots , x_{\pi(v-1)}) 
{{x_R^t} \over {\bar x_R^{t+1}}} \quad , ( t \geq 0 ) \quad .
$$
Hence,
$$
Res_{x_R, x_{\pi(u+1)} , \dots ,  x_{\pi(k)}}
\left [ {{\tilde B_i} \over {1- \bar x_i x_R}} \right ]
$$
is a sum of terms of the forms
$$
Res_{x_R,x_{\pi(u+1)} , \dots ,  x_{\pi(v-1)}}
\left [  {{ x_R^t \cdot 
RAT_t (x_{\pi(1)} , \dots , \hat x_R, \dots , x_{\pi(v-1)}) 
} \over
{x_R^{n+k-R} \bar x_R^{n+k+t+2}}} \right ] \quad ,
(t \geq 0 ) \quad .
\eqno(Gil')
$$
Applying $Res_{x_{\pi(u+1)} , \dots , x_{\pi(v-1)}}$ to
the residuand of $(Gil')$ gets rid of the dependence on
$(x_{\pi(u+1)} , \dots , x_{\pi(v-1)})$, so we have that
$$
Res_{x_R , x_{\pi(u+1)} , \dots ,  x_{\pi(k)}}
\left [ {{\tilde B_i} \over {1- \bar x_i x_R}} \right ]
$$
is a sum of terms of the form
$$
Res_{x_R}
\left [ {{ x_R^t RAT_t (x_{\pi(1)} , \dots ,x_{\pi(u-1)}) 
} \over
{x_R^{n+k-R} \bar x_R^{n+k+t+2}}} \right ] \quad ,
(t \geq 0 )  \quad .
\eqno(Gil'')
$$
But, each such term can be written as
$$
RAT_t (x_{\pi(1)} , \dots ,x_{\pi(u-1)})
Res_{x_R}
\left [ {{ x_R^{R+2t+2}
} \over
{x_R^{n+k+t+2} \bar x_R^{n+k+t+2}}} \right ] \quad ,
\quad (t \geq 0 )  \quad .
\eqno(Gil''')
$$
By crucial fact $\aleph_5'$, since $R+2t+2 \leq 2(n+k+t+1)$, we can 
perform $x_R \rightarrow \bar x_R$ in the residuand of $(Gil''')$, so
$$
Res_{x_R , x_{\pi(u+1)}, \dots ,  x_{\pi(k)}}
\left [ {{\tilde B_i} \over {1- \bar x_i x_R}} \right ]
=
Res_{x_R , x_{\pi(u+1)}, \dots ,  x_{\pi(k)}}
\left [ {{\tilde B_i} \over {1- \bar x_i x_R}} 
( x_R \rightarrow \bar x_R ) \right ] \quad.
$$
Applying $Res_{x_{\pi(1)} , \dots , x_{\pi(u-1)}}$ to both sides
of the above equation finishes up subcase Ib of the proof
of subsubsublemma $1.3.1.4$.
 
{\bf  Case II:} $i \in \{ \pi (1) , \dots , \pi (u-1) \}$ .
 
The proof is identical to Case II of $(sub)^4$lemma $1.4.1.4$
proved in act IV
(except that the form of $\tilde B_i$ is slightly more complicated
there, but the argument that reduces it to Case Ia goes verbatim).
This completes the proof
of $(sub)^3$lemma $1.3.1.4$. \halmos.
 
Going back to the proof of $(sub)^2$lemma $1.3.1$, we know by
the two $(sub)^3$lemmas $(1.3.1.3-4)$ that each part of
$(Celia)$, in the statement of $(sub)^3$lemma $1.3.1.2$, has
property $(Hadas)$, and by the additivity of this property, so does
the whole 
$f_{n,k} ( \epsilon_1(x_1) , \dots ,  x_R , \dots , \epsilon_k(x_k) )$.
But this is precisely what $(sub)^2$lemma
$1.3.1$ is saying. 
This completes the proof of $(sub)^2$lemma $1.3.1$. \halmos
 
As we noted at the very beginning of the proof of sublemma $1.3$,
subsublemma $1.3.1$, read from right to left, can be used repeatedly
to `unbar' $x_i$'s, until we get that
$$
Res_{x_{\pi(1)} , \dots , x_{\pi(k)}}
[ f_{n,k} ( \epsilon_1(x_1) , \dots , \epsilon_k(x_k) ) ] = 
Res_{x_{\pi(1)} , \dots , x_{\pi(k)}}
 [  f_{n,k} (x_1 , \dots , x_k ) ]  \quad  .
\eqno(Ktsat)
$$
 
Since  $  f_{n,k} (x_1 , \dots , x_k ) $ has a genuine
Laurent expansion, crucial fact $\aleph_4'$ assures
us that we can unscramble the  order of residue-taking on the right side
of $(Ktsat)$:
$$
Res_{x_{\pi(1)} , \dots , x_{\pi(k)}}
 [  f_{n,k} (x_1 , \dots , x_k ) ]  
=
Res_{x_1 , \dots , x_k }
 [  f_{n,k} (x_1 , \dots , x_k ) ]  
\quad  .
\eqno(Nimas)
$$
Combining $(Ktsat)$ with $(Nimas)$ yields $(1.3')$, which we
noted was equivalent to the statement of sublemma 1.3.
This completes the proof of sublemma $1.3$. \halmos
 
\vfill
\eject
%begin macros
\baselineskip=14pt
\parskip=10pt
\def\B{{\cal B}}
\def\S{{\cal S}}
\def \inv{\mathop{\rm inv} \nolimits}
\def \sgn{\mathop{\rm sgn} \nolimits} 
\def\hat{\widehat}
\def\tilde{\widetilde}
 \def\epsilon{\varepsilon}
\def\halmos{\hbox{\vrule height0.15cm width0.01cm\vbox{\hrule height
 0.01cm width0.2cm \vskip0.15cm \hrule height 0.01cm width0.2cm}\vrule
 height0.15cm width 0.01cm}}
\font\eightrm=cmr8  
\font\eighttt=cmtt8
\magnification=\magstephalf

\parindent=0pt
\overfullrule=0in
\headline={\rm  \ifodd\pageno  \RightHead  \else  \LeftHead  \fi}
\def\RightHead{\centerline{Proof  of the ASM Conjecture-Act  IV}}
\def\LeftHead{ \centerline{Doron Zeilberger}}
%\pageno=46
%end macros

\centerline{\bf Act IV. (GogTotal') IS UNFAZED BY THE SIGNED PERMUTATIONS'
ACTION} 
 
\qquad \qquad \qquad \qquad\qquad \qquad \qquad \qquad\qquad 
\qquad \qquad \qquad {\it  $\dots$  signed $\dots$-- (Daniel VI, 10) }
 
{\bf Sublemma 1.4 :} Let $F_{n,k} (x)= F_{n,k}( x_1 , \dots , x_k )$
be the residuand of $(GogTotal')$ , that is,
$$
F_{n,k} ( x_1 , \dots , x_k ):=
 { { \Phi_k ( x_1 , \dots , x_k ) } \over
   { \prod_{i=1}^k  x_i^{n+1} \bar x_i^{n+i+1} \prod_{ 1 \leq i < j \leq k}
(1-x_i x_j ) (1- \bar x_i x_j ) }
}
\qquad ,
$$
where $\Phi_k ( x_1 , \dots , x_k)$ is the polynomial defined in ($Gog_1$).
Let $g$ be any signed permutation in $W(\B_k)$, then
$$
Res_{x_1 , \dots , x_k}[ g F_{n,k} (x) ] = Res_{x_1 , \dots , x_k}
 [  F_{n,k} (x ) ]  \quad (=m_k (n)) \quad .
$$
{\eightrm  [ Type `{\eighttt S14(k,n):}' in ROBBINS, for specific k and n.]}
 
{\bf Proof}: As in the proof of $1.3$, by
renaming the variables, the statement of sublemma 1.4 is
equivalent to
 
{\bf Sublemma 1.4' :} Let $F_{n,k} (x)= F_{n,k}( x_1 , \dots , x_k )$
be the residuand of $(GogTotal')$, reproduced just above.
Let $\pi$ be any permutation of
$\{ 1 , \dots , k \}$, and let 
$\epsilon = (\epsilon_1 , \dots , \epsilon_k )$ be any sign-assignment.
Then
$$
Res_{x_{\pi(1)} , \dots , x_{\pi(k)}}
[ F_{n,k} ( \epsilon_1(x_1) , \dots , \epsilon_k(x_k) ) ] = 
Res_{x_1 , \dots , x_k} [  F_{n,k} (x_1 , \dots , x_k ) ]  \quad  .
\eqno(1.4')
$$
 
The proof will be by induction on the number of $\epsilon_i$'s
that are equal to $-1$. If there are no $-1$'s in $\epsilon$, then
$1.4'$ follows immediately from crucial fact $\aleph_4$,
since then the residuand $F_{n,k}( x_1 , \dots , x_k )$ has a Laurent
series. Otherwise, let $u$ be the smallest integer $i$ ($1 \leq i \leq k$)
such that $\epsilon_{\pi(i)} = -1$.  
The inductive step is provided by
the following subsublemma, that asserts that it is permissible to change
$\epsilon_{\pi(u)}$ from $-1$ to $+1$ 
inside the residuand of the iterated residue on the left side 
of $(1.4')$ without affecting its value.
 
{ \bf Subsublemma 1.4.1:} Let $F_{n,k} ( x_1 , \dots , x_k )$ be
the residuand of $(GogTotal')$, let $\pi$ be any permutation of
$\{ 1 , \dots , k \}$, and let 
$\epsilon = ( \epsilon_1 , \dots , \epsilon_k )$ be any sign-assignment
with at least one $\epsilon_i$ equal to $-1$.
Let $u$ be the smallest $i$ ($1 \leq i \leq k$) such that
$\epsilon_{\pi(i)} = -1$. (i.e., we have
$\epsilon_{\pi(1)}= \dots  = \epsilon_{\pi(u-1)}= +1$,
and $\epsilon_{\pi(u)}= -1$). Then:
$$
Res_{x_{\pi(1)} , \dots , x_{\pi(k)}}
[ F_{n,k} ( \epsilon_1(x_1) , \dots ,  x_{\pi(u)} , \dots ,
  \epsilon_k(x_k) ) ] = 
Res_{x_{\pi(1)} , \dots , x_{\pi(k)}}
[ F_{n,k} ( \epsilon_1(x_1) , \dots ,  \bar x_{\pi(u)} , 
\dots , \epsilon_k(x_k) ) ] 
\quad  .
$$
{\eightrm  [ Type `{\eighttt S141(k,n):}' in ROBBINS, for specific k and n.]}
 
{\bf Proof:} Let's abbreviate $R:= \pi (u)$, and
$z_i := \epsilon_i(x_i)$ ( $1 \leq i \leq k , i \neq R $).
Recall that this means that $z_i=x_i$ if $\epsilon_i=+1$,
and $z_i = \bar x_i = 1 - x_i$ if $\epsilon_i = -1$. Also put $z_R=x_R$.
By the $W(\B_k)$-antisymmetry of $\Phi_k(x)$, we have
$$
F_{n,k} ( \epsilon_1(x_1) , \dots ,  x_R , \dots , \epsilon_k(x_k) )
=
 { { \pm \Phi_k ( x_1 , \dots , x_k ) } \over
   { x_R^{n+1} \bar x_R^{n+R+1} \prod_{ {{j=1} \atop {j \neq R } }}^{k} 
  z_j^{n+1} \bar z_j^{n+j+1}
 \prod_{ 1 \leq i < j \leq k}
(1-z_i z_j ) (1- \bar z_i z_j )
    }
}
$$
$$
=
{ {1} \over { x_R^{n+1} \bar x_R^{n+R+1} } }
\prod_{ {{j=1} \atop {j \neq R } }}^{k}
{{1} \over {z_j^{n+1} \bar z_j^{n+j+1} }}
\prod_{ {{1 \leq r < s \leq k } \atop {r,s \neq R} }}
{{1} \over {(1- z_r z_s )(1- \bar z_r z_s ) }} \cdot
$$
$$
\left \{
{ {\pm \Phi_k ( x_1 , \dots , x_k ) } \over
   { \prod_{{{i=1} \atop { i \neq R}} }^{k} (1- z_i x_R ) 
\prod_{i=1}^{R-1} ( 1- \bar z_i x_R )
\prod_{i=R+1}^{k} ( 1-  z_i \bar x_R )}
} \right \} \quad .
\eqno(Jane)
$$
 
The next step is to perform a partial-fraction decomposition,
with respect to the variable $x_R$, of the fraction inside
$(Jane)$'s braces.
 
{\bf Subsubsublemma 1.4.1.1:}
$$
{ {\pm \Phi_k ( x_1 , \dots , x_k ) } \over
   { \prod_{{{i=1} \atop { i \neq R}} }^{k} 
(1- z_i x_R ) \prod_{i=1}^{R-1} ( 1- \bar z_i x_R )
\prod_{i=R+1}^{k} ( 1-  z_i \bar x_R )}
} 
$$
$$
=  P_R +
\sum_{ {{i=1} \atop {i \neq R } }}^{k} 
{  {A_i} \over {1- z_i x_R} } 
+ \sum_{i=1}^{R-1} {  {B_i} \over {1- \bar z_i x_R} }
+ \sum_{i=R+1}^{k} {  {B_i} \over {1- z_i \bar x_R } } \quad ,
\eqno(Tamar)
$$
where $P_R$ is a polynomial of degree $\leq 2k-1$ in $x_R$,
with coefficients that are rational functions of the other
variables $(x_1 , \dots , x_{R-1} , x_{R+1} , \dots , x_k )$, and
the $A_i$'s and $B_i$'s are as follows.
 
When $1 \leq i < R$,
$$
A_i =  {{\pm ( z_i -2 ) (1- z_i^2 ) (1- z_i \bar z_i ) \Phi_{k-2} (x_1 , \dots ,
\hat {x_i} , \dots , \hat x_R , \dots , x_k ) }
\over {z_i^k} } \cdot
$$
$$
\left \{ 
\prod_{ {{j=1} \atop {j \neq i,R}} }^{k}
(1- z_i z_j )(1- \bar z_i z_j ) (1- z_i \bar z_j )
\right \}
\prod_{j=R+1}^k ( z_j + z_i -1)
\prod_{ {{j=1} \atop {j \neq i}} }^{R-1} (1- \bar z_i \bar z_j )
\quad .
\eqno(Ai1)
$$
When $ R < i \leq k$, we have:
$$
A_i =  {{\pm ( 1 - 2 z_i ) (1- z_i^2 ) (1- z_i \bar z_i ) 
\Phi_{k-2} (x_1 , \dots ,
\hat x_R , \dots , \hat {x_i} , \dots , x_k ) }
\over {z_i^{k+1} } }\cdot
$$
$$
\left \{ 
\prod_{ {{j=1} \atop {j \neq i,R}} }^{k}
(1- z_i z_j )(1- \bar z_i z_j ) (1- z_i \bar z_j )
\right \}
\prod_{ {{j=R+1} \atop {j \neq i}} }^{k} ( z_j + z_i -1)
\prod_{j=1}^{R-1} (1- \bar z_i \bar z_j )
\quad .
\eqno(Ai2)
$$
As for the $B_i$, when $1 \leq i < R$, we have
$$
B_i =  {{\pm (  \bar z_i -2 ) (1- \bar z_i^2 ) (1- z_i \bar z_i ) 
\Phi_{k-2} (x_1 , \dots ,
\hat {x_i} , \dots , \hat x_R , \dots , x_k ) }
\over {\bar z_i^{k} }} \cdot
$$
$$
\left \{ 
\prod_{ {{j=1} \atop {j \neq i,R}} }^{k}
(1- z_i z_j )(1- \bar z_i z_j ) (1- \bar z_i \bar z_j )
\right \}
\prod_{j=R+1}^k ( z_j - z_i)
\prod_{ {{j=1} \atop {j \neq i}} }^{R-1} (1- z_i \bar z_j )
\quad .
\eqno(Bi1)
$$
Finally, when $R < i \leq k$, we have:
$$
B_i =  {{\pm ( 1 - 2 z_i ) (1- z_i^2 ) (1- z_i \bar z_i ) 
\Phi_{k-2} (x_1 , \dots ,
\hat x_R , \dots , \hat {x_i} , \dots , x_k ) }
\over {z_i^{k+1} }} \cdot
$$
$$
\left \{ 
\prod_{ {{j=1} \atop {j \neq i,R}} }^{k}
(1- z_i z_j )(1-  z_i \bar z_j ) ( z_i + z_j -1)
\right \}
\prod_{ {{j=R+1} \atop {j \neq i}} }^{k}  (1- \bar z_i z_j )
\prod_{j=1}^{R-1} ( z_i - z_j )
\quad .
\eqno(Bi2)
$$
{\eightrm  [ Type `{\eighttt S1411(k,R):}' in ROBBINS, for specific k and R.]}
 
{\bf Proof:} By Freshman Calculus (or rather the algebraic part of it),
we know that the left side of $(Tamar)$, when viewed as a rational function
of $x_R$, can be written as on its right side, for
{\it some } polynomial $P_R$, and {\it constants} (in the present
context, i.e. quantities  not involving $x_R$) $A_i$, $B_i$.
The degree of $\Phi_k$, in $x_R$
(or for that matter, in any variable),
from its definition in $(Gog_1)$, is
$\leq 4k-3$, and the degree of the denominator of the left side of
$(Tamar)$ is exactly $2(k-1)$. Hence the degree of $P_R$, in $x_R$, is
$\leq (4k-3)-(2k-2)=2k-1$. The coefficients of $P_R$ are certain
rational functions of the remaining variables whose exact form
does not interest us. The exact forms of $A_i$ and $B_i$ 
{\it should} interest us,  since they are used heavily in the proof.
So let's brace ourselves and embark on partial-fractioning.
 
Multiplying both sides of $(Tamar)$ by the denominator of
its left side yields the equation:
$$
\pm \Phi_k ( x_1 , \dots , x_k ) 
=
 \left \{   \prod_{{{j=1} \atop {j \neq R} }}^{k} 
(1- z_j x_R ) \prod_{j=1}^{R-1} ( 1- \bar z_j x_R )
\prod_{j=R+1}^{k} ( 1-  z_j \bar x_R )  \right \}
P_R 
$$
$$
+ \sum_{ {{i=1} \atop {i \neq R } }}^{k} 
 \left \{   \prod_{{{j=1} \atop {j \neq R,i} }}^{k} 
(1- z_j x_R ) \prod_{j=1}^{R-1} ( 1- \bar z_j x_R )
\prod_{j=R+1}^{k} ( 1-  z_j \bar x_R )  \right \}
A_i
$$
$$ 
+
\sum_{i=1}^{R-1}
 \left \{   \prod_{{{j=1} \atop {j \neq R} }}^{k} 
(1- z_j x_R ) \prod_{{{j=1} \atop {j \neq i} }}^{R-1} ( 1- \bar z_j x_R )
\prod_{j=R+1}^{k} ( 1-  z_j \bar x_R )  \right \}
B_i  
$$
$$ 
+
\sum_{i=R+1}^{k}
 \left \{   \prod_{{{j=1} \atop {j \neq R} }}^{k} 
(1- z_j x_R ) \prod_{j=1}^{R-1} ( 1- \bar z_j x_R )
\prod_{{{j=R+1} \atop {j \neq i} }}^{k} ( 1-  z_j \bar x_R )  \right \}
B_i  
\quad .
\eqno(Tamar')
$$
 
To decipher $A_i$ ($1 \leq i \leq k$, $i \neq R$), we plug in
$x_R=1/z_i$ on both sides of $(Tamar')$ to get:
$$
A_i=
{  {\pm \Phi_k ( x_1 , \dots , x_k ) \vert_{x_R=1/z_i}} \over
{ \left \{   \prod_{{{j=1} \atop {j \neq R,i} }}^{k} 
(1- z_j/z_i ) \prod_{j=1}^{R-1} ( 1- \bar z_j/z_i )
\prod_{j=R+1}^{k} ( 1-  z_j (1- 1/z_i) )  \right \} } }
$$
$$
= {{ \pm  z_i^{2k-3} \Phi_k ( x_1 , \dots , x_k ) \vert_{x_R=1/z_i}} \over
{ \left \{   \prod_{{{j=1} \atop {j \neq R,i} }}^{k} 
( z_i - z_j ) \prod_{j=1}^{R-1} ( z_i- \bar z_j )
\prod_{j=R+1}^{k} ( 1-  \bar z_i \bar z_j  )  \right \}} } \quad.
\eqno(PP1)
$$
 
To decipher $B_i$ ($1 \leq i <R$ ), we plug in
$x_R=1/ \bar z_i$ on both sides of $(Tamar')$ to get:
$$
B_i=
{  {\pm \Phi_k ( x_1 , \dots , x_k ) \vert_{x_R=1/ \bar z_i}} \over
{ \left \{   \prod_{{{j=1} \atop {j \neq R} }}^{k} 
(1- z_j/ \bar z_i ) \prod_{{{j=1} \atop {j \neq i}}}^{R-1} 
( 1- \bar z_j/ \bar z_i )
\prod_{j=R+1}^{k} ( 1-  z_j (1- 1/\bar z_i) )  \right \} } }
$$
$$
= {{ \pm \bar z_i^{2k-3} \Phi_k ( x_1 , \dots , x_k ) \vert_{x_R=1/ \bar z_i}} 
\over
{ \left \{   \prod_{{{j=1} \atop {j \neq R} }}^{k} 
( \bar z_i - z_j ) \prod_{{{j=1} \atop {j \neq i}}}
^{R-1} ( \bar z_i- \bar z_j )
\prod_{j=R+1}^{k} ( 1-   z_i \bar z_j  )  \right \}} } 
\qquad (1 \leq i <R) \quad.
\eqno(PP2)
$$
 
To decipher $B_i$ ($R < i \leq k$ ), we plug in
$\bar x_R=1/ z_i$, i.e. $x_R=1- 1/z_i = - \bar z_i/z_i$,
on both sides of $(Tamar')$ to get:
$$
B_i=
{  {\pm \Phi_k ( x_1 , \dots , x_k ) \vert_{\bar x_R=1/z_i}} \over
{ \left \{   \prod_{{{j=1} \atop {j \neq R} }}^{k} 
(1+z_j \bar z_i /z_i ) \prod_{j=1}^{R-1} ( 1+ \bar z_j \bar z_i / z_i )
\prod_{{{j=R+1} \atop {j \neq i}}}^{k} ( 1-  z_j /z_i )  \right \} } }
$$
$$
= {{ \pm  z_i^{2k-3} \Phi_k ( x_1 , \dots , x_k ) \vert_{\bar x_R=1/ z_i}} 
\over
{ \left \{   \prod_{{{j=1} \atop {j \neq R} }}^{k} 
( 1- \bar z_i \bar z_j ) \prod_{j=1}^{R-1} ( 1- \bar z_i z_j )
\prod_{{{j=R+1} \atop {j \neq i}}}^{k} ( z_i - z_j  )  \right \}} } 
\qquad ( R < i  \leq k) \quad.
\eqno(PP3)
$$
We now need:
 
{\bf Subsubsubsublemma 1.4.1.1.1:} For $1 \leq R \neq i \leq k$:
$$
\Phi_k ( x_1 , x_2 , \dots , x_k ) \vert_{x_R=1/z_i} =
\pm {{(z_i -2 ) (1- 2 z_i ) (1- z_i^2) (1- z_i \bar z_i ) 
\cdot \Phi_{k-2}( x_1, \dots,\hat x_i , \dots, \hat x_R , \dots ,
x_k ) }
\over { z_i^{3k-3} }} \cdot
$$
$$
\prod_{ {{j=1} \atop { j \neq i,R }}}^{k}
\{ (1- z_i z_j ) ( 1 - \bar z_i z_j )( 1- z_i \bar z_j )
(1- \bar z_i \bar z_j ) (z_j - z_i ) (z_j + z_i -1) \} \quad .
\eqno(1.4.1.1.1)
$$
{\eightrm  [ Type `{\eighttt S14111(k,R,i):}' in ROBBINS, for specific k,R, and 
i.]}
 
{\bf Proof:} By the anti-symmetry of $\Phi_k$ with respect to the symmetric
group, it suffices to consider $R=1$. By its anti-symmetry w.r.t. the
sign-action, it suffices to consider the case $z_i=x_i$
(i.e. $\epsilon_i = +1$). But this follows from multiplying the respective
sides of the identities of
$(sub)^4$lemma  $1.3.1.1.1$ (proved above) and subsublemma 
$1.5.2$ (proved below).
This completes the proof of
$(sub)^4$lemma 1.4.1.1.1. \halmos .
 
Substituting $(1.4.1.1.1)$ into $(PP1)$, $(PP2)$, and $(PP3)$, and cancelling
out whenever possible, finishes up the proof of $(sub)^3$lemma $1.4.1.1$.
For $(PP3)$ we must note that $\Phi_k$ remains the same, except for a sign change,
when one replaces $x_R$ by $\bar x_R$, so we can also apply $1.4.1.1.1$ with
$x_R$ replaced by $\bar x_R$.  This completes the proof of
$(sub)^3$lemma 1.4.1.1. \halmos
 
For the 
rest of the proof of $(sub)^2$lemma $1.4.1$, we need some new notation.
Let $POL$ stand for ``some polynomial of $(x_1 , \dots , x_k )$'',
and let $RAT$ stand for ``some rational function of $(x_1 , \dots , x_k )$''.
$POL( \hat x_R )$ stands for such a polynomial that does not involve $x_R$.
$POLRAT( x_R; x_1 , \dots , \hat x_R , \dots , x_k)$, or $POLRAT(x_R \,\, ; \, \, \dots )$
for short, stands for a polynomial in $x_R$ with coefficients that
are rational functions of the rest of the variables.
$POLRAT(x_i ; \hat x_R)$ stands for a polynomial in $x_i$
with coefficients that are rational functions of the remaining
variables {\it and} that do not involve $x_R$.

Inserting $(Tamar)$(in the statement of $1.4.1.1$) 
into $(Jane)$ is going to yield the following $(sub)^3$lemma:
 
{\bf Subsubsublemma 1.4.1.2:}
$$
F_{n,k} ( \epsilon_1(x_1) , \dots ,  x_R , \dots , \epsilon_k(x_k) )
$$
$$
=  {{\tilde P_R} \over {x_R^{n+1} \bar x_R^{n+R+1} } } +
\sum_{ {{i=1} \atop {i \neq R } }}^{k} 
{  {\tilde A_i} \over {1- z_i x_R} } 
+ \sum_{i=1}^{R-1} {  {\tilde B_i} \over {1- \bar z_i x_R} }
+ \sum_{i=R+1}^{k} {  {\tilde B_i} \over {1- z_i \bar x_R } } \quad ,
\eqno(Celia)
$$
 
where $\tilde P_R$ is a polynomial of degree $\leq 2k-1$ in $x_R$,
with coefficients that are rational functions of the rest of the
variables, and
$\tilde A_i$, and $\tilde B_i$ can be written as
follows.
 
For $1 \leq i < R$, we have
$$
\tilde A_i = {
{POL(\hat x_R) } \over { x_R^{n+1} \bar x_R^{n+R+1}
                        z_i^{n+k+1} \bar z_i^{n+i+1} }} 
\cdot \prod_{ {{j=1} \atop {j \neq i,R} }}
^k {{1} \over {z_j^{n+1} \bar z_j^{n+j+1} }} \cdot
\prod_{ {{1 \leq r < s \leq k } \atop {r,s \neq R,i} }}
{{1} \over {(1- z_r z_s )(1- \bar z_r z_s ) }}
$$
$$
={{POLRAT(x_i ;\hat x_R) } \over
   {x_R^{n+1} \bar x_R^{n+R+1} z_i^{n+k+1} \bar z_i^{n+i+1} } }
\qquad  \quad .
\eqno(\tilde Ai1)
$$
 
For $R < i \leq k$, we have
$$
\tilde A_i = {{POL(\hat x_R) } \over { x_R^{n+1} \bar x_R^{n+R+1}
 z_i^{n+k+2} \bar z_i^{n+i+1} }} \cdot \prod_{ {{j=1} \atop {j \neq i,R} }}
^k {{1} \over {z_j^{n+1} \bar z_j^{n+j+1} }} \cdot
\prod_{ {{1 \leq r < s \leq k } \atop {r,s \neq R,i} }}
{{1} \over {(1- z_r z_s )(1- \bar z_r z_s ) }}
$$ 
$$
={{POLRAT(x_i ; \hat x_R) } \over
   {x_R^{n+1} \bar x_R^{n+R+1} z_i^{n+k+2} \bar z_i^{n+i+1} } }
\qquad   \quad .
\eqno(\tilde Ai2)
$$
 
For $1 \leq i < R$, we have
$$
\tilde B_i = {{POL(\hat x_R) } \over { x_R^{n+1} \bar x_R^{n+R+1}
 z_i^{n+1} \bar z_i^{n+i+k+1} }} \cdot \prod_{ {{j=1} \atop {j \neq i,R} }}
^k {{1} \over {z_j^{n+1} \bar z_j^{n+j+1} } } \cdot
\prod_{ {{1 \leq r < s \leq k } \atop {r,s \neq R,i} }}
{{1} \over {(1- z_r z_s )(1- \bar z_r z_s ) }}
$$
$$
={{POLRAT(x_i ; \hat x_R) } \over
   {x_R^{n+1} \bar x_R^{n+R+1} z_i^{n+1} \bar z_i^{n+i+k+1} } }
 \quad .
\eqno(\tilde Bi1)
$$
 
Finally, for $R < i \leq k$, we have
$$
\tilde B_i = {{POL(\hat x_R )  } \over { x_R^{n+1} \bar x_R^{n+R+1}
 z_i^{n+k+2} \bar z_i^{n+i+1} }} \cdot \prod_{ {{j=1} \atop {j \neq i,R} }}
^k {{1} \over {z_j^{n+1} \bar z_j^{n+j+1} }} \cdot
\prod_{ {{1 \leq r < s \leq k } \atop {r,s \neq R,i} }}
{{1} \over {(1- z_r z_s )(1- \bar z_r z_s ) }}
$$
$$
={{POLRAT( x_i ; \hat x_R)} \over
   {x_R^{n+1} \bar x_R^{n+R+1} z_i^{n+k+2} \bar z_i^{n+i+1} } } \quad .
\eqno(\tilde Bi2)
$$
Furthermore, all the $POLRATS$ are of degree $\leq 2k+1$ in $x_i$
(we actually only use this fact for the $POLRAT$ of $(\tilde Bi1)$).
 
{\eightrm  [ Type `{\eighttt S1412(k,n,eps,R):}' in ROBBINS, 
for specific k,n,eps, and R.]}
 
{\bf Proof:} Insert $(Tamar)$ into $(Jane)$, bringing the quantity
in front of the sums inside, and for each $A_i$ and $B_i$, split
the double product:
$$
\prod_{ {{1 \leq r < s \leq k } \atop {r,s \neq R} }}
{{1} \over {(1- z_r z_s )(1- \bar z_r z_s ) }}=
\prod_{ {{1 \leq r < s \leq k } \atop {r,s \neq R,i} }}
{{1} \over {(1- z_r z_s )(1- \bar z_r z_s ) }} \cdot
\prod_{ {{j=1 } \atop {j \neq R} }}^{i-1}
{{1} \over {(1- z_j z_i )(1- \bar z_j z_i ) }}
\prod_{ {{j=i+1 } \atop {j \neq R} }}^{k}
{{1} \over {(1- z_i z_j )(1- \bar z_i z_j ) }} ,
\eqno(Split)
$$
where the condition $j \neq R$ is vacuous in exactly one of the two
last products
according to whether $i<R$ or $i>R$. Finally perform 
all possible  cancellations.
This completes  the proof of $(sub)^3$lemma
$1.4.1.2$. \halmos
 
We will say that a rational function $h=h( x_1 , \dots , x_k )$ has
property $(Hadas)$ if:
$$
Res_{x_{\pi(1)} , \dots , x_{\pi(k)}} [h] =
Res_{x_{\pi(1)} , \dots , x_{\pi(k)}} [h( x_R \rightarrow \bar x_R )] \quad ,
\eqno(Hadas)
$$
where  $h( x_R \rightarrow \bar x_R )$ is the rational function obtained from $h$
by replacing the variable $x_R$ by $\bar x_R$ ($=1- x_R $.)
 
The statement of subsublemma $1.4.1$ is equivalent to saying that
$$
F_{n,k} ( \epsilon_1(x_1) , \dots ,  x_R , \dots ,
  \epsilon_k(x_k) ) 
$$
has property $(Hadas)$. Property $(Hadas)$ is obviously additive
(since the action of taking residue, and performing $x_R \rightarrow \bar x_R$ are),
so in order to complete the proof of $1.4.1$ (and hence of $1.4$), it
would suffice to show that each term in $(Celia)$ has property
$(Hadas)$. This will be accomplished by the following series
of subsubsublemmas ($1.4.1.3-6$.)
 
{\bf Subsubsublemma 1.4.1.3:}
Let $\tilde P_R$ be as in the statement of $1.4.1.2$, i.e.
a polynomial of degree $\leq 2k-1$ in $x_R$ with coefficients
that are rational functions of the remaining variables.
Then ${{\tilde P_R} \over {x_R^{n+1} \bar x_R^{n+R+1} } }$ 
has property $(Hadas)$:
$$
Res_{x_{\pi(1)} , \dots , x_{\pi(k)}} 
\left [{{\tilde P_R} \over {x_R^{n+1} \bar x_R^{n+R+1} } } \right ] =
Res_{x_{\pi(1)} , \dots , x_{\pi(k)}} 
\left [ {{\tilde P_R} \over {x_R^{n+1} \bar x_R^{n+R+1} } } 
( x_R \rightarrow \bar x_R ) \right ] 
\quad .
\eqno(Hadas_1)
$$
{\eightrm  [ Type `{\eighttt S1413(k,n):}' in ROBBINS, for specific k and n.]}
 
{\bf Proof:} Recall that $R=\pi(u)$. Since $\tilde P_R$ is
a $POLRAT( x_R ; x_1 , \dots , \hat x_R , \dots , x_k )$, we have that
$$
Res_{x_{\pi(u+1)} , \dots , x_{\pi(k)}} 
\left [{{\tilde P_R} \over {x_R^{n+1} \bar x_R^{n+R+1} } } \right ]
=
Res_{x_{\pi(u+1)} , \dots , x_{\pi(k)}} 
\left [ {{POLRAT(x_R; x_{\pi(1)} , \dots , \hat x_R , x_{\pi(u+1)},
\dots , x_{\pi(k)})} \over {x_R^{n+1} \bar x_R^{n+R+1} } }  \right ]
$$
$$
={{POLRAT(x_R; x_{\pi(1)} , \dots , x_{\pi(u-1)})} 
\over {x_R^{n+1} \bar x_R^{n+R+1} } } 
={{x_R^R \cdot POLRAT(x_R; x_{\pi(1)} , \dots , x_{\pi(u-1)})} 
\over {x_R^{n+R+1} \bar x_R^{n+R+1} } } \quad .
$$
 
Since  the degree of
$x_R^R POLRAT(x_R \,\, ; \, \, \dots)$, in $x_R$, is 
$\leq R+ (2k-1) \leq 2(n+R)$, it follows from crucial fact $\aleph_5'$ that 
$$
Res_{x_R} \left [ {{POLRAT ( x_R \,\, ; \, \, \dots) } \over 
{x_R^{n+1} \bar x_R^{n+R+1} } } \right ] =
Res_{x_R} \left [ {{POLRAT ( x_R \,\, ; \, \, \dots)} \over 
{x_R^{n+1} \bar x_R^{n+R+1} } }
( x_R \rightarrow \bar x_R ) \right ] \quad .
$$
Applying $Res_{x_{\pi(1)} , \dots , x_{\pi(u-1)}}$ to both sides of
the above equation establishes that 
${{\tilde P_R} \over {x_R^{n+1} \bar x_R^{n+R+1} } } $ indeed
has property $(Hadas)$, so that $(Hadas_1)$ holds. This
completes the proof of subsubsublemma $1.4.1.3$. \halmos
 
{\bf Subsubsublemma 1.4.1.4:}
Let $\tilde A_i$ be as in the statement of $1.4.1.2$, i.e.
given by $(\tilde Ai1)$ and $(\tilde Ai2)$ according to whether
$i<R$ or $i > R$ respectively. Then:
$$
Res_{x_{\pi(1)} , \dots , x_{\pi(k)}} 
[ {{\tilde A_i} \over {1- z_i x_R}} ] =
Res_{x_{\pi(1)} , \dots , x_{\pi(k)}} 
[ {{\tilde A_i} \over {1- z_i x_R}} ( x_R \rightarrow \bar x_R )]
\quad.
\eqno(Hadas_2)
$$
{\eightrm  [ Type `{\eighttt S1414(k,n):}' in ROBBINS, for specific k and n.]}
 
{\bf Proof:} Recall that $u=\pi^{-1}(R)$.
$(\tilde Ai1)$ that gives $\tilde A_i$ for $i <R$, and
$(\tilde Ai2)$ that gives $\tilde A_i$ for $i >R$ are of the same form,
{\it except} that the power of $z_i$ in the denominator of
$(\tilde Ai2)$ is $1$ more (namely $n+k+2$) than its counterpart
($n+k+1)$ in $(\tilde Ai1)$. Multiplying top and bottom
of $(\tilde Ai1)$ by $z_i$ makes them exactly of the same form.
 
{\bf  Case I:} $i \in \{ \pi (u+1) , \dots , \pi (k) \}$ .
 
Let $v := \pi^{-1} (i)$, so that $i= \pi(v)$, and $v > u$.
 
We have that:
$$
{{\tilde A_i} \over {1- z_i x_R}} =
{
 {
  POLRAT( x_i ; x_{\pi(1)} , \dots , \hat x_R ,
    \dots , x_{\pi(v-1)} , 
  x_{\pi(v+1)} , \dots , x_{\pi(k)} )
  } \over
   {x_R^{n+1} \bar x_R^{n+R+1} z_i^{n+k+2} \bar z_i^{n+i+1} (1- z_i x_R) }}
\quad .
$$
 
Applying $Res_{x_{\pi(v+1)} , \dots , x_{\pi(k)}}$ to it gives:
$$
Res_{x_{\pi(v+1)} , \dots , x_{\pi(k)}} \left [
{{\tilde A_i} \over {1- z_i x_R}} \right ]
={{
POLRAT(x_i ;x_{\pi(1)} , \dots , \hat x_R,
 \dots , x_{\pi(v-1)}) } \over
   {x_R^{n+1} \bar x_R^{n+R+1} z_i^{n+k+2} \bar z_i^{n+i+1} (1- z_i x_R) }}
  \quad .
$$
Applying $Res_{x_i}$ to the above gives:
$$
Res_{x_i , x_{\pi(v+1)} , \dots , x_{\pi(k)}}
\left [ {{\tilde A_i} \over {1- z_i x_R}} \right ]
=Res_{x_i} \left [ 
{{
POLRAT(x_i ; x_{\pi(1)} , \dots , \hat x_R,
 \dots , x_{\pi(v-1)}) } \over
   {x_R^{n+1} \bar x_R^{n+R+1} z_i^{n+k+2} \bar z_i^{n+i+1} (1- z_i x_R) }}
 \right ]
$$
$$
={{1} \over {x_R^{n+1} \bar x_R^{n+R+1}}} \cdot
Res_{x_i} \left [ 
{{
POLRAT(x_i; x_{\pi(1)} , \dots , \hat x_R,
 \dots , x_{\pi(v-1)} )} \over
   { z_i^{n+k+2} \bar z_i^{n+i+1} (1- z_i x_R) }}
 \right ] \quad .
\eqno(Doron)
$$
 
We now distinguish two subcases:
 
{\bf Subcase Ia:} $\epsilon_i = +1$, that is: $z_i = x_i$.
 
In this subcase, the right side of $(Doron)$ equals:
$$
 {{1 } \over
{x_R^{n+1} \bar x_R^{n+R+1}}}
\cdot
Coeff_{x_i^{n+k+1}} \left [ {{
 POLRAT(x_i ;x_{\pi(1)} , \dots ,
 \hat x_R  , \dots , x_{\pi(v-1)} )
 } \over
   { \bar x_i^{n+i+1} (1- x_i x_R) }}
 \right ]
  \quad .
\eqno(Gil)
$$
Expanding  everything involving $x_i$ as a power series in $x_i$:
$$
POLRAT(x_i;x_{\pi(1)} , \dots , \hat x_R , \dots , x_{\pi(v-1)}) =
\sum_{t=0}^{degree}
RAT_t (x_{\pi(1)} , \dots , \hat x_R , \dots , x_{\pi(v-1)} ) x_i^t \quad ,
$$
$$
{{1} \over {\bar x_i^{n+i+1} }} = \sum_{t=0}^{\infty} NUMBER_t \cdot x_i^t
\quad ,
$$
and, most importantly:
$$
(1- x_i x_R)^{-1} = \sum_{t=0}^{\infty} x_R^t x_i^t
\quad .
$$
 
Multiplying out, and collecting the coefficient of $x_i^{n+k+1}$, gives that
$Coeff_{x_i^{n+k+1}}$ in $(Gil)$ is a certain {\it polynomial }
in $x_R$, of degree $\leq n+k+1$, with coefficients that are rational
functions of $(x_{\pi(1)} , \dots , \hat x_R , \dots , x_{\pi(v-1)})$.
In other words it is a $POLRAT(x_R; x_{\pi(1)} , \dots , x_{\pi(v-1)})$
with $deg_{x_R} \leq n+k+1$. Hence,
$$
Res_{x_R, x_{\pi(u+1)} , \dots ,  x_{\pi(k)}}
\left [ {{\tilde A_i} \over {1- x_i x_R}} \right ]
=
Res_{x_R,x_{\pi(u+1)} , \dots ,  x_{\pi(v-1)}}
\left [ 
{{POLRAT(x_R ; x_{\pi(1)} , \dots , \hat x_R, \dots , x_{\pi(v-1)}) 
} \over
{x_R^{n+1} \bar x_R^{n+R+1}}} \right ] \quad ,
\eqno(Gil')
$$
where the $POLRAT$  in $(Gil')$
has degree $\leq n+k+1$ in $x_R$.
Applying $Res_{x_{\pi(u+1)} , \dots , x_{\pi(v-1)}}$ to
the residuand of $(Gil')$ gets rid of the dependence on
$(x_{\pi(u+1)} , \dots , x_{\pi(v-1)})$, so we have:
$$
Res_{x_R , x_{\pi(u+1)} , \dots ,  x_{\pi(k)}}
\left [ {{\tilde A_i} \over {1- x_i x_R}} \right ]
=
Res_{x_R}
\left [ {{POLRAT(x_R ; x_{\pi(1)} , \dots , x_{\pi(u-1)}) 
} \over
{x_R^{n+1} \bar x_R^{n+R+1}}} \right ] \quad ,
\eqno(Gil'')
$$
where $POLRAT(x_R; x_{\pi(1)}, \dots , x_{\pi(u-1)} )$ 
has degree $\leq n+k+1$ in $x_R$.
 
But, the right of $(Gil'')$ can be rewritten as:
$$
Res_{x_R}
\left [ {{
POLRAT( x_R ; x_{\pi(1)} , \dots , x_{\pi(u-1)}) 
} \over
{x_R^{n+1} \bar x_R^{n+R+1}}} \right ] =
Res_{x_R}
\left [ {{
x_R^R \cdot POLRAT(x_R; x_{\pi(1)} , \dots , x_{\pi(u-1)}) 
} \over
{x_R^{n+R+1} \bar x_R^{n+R+1}}} \right ] \quad .
\eqno(Gil''')
$$
 
By crucial fact $\aleph_5'$, since $R+(n+k+1) \leq 2(n+R)$
(recall that always $n \geq k$ and $R \geq 1$), we can 
perform $x_R \rightarrow \bar x_R$ in the residuand of $(Gil''')$, so
$$
Res_{x_R , x_{\pi(u+1)}, \dots ,  x_{\pi(k)}}
\left [ {{\tilde A_i} \over {1- x_i x_R}} \right ]
=
Res_{x_R , x_{\pi(u+1)}, \dots ,  x_{\pi(k)}}
\left [ {{\tilde A_i} \over {1- x_i x_R}} 
( x_R \rightarrow \bar x_R ) \right ] \quad.
$$
 
Applying $Res_{x_{\pi(1)} , \dots , x_{\pi(u-1)}}$ to both sides
of the above equation finishes up subcase Ia of the proof
of $(sub)^3$lemma $1.4.1.4$.
 
{\bf Subcase Ib:} $\epsilon_i = -1$, that is: $z_i = \bar x_i$.
 
In this subcase, the right side of $(Doron)$ equals:
$$
 {{1 } \over
{x_R^{n+1} \bar x_R^{n+R+1}}}
\cdot
Coeff_{x_i^{n+i}} \left [ {{
 POLRAT(x_i ; x_{\pi(1)} , \dots ,
 \hat x_R  , \dots , x_{\pi(v-1)})
 } \over
   { \bar x_i^{n+k+2} (1- \bar x_i x_R) }}
 \right ]
  \quad .
\eqno(Gil)
$$
Expanding  everything involving $x_i$ as a power series in $x_i$:
$$
POLRAT(x_i ; x_{\pi(1)} , \dots , \hat x_R , \dots , x_{\pi(v-1)}) =
\sum_{t=0}^{degree}
RAT_t (x_{\pi(1)} , \dots , \hat x_R , \dots , x_{\pi(v-1)} ) x_i^t \quad ,
$$
$$
{{1} \over {\bar x_i^{n+k+2} }} = \sum_{t=0}^{\infty} NUMBER_t \cdot x_i^t
\quad ,
$$
and, most importantly:
$$
(1- \bar x_i x_R)^{-1} = (1- x_R + x_i x_R )^{-1} =
( \bar x_R + x_i x_R )^{-1} =
\sum_{t=0}^{\infty} (-1)^t {{x_R^t} \over {\bar x_R^{t+1}}} x_i^t \quad ,
$$
multiplying out, and collecting the coefficient of $x_i^{n+i}$, gives that
$Coeff_{x_i^{n+i}}$ in $(Gil)$ is a certain sum of terms of the form
$$
RAT_t (x_{\pi(1)} , \dots , \hat x_R, \dots , x_{\pi(v-1)}) 
{{x_R^t} \over {\bar x_R^{t+1}}} \quad , ( t \geq 0 ) \quad .
$$
 
Hence,
$$
Res_{x_R, x_{\pi(u+1)} , \dots ,  x_{\pi(k)}}
\left [ {{\tilde A_i} \over {1- \bar x_i x_R}} \right ]
$$
is a sum of terms of the forms
$$
Res_{x_R,x_{\pi(u+1)} , \dots ,  x_{\pi(v-1)}}
\left [  {{ x_R^t \cdot 
RAT_t (x_{\pi(1)} , \dots , \hat x_R, \dots , x_{\pi(v-1)}) 
} \over
{x_R^{n+1} \bar x_R^{n+R+t+2}}} \right ] \quad ,
(t \geq 0 ) \quad .
\eqno(Gil')
$$
Applying $Res_{x_{\pi(u+1)} , \dots , x_{\pi(v-1)}}$ to
the residuand of $(Gil')$ gets rid of the dependence on
$(x_{\pi(u+1)} , \dots , x_{\pi(v-1)})$, so we have that
$$
Res_{x_R , x_{\pi(u+1)} , \dots ,  x_{\pi(k)}}
\left [ {{\tilde A_i} \over {1- \bar x_i x_R}} \right ]
$$
is a sum of terms of the form
$$
Res_{x_R}
\left [ {{ x_R^t RAT_t (x_{\pi(1)} , \dots ,x_{\pi(u-1)}) 
} \over
{x_R^{n+1} \bar x_R^{n+R+t+2}}} \right ] \quad ,
(t \geq 0 )  \quad .
\eqno(Gil'')
$$
But, each such term can be written as
$$
RAT_t (x_{\pi(1)} , \dots ,x_{\pi(u-1)})
Res_{x_R}
\left [ {{ x_R^{R+2t+1}
} \over
{x_R^{n+R+t+2} \bar x_R^{n+R+t+2}}} \right ] \quad ,
\quad (t \geq 0 )  \quad .
\eqno(Gil''')
$$
 
By crucial fact $\aleph_5'$, since $R+2t+1 \leq 2(n+R+t+1)$, we can 
perform $x_R \rightarrow \bar x_R$ in the residuand of $(Gil''')$, so
$$
Res_{x_R , x_{\pi(u+1)}, \dots ,  x_{\pi(k)}}
\left [ {{\tilde A_i} \over {1- \bar x_i x_R}} \right ]
=
Res_{x_R , x_{\pi(u+1)}, \dots ,  x_{\pi(k)}}
\left [ {{\tilde A_i} \over {1- \bar x_i x_R}} 
( x_R \rightarrow \bar x_R ) \right ] \quad.
$$
 
Applying $Res_{x_{\pi(1)} , \dots , x_{\pi(u-1)}}$ to both sides
of the above equation finishes up subcase Ib of the proof
of subsubsublemma $1.4.1.4$.
 
{\bf  Case II:} $i \in \{ \pi (1) , \dots , \pi (u-1) \}$ .
 
Let $v := \pi^{-1} (i)$, so that $i= \pi(v)$, and $v < u$.
By the definition of $u$, $\epsilon_i$ must be equal to $+1$,
so that $z_i=x_i$.
 
The left side of $(Hadas_2)$  (in the statement of $(sub)^3$lemma
$1.4.1.4$, which we are currently proving,) can be written as:
$$
Res_{x_{\pi(1)} , \dots , x_{\pi(k)}} 
[ {{\tilde A_i} \over {1- x_i x_R}} ] 
= Res_{x_{\pi(1)} , \dots , x_{\pi(v-1)}, x_i , x_{\pi(v+1)} , \dots ,
x_{\pi(u-1)},x_R, x_{\pi(u+1)} , \dots , x_{\pi(k)}} 
\left [ {{\tilde A_i} \over {1- x_i x_R}} \right ] \quad .
\eqno(Ruthie)
$$
 
The next $(sub)^4$lemma assures us that, on the right side of $(Ruthie)$,
in the order of residue taking $Res_{x_{\pi(1)} , \dots , x_{\pi(k)}}$,
it is permissible to transpose $x_i$ and $x_R$, without affecting its value,
and likewise when the residuand is replaced by
${{\tilde A_i} \over {1- x_i x_R}} (x_R \rightarrow \bar x_R )$.
 
{\bf Subsubsubsublemma 1.4.1.4.1:}
$$
Res_{x_{\pi(1)} , \dots , x_{\pi(v-1)}, x_i , x_{\pi(v+1)} , \dots ,
x_{\pi(u-1)},x_R, x_{\pi(u+1)} , \dots , x_{\pi(k)}} 
\left [ {{\tilde A_i} \over {1- x_i x_R}} \right ]
$$
$$
=Res_{x_{\pi(1)} , \dots , x_{\pi(v-1)}, x_R , x_{\pi(v+1)} , \dots ,
x_{\pi(u-1)},x_i, x_{\pi(u+1)} , \dots , x_{\pi(k)}} 
\left [ {{\tilde A_i} \over {1- x_i x_R}} \right ] \quad ,
\eqno(Noam)
$$
and
$$
Res_{x_{\pi(1)} , \dots , x_{\pi(v-1)}, x_i , x_{\pi(v+1)} , \dots ,
x_{\pi(u-1)},x_R, x_{\pi(u+1)} , \dots , x_{\pi(k)}} 
\left [ {{\tilde A_i} \over {1- x_i x_R}} (x_R \rightarrow \bar x_R ) \right ]
$$
$$
=Res_{x_{\pi(1)} , \dots , x_{\pi(v-1)}, x_R , x_{\pi(v+1)} , \dots ,
x_{\pi(u-1)},x_i, x_{\pi(u+1)} , \dots , x_{\pi(k)}} 
\left [ {{\tilde A_i} \over {1- x_i x_R}} (x_R \rightarrow \bar x_R )
\right ] \quad .
\eqno(\overline {Noam})
$$
 
We will only prove $(Noam)$, since the proof of $(\overline {Noam})$ goes
verbatim. Changing the order of residue-taking from
$Res_{x_{\pi (1)} , \dots , x_i , \dots , x_R, \dots , x_{\pi(k)}}$ to
$Res_{x_{\pi (1)} , \dots , x_R , \dots , x_i , \dots , x_{\pi(k)} }$ is
permitted by crucial fact $\aleph_4$,
{\it provided}
the residuand of $Res_{x_{\pi (1)} , \dots , x_i , \dots , x_R}$
on the left side of $(Noam)$, namely,
$$
Res_{x_{\pi(u+1)} , \dots , x_{\pi(k)}}
\left [ {{\tilde A_i} \over {1- x_i x_R}} \right ] \quad ,
$$
possesses a {\it bona fide} Laurent expansion in its variables
$(x_{\pi(1)} ,   \dots , x_i , \dots , x_{\pi(u-1)},x_R )$, i.e.
that it is a rational function whose denominator is a monomial
times a polynomial with a non-zero constant term.
 
This fact would follow, by invoking the case $w=u$, from the
following $(sub)^5$lemma.
 
{\bf Subsubsubsubsublemma 1.4.1.4.1.1 :} Let $u \leq w \leq k$, 
and $\pi^{-1}(i)<u$, then
$$
Res_{x_{\pi(w+1)} , \dots , x_{\pi(k)}}
\left [ {{\tilde A_i} \over {1- x_i x_R}} \right ] 
$$
can be expressed as a sum of terms of the form
$$
{{POL(x_{\pi(1)} ,  \dots , \hat x_R , \dots , x_{\pi(w)}) } \over
{1- x_i x_R}}
\cdot \prod_{j=1}^w {{1} \over {z_{\pi(j)}^{\alpha_j} \bar 
z_{\pi(j)}^{\beta_j} }} \cdot
\prod_{ {{r,s \in \{\pi(1) , \dots , \pi(w) \} } \atop {r<s ; r,s \neq R,i} }}
{{1} \over {(1- z_r z_s )(1- \bar z_r z_s ) }} \quad .
$$ 
 
{\bf Proof:} The proof is by decreasing induction on $w$.
When $w=k$ then this follows tautologically from $(\tilde Ai1)$ and
$(\tilde Ai2)$. Suppose the statement of the $(sub)^5$lemma is true
when $w$ is replaced by $w+1$, and let's deduce from it that it is true
for $w$ proper.
 
By the inductive hypothesis, the quantity of interest is a sum of terms of
the form
$$
Res_{x_{\pi(w+1)}}
\left [
{{POL(x_{\pi(1)} ,  \dots , \hat x_R , \dots , x_{\pi(w)}, x_{\pi(w+1)} ) }
\over {1- x_i x_R }}
\cdot \prod_{j=1}^{w+1} {{1} \over {z_{\pi(j)}^{\alpha_j} \bar 
z_{\pi(j)}^{\beta_j} }} \cdot
\prod_{ {{r,s \in \{\pi(1) , \dots , \pi(w+1)\} } 
\atop { r<s ; r,s \neq R,i} }}
{{1} \over {(1- z_r z_s )(1- \bar z_r z_s ) }} 
\right ] \quad .
$$ 
 
Let's abbreviate $\pi(w+1)$ to $S$, $\alpha_{w+1}$ to $\alpha$
and $\beta_{w+1}$ to $\beta$. Then the above can be rewritten
as:
$$
{{1} \over {1- x_i x_R}}
 \prod_{j=1}^{w} {{1} \over {z_{\pi(j)}^{\alpha_j} \bar 
z_{\pi(j)}^{\beta_j} }} \cdot
\prod_{ { {r,s \in \{\pi(1) , \dots , \pi(w)\} } \atop {r<s ; r,s \neq R,i} }}
{{1} \over {(1- z_r z_s )(1- \bar z_r z_s ) }} \cdot
Res_{x_{S}} [
$$
$$
{{POL(x_{\pi(1)} ,  \dots , \hat x_R , \dots , x_{\pi(w)}, x_S ) } 
\over { z_S^{\alpha} \bar z_S^{\beta}} } 
\prod_{ {{r \in \{\pi(1) , \dots , \pi(w) \} } \atop { r \neq S,R,i} }}
{{1} \over {1- z_r z_S }}
\prod_{ {{r \in \{\pi(1) , \dots , \pi(w) \} } \atop { r<S ; r \neq R,i} }}
{{1} \over {1- \bar z_r z_S  }}
\prod_{ {{r \in \{\pi(1) , \dots , \pi(w) \} } \atop { r>S ; r \neq R,i} }}
{{1} \over {1- \bar z_S z_r  }} ].
$$
 
$Res_{x_S} [ \dots ]$ above can be rewritten as:
$$
Coeff_{x_{S}^{(\alpha \,\,\, or \,\,\, \beta )-1}} 
$$
$$
\left [
{{POL(x_{\pi(1)} ,  \dots , \hat x_R , \dots , x_{\pi(w)}, x_S ) } 
\over {  \bar x_S^{(\beta \,\,\, or  \,\,\, \alpha ) }} } 
\prod_{ {{r \in \{\pi(1) , \dots , \pi(w) \} } \atop { r \neq S,R,i} }}
{{1} \over {1- z_r z_S }}
\prod_{ {{r \in \{\pi(1) , \dots , \pi(w) \} } \atop { r<S ; r \neq R,i} }}
{{1} \over {1- \bar z_r z_S  }}
\prod_{ {{r \in \{\pi(1) , \dots , \pi(w) \} } \atop { r>S ; r \neq R,i} }}
{{1} \over {1- \bar z_S z_r  }} \right ] .
\eqno(Danny)
$$
 
Expanding everything in powers of $x_S$, using
$$
POL(x_{\pi(1)} ,  \dots , \hat x_R , \dots , x_{\pi(w)}, x_S )=
\sum_{t=0}^{finite} POL_t (x_{\pi(1)},  
\dots , \hat x_R , \dots , x_{\pi(w)}) x_S ^t \quad ,
$$
$$
{{1} \over {1- x_S z_r } } =
\sum_{t=0}^{\infty} z_r^t x_S^t  \quad , \quad 
{{1} \over {1- x_S \bar z_r } } =
\sum_{t=0}^{\infty} \bar z_r^t x_S^t  \quad ,
$$
and/or
$$
{{1} \over {1- \bar x_S z_r } } ={{1} \over {\bar z_r + z_r x_S} }=
\sum_{t=0}^{\infty} {{(-1)^t z_r^t} \over {\bar z_r^{t+1}} } x_S^t  
 \quad,\quad
{{1} \over {1- \bar x_S \bar z_r } } ={{1} \over { z_r + \bar z_r x_S} }=
\sum_{t=0}^{\infty} {{(-1)^t \bar z_r^t} \over { z_r^{t+1}} } x_S^t  \quad ,
$$
for $r \in \{\pi(1) , \dots , \pi(w)\}\,\,\backslash\,\{i,R\}$, 
as well as,
$$
{{1} \over {\bar x_S}^{number}}= \sum_{t=0}^{\infty} NUMBER_t \cdot 
 x_S^t \quad
$$ 
multiplying out, 
and then taking the appropriate coefficient of $x_S$ in $(Danny)$, 
shows that indeed the
statement of the $(sub)^5$lemma also holds with $w$, except that the
previous sum (in the case $w+1$) has been replaced by a much larger,
but still finite, sum. This completes the proof of
$(sub)^5$lemma  $1.4.1.4.1.1$. \halmos
 
Going back to the proof of $(sub)^4$lemma $1.4.1.4.1$, we employ
the above $(sub)^5$lemma with its lowest case, $w=u$.
By the definition of $u$, $\epsilon_{\pi(1)} = \dots , \epsilon_{\pi(u-1)}= +1$,
so that 
$z_{\pi(1)} =x_{\pi(1)}, \dots , z_{\pi(u-1)}= x_{\pi(u-1)}$, and we see that
$$
Res_{ x_{\pi(u+1)} , \dots , x_{\pi(k)}} 
\left [ {{\tilde A_i} \over {1- x_i x_R}} \right ]
$$
can be expressed as a sum of terms of the form
$$
{{POL(x_{\pi(1)} , \dots , x_{\pi(u-1)}) } 
\over {1-x_R x_i}}
\cdot \prod_{j=1}^u {{1} \over {x_{\pi(j)}^{\alpha_j} \bar 
x_{\pi(j)}^{\beta_j} }} \cdot
\prod_{ {{r,s \in \{\pi(1) , \dots , \pi(u-1)\} } 
\atop { r<s ; r,s \neq i} }}
{{1} \over {(1- x_r x_s )(1- \bar x_r x_s ) }} 
 \quad ,
$$ 
each of which possesses a genuine Laurent series in its arguments
$(x_{\pi(1)} , \dots , x_{\pi(u)} )$, so for each of these,
by $\aleph'_4$,
applying $Res_{ \dots , x_i , \dots , x_R}$ to them is the same
as applying $Res_{ \dots , x_R , \dots , x_i}$, and this
finishes up the proof of $(sub)^4$lemma $1.4.1.4.1$. \halmos.
 
Going back to the proof of Case II of $(sub)^3$lemma 1.4.1.4 , we use 
$(sub)^4$lemma $1.4.1.4.1$ to change the
right side of $(Ruthie)$ to:
$$
Res_{x_{\pi(1)} , \dots , x_{\pi(v-1)}, x_R , x_{\pi(v+1)} , \dots ,
x_{\pi(u-1)},x_i, x_{\pi(u+1)} , \dots , x_{\pi(k)}} 
\left [ {{\tilde A_i} \over {1- x_i x_R}} \right ] \quad .
\eqno(Ruthie')
$$
 
But this is identical to the form treated by subcase Ia, which
says that this is the same as:
$$
Res_{x_{\pi(1)} , \dots , x_{\pi(v-1)}, x_R , x_{\pi(v+1)} , \dots ,
x_{\pi(u-1)},x_i, x_{\pi(u+1)} , \dots , x_{\pi(k)}} 
\left [ {{\tilde A_i} \over {1- x_i x_R}} (x_R \rightarrow \bar x_R )
\right ] \quad .
\eqno(Ruthie'')
$$
 
Now, we use the second part of $(sub)^4$lemma $1.4.1.4.1$, i.e.
$(\overline {Noam})$ (whose proof was omitted since it was identical to
that of $(Noam)$), to change this into
$$
Res_{x_{\pi(1)} , \dots , x_{\pi(v-1)}, x_i , x_{\pi(v+1)} , \dots ,
x_{\pi(u-1)},x_R, x_{\pi(u+1)} , \dots , x_{\pi(k)}} 
\left [ {{\tilde A_i} \over {1- x_i x_R}} (x_R \rightarrow \bar x_R )
\right ] \quad ,
$$
which was exactly what we had to show in order to complete the proof
of Case II of $(sub)^3$lemma $1.4.1.4$. This completes the proof
of $(sub)^3$lemma $1.4.1.4$. \halmos.
 
{\bf Subsubsublemma 1.4.1.5:}
Let $\tilde B_i$, ($1 \leq i <R$), be as in the statement of $1.4.1.2$, i.e.
given by $(\tilde Bi1)$. Then
$$
Res_{x_{\pi(1)} , \dots , x_{\pi(k)}} 
\left [ {{\tilde B_i} \over {1-  \bar z_i x_R}} \right ] =
Res_{x_{\pi(1)} , \dots , x_{\pi(k)}} 
\left [ {{\tilde B_i} \over {1-  \bar z_i x_R}} ( x_R \rightarrow \bar x_R )
\right ]
\quad.
\eqno(Hadas_3)
$$
{\eightrm  [ Type `{\eighttt S1415(k,n):}' in ROBBINS, for specific k and n.]}
 
{\bf Proof:}
 
{\bf  Case I:} $i \in \{ \pi (u+1) , \dots , \pi (k) \}$ .
 
Let $v := \pi^{-1} (i)$, so that $i= \pi(v)$, and $v > u$.
We have that:
$$
{{\tilde B_i} \over {1- \bar z_i x_R}} =
{
 {
  POLRAT( x_i ; x_{\pi(1)} , \dots , \hat x_R ,
    \dots , x_{\pi(v-1)} , 
  x_{\pi(v+1)} , \dots , x_{\pi(k)} )
  } \over
  {x_R^{n+1} \bar x_R^{n+R+1} z_i^{n+1} \bar z_i^{n+i+k+1} (1- \bar z_i x_R) }}
\quad .
$$
Applying $Res_{x_{\pi(v+1)} , \dots , x_{\pi(k)}}$ to it gives:
$$
Res_{x_{\pi(v+1)} , \dots , x_{\pi(k)}} \left [
{{\tilde B_i} \over {1- \bar z_i x_R}} \right ]
={{
POLRAT(x_i ;x_{\pi(1)} , \dots , \hat x_R,
 \dots , x_{\pi(v-1)}) } \over
   {x_R^{n+1} \bar x_R^{n+R+1} z_i^{n+1} \bar z_i^{n+i+k+1} (1- \bar z_i x_R) }}
  \quad .
$$
Applying $Res_{x_i}$ to the above gives:
$$
Res_{x_i , x_{\pi(v+1)} , \dots , x_{\pi(k)}}
\left [ {{\tilde B_i} \over {1- \bar z_i x_R}} \right ]
=Res_{x_i} \left [ 
{{
POLRAT(x_i ; x_{\pi(1)} , \dots , \hat x_R,
 \dots , x_{\pi(v-1)}) } \over
   {x_R^{n+1} \bar x_R^{n+R+1} z_i^{n+1} \bar z_i^{n+i+k+1} (1- \bar z_i x_R) }}
 \right ]
$$
$$
={{1} \over {x_R^{n+1} \bar x_R^{n+R+1}}} \cdot
Res_{x_i} \left [ 
{{
POLRAT(x_i; x_{\pi(1)} , \dots , \hat x_R,
 \dots , x_{\pi(v-1)} )} \over
   { z_i^{n+1} \bar z_i^{n+i+k+1} (1- \bar z_i x_R) }}
 \right ] \quad .
\eqno(Doron)
$$
 
We now distinguish two subcases:
 
{\bf Subcase Ia:} $\epsilon_i = -1$, that is: $\bar z_i = x_i$
(i.e. $z_i=\bar x_i$).
 
In this subcase, the right side of $(Doron)$ equals:
$$
 {{1 } \over
{x_R^{n+1} \bar x_R^{n+R+1}}}
\cdot
Coeff_{x_i^{n+i+k}} \left [ {{
 POLRAT(x_i ;x_{\pi(1)} , \dots ,
 \hat x_R  , \dots , x_{\pi(v-1)} )
 } \over
   { \bar x_i^{n+1} (1- x_i x_R) }}
 \right ]
  \quad .
\eqno(Gil)
$$
 
Expanding  everything involving $x_i$ as a power series in $x_i$:
$$
POLRAT(x_i;x_{\pi(1)} , \dots , \hat x_R , \dots , x_{\pi(v-1)}) =
\sum_{t=0}^{degree}
RAT_t (x_{\pi(1)} , \dots , \hat x_R , \dots , x_{\pi(v-1)} ) x_i^t \quad ,
$$
$$
{{1} \over {\bar x_i^{n+1} }} = \sum_{t=0}^{\infty} NUMBER_t \cdot x_i^t
\quad ,
$$
and, most importantly:
$$
(1- x_i x_R)^{-1} = \sum_{t=0}^{\infty} x_R^t x_i^t
\quad ,
$$
multiplying out, and collecting the coefficient of $x_i^{n+i+k}$, gives that
$Coeff_{x_i^{n+i+k}}$ in $(Gil)$ is a certain {\it polynomial }
in $x_R$, of degree $\leq n+i+k$, with coefficients that are rational
functions of $(x_{\pi(1)} , \dots , \hat x_R , \dots , x_{\pi(v-1)})$.
In other words it is a $POLRAT(x_R; x_{\pi(1)} , \dots , x_{\pi(v-1)})$
with $deg_{x_R} \leq n+i+k$. Hence,
$$
Res_{x_R, x_{\pi(u+1)} , \dots ,  x_{\pi(k)}}
\left [ {{\tilde B_i} \over {1- x_i x_R}} \right ]
=
Res_{x_R,x_{\pi(u+1)} , \dots ,  x_{\pi(v-1)}}
\left [ 
{{POLRAT(x_R ; x_{\pi(1)} , \dots , \hat x_R, \dots , x_{\pi(v-1)}) 
} \over
{x_R^{n+1} \bar x_R^{n+R+1}}} \right ] \quad ,
\eqno(Gil')
$$
where the $POLRAT$  in $(Gil')$
has degree $\leq n+i+k$ in $x_R$.
Applying $Res_{x_{\pi(u+1)} , \dots , x_{\pi(v-1)}}$ to
the residuand of $(Gil')$ gets rid of the dependence on
$(x_{\pi(u+1)} , \dots , x_{\pi(v-1)})$, so we have:
$$
Res_{x_R , x_{\pi(u+1)} , \dots ,  x_{\pi(k)}}
\left [ {{\tilde B_i} \over {1- x_i x_R}} \right ]
=
Res_{x_R}
\left [ {{POLRAT(x_R ; x_{\pi(1)} , \dots , x_{\pi(u-1)}) 
} \over
{x_R^{n+1} \bar x_R^{n+R+1}}} \right ] \quad ,
\eqno(Gil'')
$$
where $POLRAT(x_R; x_{\pi(1)}, \dots , x_{\pi(u-1)} )$ 
has degree $\leq n+i+k$ in $x_R$.
 
But, the right of $(Gil'')$ can be rewritten as:
$$
Res_{x_R}
\left [ {{
x_R^R \cdot POLRAT(x_R; x_{\pi(1)} , \dots , x_{\pi(u-1)}) 
} \over
{x_R^{n+R+1} \bar x_R^{n+R+1}}} \right ] \quad .
\eqno(Gil''')
$$
 
By crucial fact $\aleph_5'$, since $R+(n+i+k) \leq 2(n+R)$
(recall that always $n \geq k$ and $i < R$), we can 
perform $x_R \rightarrow \bar x_R$ in the residuand of $(Gil''')$, so
$$
Res_{x_R , x_{\pi(u+1)}, \dots ,  x_{\pi(k)}}
\left [ {{\tilde B_i} \over {1- x_i x_R}} \right ]
=
Res_{x_R , x_{\pi(u+1)}, \dots ,  x_{\pi(k)}}
\left [ {{\tilde B_i} \over {1- x_i x_R}} 
( x_R \rightarrow \bar x_R ) \right ] \quad.
$$
Applying $Res_{x_{\pi(1)} , \dots , x_{\pi(u-1)}}$ to both sides
of the above equation finishes up subcase Ia of the proof
of subsublemma $1.4.1.5$.
 
{\bf Subcase Ib:} $\epsilon_i = 1$, that is: $\bar z_i = \bar x_i$
(i.e. $z_i=x_i$).
 
In this subcase, the right side of $(Doron)$ equals:
$$
 {{1 } \over
{x_R^{n+1} \bar x_R^{n+R+1}}}
\cdot
Coeff_{x_i^{n}} \left [ {{
 POLRAT(x_i ; x_{\pi(1)} , \dots ,
 \hat x_R  , \dots , x_{\pi(v-1)})
 } \over
   { \bar x_i^{n+i+k+1} (1- \bar x_i x_R) }}
 \right ]
  \quad .
\eqno(Gil)
$$
Expanding  everything involving $x_i$ as a power series in $x_i$:
$$
POLRAT(x_i ; x_{\pi(1)} , \dots , \hat x_R , \dots , x_{\pi(v-1)}) =
\sum_{t=0}^{degree}
RAT_t (x_{\pi(1)} , \dots , \hat x_R , \dots , x_{\pi(v-1)} ) x_i^t \quad ,
$$
$$
{{1} \over {\bar x_i^{n+i+k+1} }} = \sum_{t=0}^{\infty} NUMBER_t \cdot x_i^t
\quad ,
$$
and, most importantly:
$$
(1- \bar x_i x_R)^{-1} = (1- x_R + x_i x_R )^{-1} =
( \bar x_R + x_i x_R )^{-1} =
\sum_{t=0}^{\infty} (-1)^t {{x_R^t} \over {\bar x_R^{t+1}}} x_i^t \quad ,
$$
multiplying out, and collecting the coefficient of $x_i^n$, gives that
$Coeff_{x_i^n}$ in $(Gil)$ is a certain sum of terms of the form
$$
RAT_t (x_{\pi(1)} , \dots , \hat x_R, \dots , x_{\pi(v-1)}) 
{{x_R^t} \over {\bar x_R^{t+1}}} \quad , ( t \geq 0 ) \quad .
$$
 
Hence,
$$
Res_{x_R, x_{\pi(u+1)} , \dots ,  x_{\pi(k)}}
\left [ {{\tilde B_i} \over {1- \bar x_i x_R}} \right ]
$$
is a sum of terms of the forms
$$
Res_{x_R,x_{\pi(u+1)} , \dots ,  x_{\pi(v-1)}}
\left [  {{ x_R^t \cdot 
RAT_t (x_{\pi(1)} , \dots , \hat x_R, \dots , x_{\pi(v-1)}) 
} \over
{x_R^{n+1} \bar x_R^{n+R+t+2}}} \right ] \quad ,
(t \geq 0 ) \quad .
\eqno(Gil')
$$
 
Applying $Res_{x_{\pi(u+1)} , \dots , x_{\pi(v-1)}}$ to
the residuand of $(Gil')$ gets rid of the dependence on
$(x_{\pi(u+1)} , \dots , x_{\pi(v-1)})$, so we have that
$$
Res_{x_R , x_{\pi(u+1)} , \dots ,  x_{\pi(k)}}
\left [ {{\tilde B_i} \over {1- \bar x_i x_R}} \right ]
$$
is a sum of terms of the form
$$
Res_{x_R}
\left [ {{ x_R^t RAT_t (x_{\pi(1)} , \dots ,x_{\pi(u-1)}) 
} \over
{x_R^{n+1} \bar x_R^{n+R+t+2}}} \right ] \quad ,
(t \geq 0 )  \quad .
\eqno(Gil'')
$$
 
But, each such term can be written as
$$
RAT_t (x_{\pi(1)} , \dots ,x_{\pi(u-1)}) \cdot
Res_{x_R}
\left [ {{ x_R^{R+2t+1}
} \over
{x_R^{n+R+t+2} \bar x_R^{n+R+t+2}}} \right ] \quad ,
\quad (t \geq 0 )  \quad .
\eqno(Gil''')
$$
 
By crucial fact $\aleph_5'$, since $R+2t+1 \leq 2(n+R+t+1)$, we can 
perform $x_R \rightarrow \bar x_R$ in the residuand of $(Gil''')$, so
$$
Res_{x_R , x_{\pi(u+1)}, \dots ,  x_{\pi(k)}}
\left [ {{\tilde B_i} \over {1- \bar x_i x_R}} \right ]
=
Res_{x_R , x_{\pi(u+1)}, \dots ,  x_{\pi(k)}}
\left [ {{\tilde B_i} \over {1- \bar x_i x_R}} 
( x_R \rightarrow \bar x_R ) \right ] \quad.
$$
 
Applying $Res_{x_{\pi(1)} , \dots , x_{\pi(u-1)}}$ to both sides
of the above equation finishes up subcase Ib of the proof
of subsubsublemma $1.4.1.5$.
 
{\bf  Case II:} $i \in \{ \pi (1) , \dots , \pi (u-1) \}$ .
 
Let $v := \pi^{-1} (i)$, so that $i= \pi(v)$, and $v < u$.
By the definition of $u$, $\epsilon_i$ must be equal to $+1$,
so that $z_i=x_i$.
 
The left side of $(Hadas_3)$ can be written as:
$$
Res_{x_{\pi(1)} , \dots , x_{\pi(k)}} 
[ {{\tilde B_i} \over {1- \bar x_i x_R}} ] 
= Res_{x_{\pi(1)} , \dots , x_{\pi(v-1)}, x_i , x_{\pi(v+1)} , \dots ,
x_{\pi(u-1)},x_R, x_{\pi(u+1)} , \dots , x_{\pi(k)}} 
\left [ {{\tilde B_i} \over {1- \bar x_i x_R}} \right ] \quad .
\eqno(Ruthie)
$$
 
The first step is to `unbar' the $x_i$ in $(Ruthie)$.
 
{\bf Subsubsubsublemma 1.4.1.5.1:}
Let $\tilde B_i$, ($1 \leq i < R$), be as in the statement of $1.4.1.2$, i.e.
given by $(\tilde Bi1)$. (Recall that $\pi$ is an arbitrary permutation,
$\pi(u)=R,\pi(v)=i$ and $v<u$.)
Then
$$
Res_{x_{\pi(1)} , \dots , x_{\pi(k)}} 
\left [ {{\tilde B_i ( x_i \rightarrow \bar x_i ) } \over 
{1-  x_i x_R}} \right ] =
Res_{x_{\pi(1)} , \dots , x_{\pi(k)}} 
\left [ {{\tilde B_i} \over {1- \bar x_i x_R}} 
\right ]
\quad.
\eqno(Michal)
$$
$$
Res_{x_{\pi(1)} , \dots , x_{\pi(k)}} 
\left [ {{\tilde B_i ( x_i \rightarrow \bar x_i , x_R \rightarrow \bar x_R)} 
\over 
{1-  x_i \bar x_R}} \right ] =
Res_{x_{\pi(1)} , \dots , x_{\pi(k)}} 
\left [ {{\tilde B_i (x_R \rightarrow \bar x_R )}
\over {1-  \bar x_i \bar x_R}}
\right ]
\quad.
\eqno(Michal')
$$
 
{\bf Proof:} The proof of $(Michal)$ and $(Michal')$ are very
similar to subcases Ia and Ib respectively, except that the
roles of $i$ and $R$ are reversed. We will only prove
$(Michal)$, since $(Michal')$ is almost identical to the
proof of subcase Ib (with $i$ and $R$ interchanged.)
 
We have that:
$$
{{\tilde B_i (x_i \rightarrow \bar x_i )} \over {1-  x_i x_R}} =
{
 {
  POLRAT( x_i ; x_{\pi(1)} , \dots , \hat x_R ,
    \dots , x_{\pi(v-1)} , 
  x_{\pi(v+1)} , \dots , x_{\pi(k)} )
  } \over
  {x_R^{n+1} \bar x_R^{n+R+1} \bar x_i^{n+1} x_i^{n+i+k+1} (1- x_i x_R) }}
\quad .
$$
 
Applying $Res_{x_{\pi(u+1)} , \dots , x_{\pi(k)}}$ to it gives:
$$
Res_{x_{\pi(u+1)} , \dots , x_{\pi(k)}} 
\left [ {{\tilde B_i ( x_i \rightarrow \bar x_i )} 
\over 
{1-  x_i x_R}} \right ]
={{
POLRAT(x_i ;x_{\pi(1)} , \dots , \hat x_i,
 \dots , x_{\pi(u-1)}) } \over
   {x_R^{n+1} \bar x_R^{n+R+1} \bar x_i^{n+1} x_i^{n+i+k+1} 
(1- x_i x_R) }}
  \quad .
$$
 
Applying $Res_{x_R}$ to the above gives:
$$
Res_{x_R , x_{\pi(u+1)} , \dots , x_{\pi(k)}}
\left [ {{\tilde B_i ( x_i \rightarrow \bar x_i ) } 
\over 
{1-  x_i  x_R}} \right ]
=Res_{x_R} \left [ 
{{
POLRAT(x_i ; x_{\pi(1)} , \dots , \hat x_i,
 \dots , x_{\pi(u-1)}) } \over
   {x_R^{n+1} \bar x_R^{n+R+1} \bar x_i^{n+1} x_i^{n+i+k+1} (1- x_i x_R) }}
 \right ]
$$
$$
={{POLRAT(x_i; x_{\pi(1)} , \dots , \hat x_i,
 \dots , x_{\pi(u-1)} )} \over {x_i^{n+i+k+1} \bar x_i^{n+1}}} \cdot
Res_{x_R} \left [ 
{{1} \over
   { x_R^{n+1} \bar x_R^{n+R+1} (1- x_i x_R) }}
 \right ] \quad .
\eqno(Herb)
$$
 
The right side of $(Herb)$ equals:
$$
{{POLRAT(x_i; x_{\pi(1)} , \dots , \hat x_i,
 \dots , x_{\pi(u-1)} )} \over {x_i^{n+i+k+1} \bar x_i^{n+1}}} 
\cdot
Coeff_{x_R^{n}} \left [ {{1}
 \over
   { \bar x_R^{n+R+1} (1- x_i x_R) }}
 \right ]
  \quad .
\eqno(David)
$$
 
Expanding  the two quantities involving $x_R$ 
above as a power series in $x_R$:
$$
{{1} \over {\bar x_R^{n+R+1} }} = \sum_{t=0}^{\infty} NUMBER_t \cdot x_R^t
\quad ,
$$
and, more importantly:
$$
(1- x_i x_R)^{-1} = \sum_{t=0}^{\infty} x_i^t x_R^t
\quad ,
$$
multiplying out, and collecting the coefficient of $x_R^{n}$, gives that
$Coeff_{x_R^n}$ in $(David)$ is a certain {\it polynomial }
in $x_i$, of degree $\leq n$.
Combined with the $POLRAT$ in the front of
$(David)$ (whose degree in $x_i$ is $\leq 2k+1$, by the last
sentence of the statement of $(sub)^3$lemma $1.4.1.2$)
we have a $POLRAT$ in $x_i$ of degree $\leq n+2k+1$ in $x_i$.
 
Hence,
$$
Res_{x_i, x_{\pi(v+1)} , \dots ,  x_{\pi(k)}}
\left [ {{\tilde B_i (x_i \rightarrow \bar x_i )} \over {1- x_i x_R}} \right ]
=
Res_{x_i,x_{\pi(v+1)} , \dots ,  x_{\pi(u-1)}}
\left [ 
{{POLRAT(x_i ; x_{\pi(1)} , \dots , \hat x_i, \dots , x_{\pi(u-1)}) 
} \over
{x_i^{n+i+k+1} \bar x_i^{n+1}}} \right ] \quad ,
\eqno(David')
$$
 
where the $POLRAT$  in $(David')$
has degree $\leq n+2k+1$ in $x_i$.
Applying $Res_{x_{\pi(v+1)} , \dots , x_{\pi(u-1)}}$ to
the residuand of the right side of $(David')$ gets rid of the dependence on
$(x_{\pi(v+1)} , \dots , x_{\pi(u-1)})$, so we have:
$$
Res_{x_i , x_{\pi(v+1)} , \dots ,  x_{\pi(k)}}
\left [ {{\tilde B_i (x_i \rightarrow \bar x_i )} \over {1- x_i x_R}} \right ]
=
Res_{x_i}
\left [ {{POLRAT(x_i ; x_{\pi(1)} , \dots , x_{\pi(v-1)}) 
} \over
{x_i^{n+i+k+1} \bar x_i^{n+1}}} \right ] \quad ,
\eqno(David'')
$$
where $POLRAT(x_i; x_{\pi(1)}, \dots , x_{\pi(v-1)} )$ 
has degree $\leq n+2k+1$ in $x_i$.
 
But, the right of $(David'')$ can be rewritten as:
$$
Res_{x_i}
\left [ {{
\bar x_i^{i+k} \cdot POLRAT(x_i; x_{\pi(1)} , \dots , x_{\pi(v-1)}) 
} \over
{x_i^{n+i+k+1} \bar x_i^{n+i+k+1}}} \right ] \quad .
\eqno(David''')
$$
 
By crucial fact $\aleph_5'$, since $(i+k)+n+(2k+1) \leq 2(n+i+k)$
(recall that always $n \geq k$ and $i \geq 1$), we can 
perform $x_i \rightarrow \bar x_i$ in the residuand of $(David''')$, so
$$
Res_{x_i , x_{\pi(v+1)}, \dots ,  x_{\pi(k)}}
\left [ {{\tilde B_i (x_i \rightarrow \bar x_i )} \over {1- x_i x_R}} \right ]
=
Res_{x_i , x_{\pi(v+1)}, \dots ,  x_{\pi(k)}}
\left [ {{\tilde B_i } \over {1- \bar x_i x_R}}  \right ] \quad.
$$
 
Applying $Res_{x_{\pi(1)} , \dots , x_{\pi(v-1)}}$ to both sides
of the above equation finishes up the proof of $(Michal)$.
As we said, the proof of $(Michal')$ is left to the reader.
This finishes the proof of $(sub)^4$lemma $1.4.1.5.1$. \halmos
 
Going back to the proof of case II of $(sub)^3$lemma $1.4.1.5$,
we have to prove that the right sides of $(Michal)$ and $(Michal')$
are equal to each other. By the above $(sub)^4$lemma,
this would follow from their left sides
being equal to each other. But the proof of this is exactly
as that of case II of $(sub)^3$lemma $1.4.1.4$: one proves
the analog of $(Noam)$ and $(\overline{Noam})$, which reduce
the statements to that of subcase Ia. This completes the proof of
$(sub)^3$lemma $1.4.1.5$. \halmos
 
We have one more $(sub)^3$lemma to go before we can complete the
proof of $(sub)^2$lemma $1.4.1$, which is to show that
the last sum on the right side of $(Celia)$ (in the statement of $(sub)^3$lemma
$1.4.1.2$, also has property $(Hadas)$.)
 
{\bf Subsubsublemma 1.4.1.6:}
Let $\tilde B_i$, ($R < i \leq k$), be as in the statement of $1.4.1.2$, i.e.
given by $(\tilde Bi2)$. Then
$$
Res_{x_{\pi(1)} , \dots , x_{\pi(k)}} 
\left [ {{\tilde B_i} \over {1-  z_i \bar x_R}} \right ] =
Res_{x_{\pi(1)} , \dots , x_{\pi(k)}} 
\left [ {{\tilde B_i} \over {1-  z_i \bar x_R}} ( x_R \rightarrow \bar x_R )
\right ]
\quad.
\eqno(Hadas_4)
$$
{\eightrm  [ Type `{\eighttt S1416(k,n):}' in ROBBINS, for specific k and n.]}
 
{\bf Proof:} The form of $\tilde B_i$, for $R < i \leq k$, given by
$(\tilde Bi2)$ is identical to that of $\tilde A_i$, for
$R < i \leq k$, given by $(\tilde Ai2)$. The only difference between
the statement of $(sub)^3$lemma $1.4.1.4$ and that of $1.4.1.6$ is that
the $x_R$ in the denominator of ${ {\tilde A_i} \over {1- z_i x_R}}$
is changed into $\bar x_R$ : ${ {\tilde B_i} \over {1- z_i \bar x_R}}$.
The proof goes verbatim. This completes the proof of $(sub)^3$lemma
$1.4.1.6$. \halmos
 
Going back to the proof of $(sub)^2$lemma $1.4.1$, we know by
$(sub)^3$lemmas $(1.4.1.3-6)$ that each part of
$(Celia)$, in the statement of $(sub)^3$lemma $1.4.1.2$, has
property $(Hadas)$, and by the additivity of this property, so does
the whole 
$F_{n,k} ( \epsilon_1(x_1) , \dots ,  x_R , \dots , \epsilon_k(x_k) )$.
But for the latter to have property $(Hadas)$ is precisely what $(sub)^2$lemma
$1.4.1$ is saying. 
This completes the proof of $(sub)^2$lemma $1.4.1$. \halmos
 
As we noted at the very beginning of the proof of sublemma $1.4$,
subsublemma $1.4.1$, read from right to left, can be used repeatedly
to `unbar' $x_i$'s, until we get that
$$
Res_{x_{\pi(1)} , \dots , x_{\pi(k)}}
[ F_{n,k} ( \epsilon_1(x_1) , \dots , \epsilon_k(x_k) ) ] = 
Res_{x_{\pi(1)} , \dots , x_{\pi(k)}}
 [  F_{n,k} (x_1 , \dots , x_k ) ]  \quad  .
\eqno(SofTov)
$$
 
Since  $  F_{n,k} (x_1 , \dots , x_k ) $ has a genuine
Laurent expansion, crucial fact $\aleph_4'$ reassures
us that we can unscramble the  order of residue-taking on the right side
of $(SofTov)$:
 
$$
Res_{x_{\pi(1)} , \dots , x_{\pi(k)}}
 [  F_{n,k} (x_1 , \dots , x_k ) ]  
=
Res_{x_1 , \dots , x_k }
 [  F_{n,k} (x_1 , \dots , x_k ) ]  
\quad  .
\eqno(HakolTov)
$$
 
Combining $(SofTov)$ with $(HakolTov)$ yields $(1.4')$, which we
noted was equivalent to the statement of sublemma 1.4.
This completes the proof of sublemma $1.4$. \halmos
 
\medskip
\centerline{\bf INTERACT, Part II} 
 
At this juncture, we have to introduce 
the rational function $T_{k,n}(x)=T_{k,n}(x_1 , \dots , x_k)$ defined by:
$$
T_{k,n} (x) := \prod_{i=1}^{k} 
{ {1} \over { ( \bar  x_i x_i )^{n+k+1}  }}
\prod_{ 1 \leq i < j \leq k}
{ {1} \over {(1- x_i x_j ) (1- \bar x_i x_j )(1- x_i \bar x_j )
             (1- \bar x_i \bar x_j )} } \qquad .
$$
 
We can rewrite Sublemmas 1.1' and 1.2' more compactly as follows.
 
{\bf Sublemma 1.1''}  The total number of $n \times k$-Magog trapezoids, 
$b_k (n)$, is given by:
 
$$
Res_{x_1 , \dots , x_k} \left [
\Delta_k (x) \, T_{k,n} (x) \,
\prod_{i=1}^k x_i^{i+1}
\prod_{1 \leq i < j \leq k }
(1- \bar x_i x_j ) (1- x_i \bar x_j ) ( 1- \bar x_i \bar x_j ) \right ]
\qquad .
\eqno(MagogTotal'')
$$
 
{\bf Proof:} \halmos .
 
{\bf Sublemma 1.2'':} The total number of $n \times k -$ Gog 
trapezoids, $m_k (n)$, is given by:
$$
Res_{x_1 , \dots , x_k}
\left [ \Phi_k ( x ) T_{k,n} (x)
 \prod_{i=1}^k  (\bar x_i)^{k-i} x_i^k \prod_{ 1 \leq i < j \leq k}
(1-x_i \bar x_j ) (1- \bar x_i \bar x_j )  \right ]
\qquad ,
\eqno(GogTotal'')
$$
where $\Phi_k (x)= \Phi_k ( x_1 , \dots , x_k)$ 
is the polynomial defined in ($Gog_1$).
 
{\bf Proof:} \halmos.
  
Using  sublemmas $1.3$ and $1.4$, 
and the fact that $T_{k,n}(x)$  is symmetric while $\Delta_k$  and $\Phi_k$ are
anti-symmetric with respect to the 
action of the group of signed permutations $W( \B_k )$,
we can average Sublemmas 1.1'' and 1.2'',
to get yet other residue expressions for the desired quantities
$b_k (n)$ and $m_k (n)$, as follows.
 
{\bf Sublemma 1.1''':} The total number of $n \times k$-Magog trapezoids, 
$b_k (n)$, is given by:
$$
{ {1} \over {2^k k!}}
Res_{x_1 , \dots , x_k} \left \{
 T_{k,n} (x) \, \Delta_k (x) \, \cdot
\sum_{g \in W( \B_k ) }
\sgn(g) \cdot
g \left [
\prod_{i=1}^k x_i^{i+1}
\prod_{1 \leq i < j \leq k }
(1- \bar x_i x_j ) (1- x_i \bar x_j ) ( 1- \bar x_i \bar x_j )
  \right ] \, \right \}
 \,     \qquad .
\eqno(MagogTotal''')
$$
 
{\bf Proof:} \halmos.
 
{\bf Sublemma 1.2''':} The total number of $n \times k$-Gog
trapezoids, $m_k (n)$, is given by ,
$$
{ {(-1)^k} \over {2^k k!} } 
Res_{x_1 , \dots , x_k} \left [
T_{k,n}(x) \Phi_k(x)^2 \right ]
\qquad ,
\eqno(GogTotal''')
$$
where $\Phi_k (x)= \Phi_k ( x_1 , \dots , x_k)$ 
is the polynomial defined by ($Gog_1$), in the statement of sublemma $1.2$.
 
{\bf Proof:} See the definition of $\Phi_k (x)$ given in ($Gog_1$) \halmos.
 
\vfill
\eject
 
%begin macros
\def \inv{\mathop{\rm inv} \nolimits}
\def \sgn{\mathop{\rm sgn} \nolimits} 
\def\hat{\widehat}
\def\tilde{\widetilde}
\baselineskip=14pt
\parskip=10pt
\def\B{{\cal B}}
\def\S{{\cal S}}
\def\halmos{\hbox{\vrule height0.15cm width0.01cm\vbox{\hrule height
 0.01cm width0.2cm \vskip0.15cm \hrule height 0.01cm width0.2cm}\vrule
 height0.15cm width 0.01cm}}
\font\eightrm=cmr8  
\font\eighttt=cmtt8
\magnification=\magstephalf

\parindent=0pt
\overfullrule=0in
\headline={\rm  \ifodd\pageno  \RightHead  \else  \LeftHead  \fi}
\def\RightHead{\centerline{Proof  of the ASM Conjecture-Act  V}}
\def\LeftHead{ \centerline{Doron Zeilberger}}
%\pageno=66
%end macros
\centerline{\bf Act V. GOG=MAGOG}
 
\indent{
{ \it
$\dots$ Gog and Magog, to gather them 
together to battle:
the number of whom is as the sand of the sea.} 
\qquad\qquad---(Revelations XX,8)}
 
We have one more sublemma to go, that asserts that not only are the
{\it iterated-residues} expressing $b_k(n)$ and $m_k(n)$, as given in
$(MagogTotal''')$ and $(GogTotal''')$ respectively, equal to each other,
but the actual {\it residue-ands} are.
 
{\bf Sublemma 1.5:} Let $\Phi_k(x)=\Phi_k(x_1, \dots , x_k)$ and 
$\Delta_k(x)=\Delta_k (x_1 , \dots , x_k)$ be defined as above, i.e. by
$$
\Phi_k ( x_1 , \dots , x_k ) := 
(-1)^k \sum_{ g \in W( \B_k ) } \sgn(g) \cdot
g \left [ \prod_{i=1}^{k} \bar x_i^{k-i} x_i^{k} 
\prod_{1 \leq i < j \leq k} (1- x_i \bar x_j )(1- \bar x_i \bar x_j ) 
\right ]
\qquad,
\eqno(Gog_1)
$$
$$
\Delta_k ( x_1 , \dots , x_k ) := 
\prod_{i=1}^{k} ( 1- 2 x_i ) \prod_{ 1 \leq i < j \leq k}
( x_j - x_i ) ( x_j + x_i -1) \qquad .
\eqno(Delta)
$$
 
The following identity holds:
$$
\Delta_k (x_1 , \dots , x_k ) 
\sum_{g \in W( \B_k ) }
\sgn(g) \cdot
g  \left [
\prod_{i=1}^k x_i^{i+1}
\prod_{1 \leq i < j \leq k }
(1- \overline x_i x_j ) (1- x_i \overline x_j ) ( 1- \overline x_i \overline x_j )
\right  ] \, = \, (-1)^k \,
\Phi_k ( x_1 , \dots , x_k )^2 \qquad .
\eqno(Gog=Magog)
$$
{\eightrm  [ Type `{\eighttt S15(k):}' in ROBBINS, for specific k.]}
 
{\bf Proof:} Since $\Phi_k (x)$ is $W( \B_k )$-antisymmetric, it can be written
as $\Phi_k (x)= \Delta_k (x) \Omega_k (x)$, for some polynomial 
$\Omega_k (x) = \Omega_k (x _1 , \dots , x_k )$ 
that is {\it symmetric} with respect
to $W( \B_k )$. This means that $\Omega_k (x)$ is a symmetric polynomial in
$( x_1 \overline x_1 , \dots , x_k \overline x_k )$. Since the degree of 
$\Phi_k (x)$,
viewed as a polynomial in $x_1$, is $\leq (2k-1) +2(k-1)$, and that of
$\Delta_k (x)$ is $1+ 2(k-1)$, it follows that the degree of
$\Omega_k (x)$, as a polynomial in $x_1$ is $\leq 2k-2$, and hence as
a polynomial in $x_1 \overline x_1$ is $\leq k-1$. Similarly, the left side of
(Gog=Magog) can be written as $(-1)^k \Delta_k (x)^2 L_k (x)$, for some
polynomial $L_k (x)$,
that is symmetric in $( x_1 \overline x_1 , \dots , x_k \overline x_k )$.
The degree of $L_k$, as a polynomial in $x_1$, is $\leq$
$(k+1)+3(k-1)-(1+2(k-1))=2k-1$, and hence its degree in $x_1 \overline x_1$
is $\leq k-1$. We have to prove that $L_k - \Omega_k^2$ is identically zero.
Since $L_k -\Omega_k^2$ is a polynomial of degree $\leq 2k-2$ in $x_1 \overline x_1$,
it suffices to prove that it is zero at $2k-1$ distinct values of
$x_1 \overline x_1$. This will follow by induction on $k$ from the following
two subsublemmas.

{\bf Subsublemma 1.5.1:} Let $L_k (x)$ be the left side of (Gog=Magog)
divided by $\Delta_k (x)^2$. We have
$$
L_k ( x_1 , \dots , x_k ) \vert_{x_1 = x_i^{-1} }=
\left (
 \prod_{{j=2} \atop {j \not= i} }^k
{
  {(1- \overline x_i x_j ) (1- \overline x_i \overline x_j ) }
  \over
   {x_i}
}
\right  ) ^2 L_{k-2} ( x_2 , \dots , \widehat x_i , \dots , x_k ) \quad ,
\quad (i= 2 , \dots , k )
\eqno(1.5.1a)
$$
$$
L_k ( x_1 , \dots , x_k ) \vert_{x_1 = \overline x_i^{-1} }=
\left (
 \prod_{{j=2} \atop {j \not= i} }^k
{
  {(1- x_i x_j ) (1-  x_i \overline x_j ) }
  \over
   {\overline x_i}
}
\right ) ^2 L_{k-2} ( x_2 , \dots , \widehat x_i , \dots , x_k ) \quad ,
\quad (i= 2 , \dots , k )
\eqno(1.5.1b)
$$
$$
L_k ( x_1 , \dots , x_k ) \vert_{x_1 = 0 }=
L_{k-1} ( x_2 , \dots , x_k ) \quad .
\eqno(1.5.1c)
$$
{\eightrm  [ Type `{\eighttt S151(k):}' in ROBBINS, for specific k.]}
 
{\bf Proof:} By definition, $L_k(x)$ is the left side of (Gog=Magog)
divided by $\Delta_k(x)^2$. By using the definition $(Delta)$ of $\Delta_k(x)$,
as well as its $W(\B_k)$-antisymmetry, we get,
$$
L_k (x_1 , \dots , x_k ) =
\sum_{ g \in W( \B_k ) }
g  \left [
\prod_{i=1}^k { {x_i^{i+1}} \over {( 1 - 2 x_i )} }
\prod_{ 1 \leq i < j \leq k }
{
  {(1- \overline x_i x_j ) (1- x_i \overline x_j ) (1- \overline x_i \overline x_j ) }
   \over
   { ( x_j - x_i ) ( x_j + x_i -1 ) }
}
\right ] =
$$
$$
\sum_{ \epsilon \in \{+1 , -1 \}^k }
\epsilon  \left \{
\sum_{ \pi \in \S_k }
\pi \left [
\prod_{i=1}^k { {x_i^{i+1}} \over {( 1-2 x_i )} }
\prod_{ 1 \leq i < j \leq k }
{
  {(1- \overline x_i x_j ) (1- x_i \overline x_j ) (1- \overline x_i \overline x_j ) }
   \over
   { ( x_j - x_i ) ( x_j + x_i -1 ) }
}
 \right ]  \right \} \quad .
$$
Since the argument of $\pi$ is antisymmetric w.r.t. to $\S_k$, {\it except}
for $\prod_{i=1}^k x_i^{i+1} = \prod_{i=1}^k x_i^2 \prod_{i=1}^k x_i^{i-1} $,
it follows, by Vandermonde's determinant identity (see ``Crucial Facts"), that
 
$$
L_k ( x_1 , \dots , x_k ) \,  = \,
\sum_{ \epsilon \in \{+1 , -1 \}^k }
\epsilon  \left \{
\prod_{i=1}^k { {x_i^{2}} \over {( 1 - 2 x_i )} }
\prod_{ 1 \leq i < j \leq k }
{
  {(1- \overline x_i x_j ) (1- x_i \overline x_j ) (1- \overline x_i \overline x_j ) }
   \over
   { ( x_j + x_i -1 ) }
}
\right \}  \quad  .
\eqno(1.5.1d)
$$
 
To prove (1.5.1a) for $i=2, \dots , k$, it suffices, by symmetry, to 
consider $i=2$. This would also imply, by the symmetry 
$x \leftrightarrow \overline x$, built into $W( \B_k )$,
that (1.5.1b) is true. So let's plug in $x_1 = x_2^{-1}$ in (1.5.1d).
The only way that $\epsilon [ (1- \overline x_1 x_2 ) (1- x_1 \overline x_2 ) 
(1- \overline x_1 \overline x_2 ) ]$ cannot be zero then, is for
$\epsilon_1$ and $\epsilon_2$ to be both $+1$. Hence the
$2^k$-summand sum (1.5.1d) shrinks to a $2^{k-2}$-summand sum:
$$
{
      {( x_2^{-1} )^2}  
   \over 
      {( 1 - 2 x_2^{-1} )}
 }
 { 
   {( x_2 )^2}  
     \over 
   {( 1- 2 x_2  )} 
}
{ 
  {(1- \overline {x_2^{-1}} x_2 ) (1- x_2^{-1} \overline x_2 )
    (1- \overline {x_2^{-1}} \overline x_2 ) 
  }
    \over
  {x_2 + x_2^{-1} -1}
}
\sum_{{ {\epsilon \in \{ +1,-1 \}^k }
       \atop
        {\epsilon_1 =1 , \epsilon_2 =1}
     }}
\epsilon [
\prod_{j=3}^k
 (
{ 
  {(1- \overline {x_2^{-1}} x_j ) (1- x_2^{-1} \overline x_j )
    (1- \overline {x_2^{-1}} \overline x_j ) 
  }
    \over
  {x_j + x_2^{-1} -1}
}
\cdot
$$
$$
{ 
  {(1- \overline {x_2} x_j ) (1- x_2 \overline x_j )
    (1- \overline x_2 \overline x_j ) 
  }
    \over
  {x_j + x_2 -1}
}  )
\prod_{j=3}^k
{ {x_j^2} \over {1-2 x_j }}
\prod_{ 3 \leq i < j \leq k}
{
 { (1- \overline x_i x_j ) (1- x_i \overline x_j ) (1- \overline x_i \overline x_j ) }
\over
 { ( x_j + x_i -1 ) }
}  ] \quad .
 \eqno(1.5.1e)
$$
 
We now need the following {\it utterly} trivial subsubsublemmas.
 
{\bf Subsubsublemma 1.5.1.1:}
$$
{ 
  {(1- \overline {x_2^{-1}} t ) (1- x_2^{-1} \overline t )
    (1- \overline {x_2^{-1}} \overline t ) 
  }
    \over
  {t + x_2^{-1} -1}
}
\cdot
{ 
  {(1- \overline {x_2} t ) (1- x_2 \overline t )
    (1- \overline x_2 \overline t ) 
  }
    \over
  {t + x_2 -1}
}
=
\left (
{
   {(1- \overline x_2 \overline t ) (1- \overline x_2 t ) }
    \over 
   {x_2}
}
 \right )^2 \quad.
\eqno(1.5.1.1)
$$
{\eightrm  [ Type `{\eighttt S1511():}' in ROBBINS, to get a rigorous proof.]}
 
{\bf Proof:} \halmos.
 
{\bf Subsubsublemma 1.5.1.2:}
$$
 { {( x_2^{-1} )^2}  \over {( 1- 2 x_2^{-1} )} }
 { {( x_2 )^2}  \over {( 1- 2 x_2 )} }
{ 
  {(1- \overline {x_2^{-1}} x_2 ) (1- x_2^{-1} \overline x_2 )
    (1- \overline {x_2^{-1}} \overline x_2 ) 
  }
    \over
  {x_2 + x_2^{-1} -1}
}
= 1 \qquad .
\eqno(1.5.1.2)
$$
{\eightrm  [ Type `{\eighttt S1512():}' in ROBBINS, to get a rigorous proof.]}
 
{\bf Proof:} \halmos.
 
Substituting (1.5.1.1) and (1.5.1.2) into (1.5.1e) yields:
$$
L_k ( x_2^{-1} , x_2 , x_3 , \dots , x_k ) =
\sum_{{ {\epsilon \in \{ +1,-1 \}^k }
       \atop
        {\epsilon_1 =1 , \epsilon_2 =1}
     }}
\epsilon  [
\prod_{j=3}^k
\left (
{
   {(1- \overline x_2 \overline x_j ) (1- \overline x_2 x_j ) }
    \over 
   {x_2}
}
\right )^2 
\cdot
$$
$$
\prod_{j=3}^k
{ {x_j^2} \over {2 x_j -1}}
\prod_{ 3 \leq i < j \leq k}
{
 { (1- \overline x_i x_j ) (1- x_i \overline x_j ) (1- \overline x_i \overline x_j ) }
\over
 { ( x_j + x_i -1 ) }
}
 ] \, = \,
\left (  \prod_{j=3}^k
{
   {(1- \overline x_2 \overline x_j ) (1- \overline x_2 x_j ) }
    \over 
   {x_2}
}
\right )^2 
L_{k-2} ( x_3 , \dots , x_k ) \quad . 
\eqno(1.5.1f)
$$
This completes the proof of (1.5.1a) and (1.5.1b).  To prove
(1.5.1c), note that when $x_1 =0$, the only non-zero summands
of (1.5.1d) are those for which $\epsilon_1 = -1$. The $2^k$-summand
sum (1.5.1d) collapses then to a $2^{k-1}$-summand sum as follows:
$$
\sum_{{ {\epsilon \in \{ +1,-1 \}^k }
       \atop
        {\epsilon_1 = -1}
     }}
\epsilon
\left (
\prod_{i=2}^k 
 { {x_i^2} \over {1- 2 x_i } }
\prod_{j=2}^k {{(1- \overline x_j )} \over {x_j} }
\prod_{ 2 \leq i < j \leq k}
{
{ (1- \overline x_i x_j ) (1- x_i \overline x_j ) ( 1- \overline x_i 
\overline x_j ) }
\over
 { (x_j + x_i -1) }
}
\right )
\, = \,
L_{k-1} ( x_2 , \dots , x_k ) \quad.
$$
This completes the proof of subsublemma 1.5.1 \halmos .
 
{\bf Subsublemma 1.5.2:} $\Omega_k (x):= \Phi_k (x) / \Delta_ k (x)$ (see
$(Gog_1 )$ and (Delta), in the statement of sublemma 1.5
for definitions), satisfies
$$
\Omega_k ( x_1 , \dots , x_k ) \vert_{x_1 = x_i^{-1} }\, = \,
\left ( \prod_{{j=2} \atop {j \not= i} }^k
{
  {(1- \overline x_i x_j ) (1- \overline x_i \overline x_j ) }
  \over
   {x_i}
} \right )
 \Omega_{k-2} ( x_2 , \dots , \widehat x_i , \dots , x_k ) \quad ,
\quad (i= 2 , \dots , k )
\eqno(1.5.2a)
$$
$$
\Omega_k ( x_1 , \dots , x_k ) \vert_{x_1 = \overline x_i^{-1} }=
 \left ( \prod_{{j=2} \atop {j \not= i} }^k
{
  {(1- x_i x_j ) (1-  x_i \overline x_j ) }
  \over
   {\overline x_i}
}
\right )
\Omega_{k-2} ( x_2 , \dots , \widehat x_i , \dots , x_k ) \quad ,
\quad (i= 2 , \dots , k )
\eqno(1.5.2b)
$$
$$
\Omega_k ( x_1 , \dots , x_k ) \vert_{x_1 = 0 }=
\Omega_{k-1} ( x_2 , \dots , x_k ) \qquad .
\eqno(1.5.2c)
$$
{\eightrm  [ Type `{\eighttt S152(k):}' in ROBBINS, for specific k.]}
 
{\bf Proof:} By $W(\B_k)$ symmetry, it is enough to prove
(1.5.2a) for $x_1 = x_2^{-1}$, which also implies $(1.5.2b)$. 
To wit, we have to prove that:
$$
x_2^{k-2} \Omega_k ( x_1 , \dots , x_k ) \vert_{x_1 = x_2^{-1} }=
 \left ( \prod_{j=3}^k
  {(1- \overline x_2 x_j ) (1- \overline x_2 \overline x_j ) }
 \right )
 \Omega_{k-2} ( x_3 , \dots , x_k ) \quad .
\eqno(Dominique)
$$
 
We now need:
 
{\bf Subsubsublemma 1.5.2.1:} Let $\Omega_k(x):=\Phi_k(x)/\Delta_k(x)$
be as in the statement of subsublemma $1.5.2$, then
$$
\Omega_k( x_1 , \dots , x_k) \vert_{x_1=x_2^{-1}}
$$
is a Laurent polynomial in $x_2$ of degree $\leq (k-2)$ and
of low-degree $\geq -(k-2)$.
 
[The low-degree, in $x$, of a Laurent polynomial $f(x)$ is 
the degree of $f(x^{-1})$.]
 
{\eightrm  [ Type `{\eighttt S1521(k):}' in ROBBINS, for specific k.]}
 
{\bf Proof:} It is readily checked that
$\Delta_k( x_1 , \dots , x_k) \vert_{x_1=x_2^{-1}}$
is a Laurent polynomial in $x_2$ of degree $(2k-1)$ and
low-degree $-(2k-1)$, the statement of the $(sub)^3$lemma is equivalent
to the statement that $\Phi_k( x_1 , \dots , x_k) \vert_{x_1=x_2^{-1}}$
is a Laurent polynomial in $x_2$ of degree $\leq 3k-3$ and low-degree
$\geq -(3k-3)$. We will prove that this is even true for each and every
single term in the defining sum $(Gog_1)$ of $\Phi_k$.
Let $g=(\pi, \epsilon)$ be a typical signed permutation.
Let $\alpha:=\pi^{-1}(1)$ and $\beta:=\pi^{-1}(2)$, so that
$\pi(\alpha)=1$  and $\pi(\beta)=2$. Without loss of generality
we can assume that $\alpha<\beta$, since in the opposite case
we can exchange $x_1$ and $x_2$.
Let $POL(\hat x_1 , \hat x_2)$ denote a polynomial in $(x_3 , \dots , x_k)$.
Let $z_i:=\epsilon_i(x_i)$ (i.e. $z_i=x_i$ if $\epsilon_i=1$ and
$z_i=\bar x_i (=1- x_i)$ if $\epsilon_i=-1$.)
 
The typical term corresponding to $g=(\pi, \epsilon)$ in the
definition $(Gog_1)$ of $\Phi_k$ can be written as follows:
$$
\displaylines{
g \left [ \prod_{i=1}^{k} \bar x_i^{k-i} x_i^{k} 
\prod_{1 \leq i < j \leq k} (1- x_i \bar x_j )(1- \bar x_i \bar x_j ) 
\right ]
=
 \prod_{i=1}^{k} \bar z_{\pi(i)}^{k-i} z_{\pi(i)}^{k} 
\prod_{1 \leq i < j \leq k} (1- z_{\pi(i)} \bar z_{\pi(j)} )
(1- \bar z_{\pi(i)} \bar z_{\pi(j)} ) 
\cr
=POL( \hat x_1 , \hat x_2 ) \cdot
(\bar z_1^{k- \alpha} z_1^k) \cdot ( \bar z_2^{k- \beta} z_2^k  )
(1- z_1 \bar z_2) (1- \bar z_1 \bar z_2)
\prod_{i=1}^{\alpha-1} (1- z_{\pi(i)} \bar z_1 )(1- \bar z_{\pi(i)} \bar z_1 )
\prod_{{{j=\alpha+1} \atop {j \neq \beta} }}^{k} 
(1- z_1 \bar z_{\pi(j)} )(1- \bar z_1 \bar z_{\pi(j)})
\cdot
\cr
\prod_{{{i=1} \atop {i \neq \alpha} }}
^{\beta-1} (1- z_{\pi(i)} \bar z_2 )(1- \bar z_{\pi(i)} \bar z_2 )
\prod_{j=\beta +1}^{k} 
(1- z_2 \bar z_{\pi(j)} )(1- \bar z_2 \bar z_{\pi(j)}) \quad .
\cr}
$$
 
If $\epsilon_2=-1$, i.e. $z_2=\bar x_2$, the above term must 
vanish when $x_1 x_2=1$, because of 
$(1- z_1 \bar z_2) (1- \bar z_1 \bar z_2)$.
So we can assume that $\epsilon_2=+1$, and so $z_2=x_2$.
We now distinguish two cases:
 
{\bf Case I:} $\epsilon_1=+1$, i.e. $z_1=x_1$, then the term is
$$
\displaylines{
POL( \hat x_1 , \hat x_2 ) \cdot
(\bar x_1^{k- \alpha} x_1^k) \cdot ( \bar x_2^{k- \beta} x_2^k  )
(1- x_1 \bar x_2) (1- \bar x_1 \bar x_2)
\prod_{i=1}^{\alpha-1} (1- z_{\pi(i)} \bar x_1 )(1- \bar z_{\pi(i)} \bar x_1 )
\prod_{{{j=\alpha+1} \atop {j \neq \beta} }}^{k} 
(1- x_1 \bar z_{\pi(j)} )(1- \bar x_1 \bar z_{\pi(j)})
\cdot
\cr
\prod_{{{i=1} \atop {i \neq \alpha} }}
^{\beta-1} (1- z_{\pi(i)} \bar x_2 )(1- \bar z_{\pi(i)} \bar x_2 )
\prod_{j=\beta +1}^{k} 
(1- x_2 \bar z_{\pi(j)} )(1- \bar x_2 \bar z_{\pi(j)}) \quad .
\cr}
$$
 
Now plug in $x_1=x_2^{-1}$. The product of the four products above
times $(1- x_2^{-1} \bar x_2) (1-  \overline{x_2^{-1}} \bar x_2)$
($(1- x_1 \bar x_2) (1- \bar x_1 \bar x_2)$ evaluated at $x_1=x_2^{-1}$)
is a Laurent polynomial
in $x_2$ of degree $\leq 2(k-2)+2=2k-2$ and low-degree
$\geq -(2k-2)$. On the other hand
$(1- x_2^{-1})^{k- \alpha} (1-x_2)^{k- \beta}$ has
degree $\leq k-\beta \leq k-1$ and low-degree 
$\geq -(k-\alpha) \geq -(k-1)$.
So the degree in $x_2$ of the whole term is
$\leq 2k-2+k -1 \leq 3k-3$, and its low-degree is 
$\geq -[(2k-2)+k-1] \geq -(3k-3)$.  This completes
case I of subsubsublemma  $1.5.2.1$.
 
{\bf Case II:} $\epsilon_1=-1$, i.e. $z_1= \bar x_1$, then the term is
$$
\displaylines{
POL( \hat x_1 , \hat x_2 ) \cdot
( x_1^{k- \alpha} \bar x_1^k) \cdot ( \bar x_2^{k- \beta} x_2^k  )
(1- \bar x_1 \bar x_2) (1-  x_1 \bar x_2)
\prod_{i=1}^{\alpha-1} (1- z_{\pi(i)}  x_1 )(1- \bar z_{\pi(i)} x_1 )
\prod_{{{j=\alpha+1} \atop {j \neq \beta} }}^{k} 
(1- \bar x_1 \bar z_{\pi(j)} )(1- x_1 \bar z_{\pi(j)})
\cdot
\cr
\prod_{{{i=1} \atop {i \neq \alpha} }}
^{\beta-1} (1- z_{\pi(i)} \bar x_2 )(1- \bar z_{\pi(i)} \bar x_2 )
\prod_{j=\beta +1}^{k} 
(1-  x_2 \bar z_{\pi(j)} )(1-  \bar x_2 \bar z_{\pi(j)}) \quad .
\cr}
$$
 
Plugging in $x_1=x_2^{-1}$ above yields a Laurent polynomial
in $x_2$ of degree $\leq 2(k-2)+2=2k-2$ and low-degree
$\geq -(2k-2)$, times 
$$
(x_1^{k- \alpha} \bar x_1^k) \cdot ( \bar x_2^{k- \beta} x_2^k  )
=(x_1 x_2)^{k- \alpha} \bar x_1^k \bar x_2^{k- \beta} x_2^{\alpha}=
(1- x_2^{-1})^k (1-x_2)^{k- \beta} x_2^{\alpha}
 \quad .
$$
 
But $(1- x_2^{-1})^k (1-x_2)^{k- \beta} x_2^{\alpha}$ has
degree $\leq k-\beta +\alpha \leq k-1$ and low-degree 
$\geq -(k-\alpha) \geq -(k-1)$.
So the degree in $x_2$ of the whole term is
$\leq 2k-2+k -1 \leq 3k-3$, and its low-degree is 
$\geq -[(2k-2)+k-1] \geq -(3k-3)$.
This completes the proof of subsubsublemma $1.5.2.1$. \halmos
 
Going back to the proof of subsublemma $1.5.2$, 
both sides of $(Dominique)$ are polynomials of degree $\leq 2(k-2)$ 
in $\overline x_2$. 
(By subsubsublemma $1.5.2.1$, The left side is a polynomial in $x_2$ 
of degree $\leq 2(k-2)$, and hence also in $\bar x_2$,
since $x_2=1- \bar x_2 $ .)
The right side vanishes at the $2(k-2)$ values
$\overline x_2 = x_j^{-1} , \overline x_2 = {\overline x_j}^{-1}$,
$j=3, \dots , k$. The next subsubsublemma asserts that
this is true as well for the left side of $(Dominique)$.
 
{\bf Subsubsublemma 1.5.2.2:} The left side of $(Dominique)$,
$$
x_2^{k-2} \Omega_k ( x_1 , \dots , x_k ) \vert_{x_1=x_2^{-1}} \quad ,
$$
when viewed as a polynomial in $\bar x_2$, vanishes at the
$2(k-2)$ values $\bar x_2=x_j^{-1}, \bar x_2=\bar x_j^{-1}$, 
for $j=3, \dots , k$.
 
{\eightrm  [ Type `{\eighttt S1522(k):}' in ROBBINS, for specific k$>$2.]}
 
{\bf Proof:}  By the $W(\B_k)$-symmetry of $\Omega_k(x_1 , \dots , x_k)$,
it suffices to prove that the left side of $(Dominique)$ vanishes
when $\bar x_2=x_3^{-1}$. This is equivalent to proving that
$\Omega_k( x_1 , \dots , x_k  )$ vanishes whenever both
$1-x_1 x_2 =0$ and $1- \bar x_2 x_3=0$. Note that these two last
conditions also imply that
$1- \bar x_1 \bar x_3 =0$. Since $\Delta_k (x)$ does
not vanish there, it suffices to prove that $\Phi_k (x)$ vanishes 
when these relations hold. Putting $z_i:=\epsilon_i(x_i)$,
$\Phi_k( x_1 , \dots , x_k)$ can be written as:
$$
\Phi_k ( x_1 , \dots , x_k ) = 
(-1)^k \sum_{\epsilon \in \{ -1,+1 \}^k } \sum_{ \pi \in \S_k } \sgn(\epsilon)
\sgn(\pi) \,\,
 \prod_{i=1}^{k} \bar z_{\pi(i)}^{k-i} z_{\pi(i)}^{k} 
\prod_{1 \leq i < j \leq k} (1- z_{\pi(i)} \bar z_{\pi(j)} )
(1- \bar z_{\pi(i)} \bar z_{\pi(j)} ) \quad .
\eqno(Bruce)
$$
 
We will prove that each and every term of $(Bruce)$ vanishes when
$1-x_1 x_2=0$ and $1-\bar x_2 x_3=0$. We now split the sum in $(Bruce)$
as follows:
$$
\sum_{\epsilon \in \{ -1,+1\}^k}  \,\,
\sum_{1  \leq \alpha < \beta  < \gamma \leq k}
\sum_{ {{\pi \in \S_k} \atop  {
\{ \pi(\alpha), \pi(\beta), \pi(\gamma) \}=\{1,2,3\}
                             }}
 } \sgn(\epsilon)
\sgn(\pi) \,\,
 \prod_{i=1}^{k} \bar z_{\pi(i)}^{k-i} z_{\pi(i)}^{k} 
\prod_{
       {{1 \leq i < j \leq k}   \atop {i,j \not \in  
                                      \{ \alpha  , \beta, \gamma \}}
        }
       } 
(1- z_{\pi(i)} \bar z_{\pi(j)} )
(1- \bar z_{\pi(i)} \bar z_{\pi(j)} ) \cdot
$$
$$
(1- z_{\pi(\alpha)} \bar z_{\pi(\beta)} )
(1- \bar z_{\pi(\alpha)} \bar z_{\pi(\beta)} )
(1- z_{\pi(\alpha)} \bar z_{\pi(\gamma)} )
(1- \bar z_{\pi(\alpha)} \bar z_{\pi(\gamma)} )
(1- z_{\pi(\beta)} \bar z_{\pi(\gamma)} )
(1- \bar z_{\pi(\beta)} \bar z_{\pi(\gamma)} )
\qquad .
\eqno(Bruce')
$$
 
For each such term in $(Bruce')$, let $\sigma  \in  \S_3$ be defined
by $\sigma(1)=\pi(\alpha)$, $\sigma(2)=\pi(\beta)$, and  
$\sigma(3)=\pi(\gamma)$.
Also let  $\epsilon'$ be the first three components of $\epsilon$:
i.e. $\epsilon'(1)=\epsilon(1)$,
$\epsilon'(2)=\epsilon(2)$ and
$\epsilon'(3)=\epsilon(3)$. The fact that each and every
single term in $(Bruce')$ vanishes when $1-x_1 x_2=0$ and 
$1-\bar x_2 x_3=0$ 
follows from the next $(sub)^4$ lemma:
 
{\bf Subsubsubsublemma 1.5.2.2.1:} For each signed
permutation $(\sigma, \epsilon')$ of $W(\B_3)$,
$$
(\sigma,\epsilon') [ (1- x_1 \bar x_2) (1- \bar x_1 \bar x_2)
(1- x_1 \bar x_3) (1- \bar x_1 \bar x_3)
(1- x_2 \bar x_3) (1- \bar x_2 \bar x_3) ]
$$
vanishes when $x_3=(\bar x_2)^{-1}$ and $x_1=(x_2)^{-1}$.
 
{\eightrm  [ Type `{\eighttt S15221():}' in ROBBINS, to get a rigorous proof.]}
 
{\bf Proof:} \halmos.
 
It follows that indeed each and every single term of $(Bruce')$
vanishes when $x_3=(\bar x_2)^{-1}$ and $x_1=(x_2)^{-1}$, and since
$(2^k k!)\cdot 0=0$, the whole thing vanishes as well. This
completes the proof of subsubsublemma $1.5.2.2$. \halmos
 
Going back to the proof of subsublemma $1.5.2$, 
in order to complete the proof of $(Dominique)$, we need to compare
at one more value of $\bar x_2$. The value to pick is $\bar x_2=0$,
which implies $x_1=x_2=1$. And, indeed
$\Omega_k(1,1,x_3, \dots , x_k)=\Omega_{k-2}(x_3, \dots , x_k)$,
as follows from 
$(\overline{Walt})$  in subsublemma  $1.5.2.3$ below, applied twice.
 
Since both sides of $(Dominique)$ are polynomials of degree $\leq 2(k-2)$
in $\bar x_2$, that agree in $2k-3$ distinct values of $\bar x_2$, they
must be identically equal. This completes the proof of $(Dominique)$
and hence of $(1.5.2a-b)$.
 
By symmetry, the statement of $(1.5.2c)$ is equivalent to $(Walt)$
in the following subsubsublemma.
 
{\bf Subsubsublemma 1.5.2.3:} 
$\Omega_k(x)=\Phi_k(x)/\Delta_k(x)$ ($k \geq 1$) satisfies: 
$$
\Omega_k ( x_1 , \dots , x_{k-1},0 ) =
\Omega_{k-1} ( x_1 , \dots , x_{k-1} ) \quad .
\eqno(Walt)
$$
$$
\Omega_k ( x_1 , \dots , x_{k-1}, 1 ) =
\Omega_{k-1} ( x_1 , \dots , x_{k-1} ) \quad .
\eqno(\overline{Walt})
$$
{\eightrm  [ Type `{\eighttt S1523(k):}' in ROBBINS, for specific k.]}
 
{\bf Proof:} By the   $W(\B_k)$-symmetry  of $\Omega_k$ it suffices
to prove $(Walt)$.
When $x_k =0$ all the summands, in the definition ($Gog_1$) of
$\Phi_k (x)$, vanish except those coming from $(\pi , \epsilon)$
for which $\pi (k) =k$, $\epsilon_k= -1$. Thus the sum reduces
to a sum over $W(\B_{k-1})$ that is easily seen to be equal to
$$
{
 {\Delta_k(x_1, \dots , x_{k-1} ,0) }
 \over
  {\Delta_{k-1}(x_1, \dots , x_{k-1})}
} \Phi_{k-1} (x_1, \dots , x_{k-1} ) \quad .
$$
 
This completes the proof of subsubsublemma 1.5.2.3. \halmos
 
This completes the proof of subsublemma 1.5.2 . \halmos.
 
Combining subsublemmas 1.5.1 and 1.5.2, and using induction on $k$,
we see that indeed $L_k(x) = \Omega_k(x)^2$, for any $k \geq 2$, since 
this is true for $k=0,1$. This completes the proof of
sublemma 1.5. \halmos 
 
Combining Sublemmas 1.1''', 1.2''' , and 1.5 finishes the proof
of Lemma 1. \halmos  \hbox{  } \halmos.
%%end of 1.5 
\vfill
\eject
%begin macros
\baselineskip=14pt
\parskip=10pt
\def\epsilon{\varepsilon}
\def\halmos{\hbox{\vrule height0.15cm width0.01cm\vbox{\hrule height
 0.01cm width0.2cm \vskip0.15cm \hrule height 0.01cm width0.2cm}\vrule
 height0.15cm width 0.01cm}}
\font\eightrm=cmr8  
\font\eighttt=cmtt8
\magnification=\magstephalf

\parindent=0pt
\overfullrule=0in
\headline={\rm  \ifodd\pageno  \RightHead  \else  \LeftHead  \fi}
\def\RightHead{\centerline{Proof  of the ASM Conjecture-References}}
\def\LeftHead{ \centerline{Doron Zeilberger}}
%\pageno=73
%end macros
 
\medskip
 
\centerline{\bf REFERENCES}
 
[A1] G.E. Andrews, {\it Plane Partitions III: the weak
Macdonald conjecture} , Invent. Math. {\bf 53}(1979), 193-225.
 
[A2] G.E. Andrews, {\it Plane Partitions V: the T.S.S.C.P.P
conjecture} , J. Combin. Theory Ser. A, {\bf 66}(1994), 28-39.
 
[D] R.J. Duffin, {\it Basic properties of discrete analytic functions},
Duke Math. J. {\bf 23}(1956), 335-363.
 
[G] I. Gessel, {\it Tournaments and Vandermonde's
 determinant} ,
J. Graph Theory {\bf 3}(1979), 305-307.
 
[Ma] I.G. Macdonald, {\it SYMMETRIC FUNCTIONS AND HALL POLYNOMIALS} ,
Oxford University Press, Oxford, 1979.
 
[MRR1] W.H. Mills, D.P. Robbins, and H.Rumsey,
{\it Proof of the Macdonald conjecture} , Invent. Math. {\bf 66}(1982),
73-87.
 
[MRR2] W.H. Mills, D.P. Robbins, and H.Rumsey, {\it Alternating sign matrices
and descending plane partitions}, J. Combin. Theo. Ser. A, {\bf 34}(1983),
340-359.
 
[MRR3] W.H. Mills, D.P. Robbins, and H.Rumsey,
{\it Self complementary totally symmetric plane partitions} , J.
Combin. Theory Ser. A {\bf 42}(1986), 277-292.
 
[PS] G. P\'olya and G. Sz\"ego, {\it PROBLEMS AND THEOREMS IN ANALYSIS[1]},
Springer, 1976.
 
[R] D.P. Robbins, {\it The story of 1, 2, 7, 42, 429, 7436, ...} ,
Math. Intell {\bf v.13, \#2}(Spring 1991), 12-19.
 
[RR] D.P. Robbins and H. Rumsey, Jr., {\it Determinants and alternating sign
matrices}, Advances in Mathematics {\bf 62}(1986), 169-184.
 
[Stanl] R.P. Stanley {\it A baker's dozen of conjectures concerning
plane partitions} , in: COMBINATOIRE ENUMERATIVE, ed. by G. Labelle
and P. Leroux, Lecture Notes in Mathematics {\bf 1234}, Springer,
Berlin, 1986.
 
[Stant] D. Stanton, {\it Sign variations of the Macdonald identities},
SIAM J. Math. Anal. {\bf 17}(1986), 1454-1460.
 
[Ste] J. Stembridge, {\it A short proof of Macdonald's conjecture
for the root system of type A}, Proc. Amer. Math. Soc. {\bf 102}(1988),
777-786.
 
[WZ] H.S. Wilf and D. Zeilberger,
 {\it An algorithmic proof theory for hypergeometric
(ordinary and ``q") multisum/integral identities}, Invent. Math. 
{\bf 108}(1992), 575-633.
 
[Z1] D. Zeilberger, {\it  The algebra of linear partial difference 
operators and its
applications} , SIAM J. Math. Anal. {\bf 11}(1980) , 919-934.
 
[Z2] D. Zeilberger {\it  Partial difference equations in $ m_1 \geq ... \geq m_n
\geq 0 $ and their applications to combinatorics} , Discrete Math
{\bf 31}(1980) , 65-77.
 
[Z3] D. Zeilberger,
{\it A constant term identity featuring the ubiquitous (and mysterious)
Andrews-Mills-Robbins-Rumsey numbers $\{1,2,7,42,429, ...\}$} , 
J. Combin. Theory Ser. A {\bf 66}(1994), 17-27.
 
[Z4] D. Zeilberger, 
{\it A unified approach to Macdonald's root-system conjectures} ,
SIAM J. Math. Anal. {\bf 19}(1988), 987-1013.
 
[Z5] D. Zeilberger, 
{\it A Stembridge-Stanton style proof of the Habsieger-Kadell
q-Morris identity} , Discrete Math. {\bf 79}(1989/90) , 313-322.
 
[ZB] D. Zeilberger and D.M. Bressoud,
{\it A proof of Andrews' q-Dyson conjecture} ,
Discrete Math. {\bf 54}(1985) , 201-224.
 
Original Version: Kislev 5753; This Version: Nisan  5755.
\vfill
\eject
%%Proof of the Alternating Sign  Matrix Conjecture by D. Zeilberger, Exodion
\baselineskip=14pt
\parskip=10pt
\def\B{{\cal B}}
\def\S{{\cal S}}
\def\halmos{\hbox{\vrule height0.15cm width0.01cm\vbox{\hrule height
 0.01cm width0.2cm \vskip0.15cm \hrule height 0.01cm width0.2cm}\vrule
 height0.15cm width 0.01cm}}
\font\eightrm=cmr8  
\font\et=cmtt8
\font\ebf=cmbx8
\font\eit=cmti8
\magnification=\magstephalf

\parindent=0pt
\overfullrule=0in
\headline={\rm  \ifodd\pageno  \RightHead  \else  \LeftHead  \fi}
\def\RightHead{\centerline{Proof  of the ASM Conjecture-Exodion}}
\def\LeftHead{ \centerline{Doron Zeilberger}}
%\pageno=75
%end macros
 
\centerline{\bf EXODION}
 
\centerline{\bf The Skeleton of the Proof with Checkers' Names}
 
A person's name after the appropriate
$(sub)^i$lemma ($0 \leq i \leq 6$), means that she or he has checked that
its proof is formally correct, provided that all the direct 
descendent-sublemmas are correct. They assume no responsibility whatsoever
to the actual correctness of the offsprings of the $(sub)^i$lemma that they
have been charged with, that lies solely with their subcheckers.
 
1 (G. Andrews, D. Foata, W. Johnson)
 
\quad 1.1 (N.  Alon,  H. Wilf)
 
\qquad \qquad 1.1.1 (A. Bj\"orner)
 
\qquad \qquad \qquad \qquad 1.1.1.1 (C. Orr)
 
\qquad \qquad \qquad \qquad 1.1.1.2 (D. Bar-Natan)
 
\qquad \qquad \qquad \qquad \qquad 1.1.1.2.1 (J. Noonan)
 
\qquad \qquad \qquad \qquad \qquad 1.1.1.2.2 (J. Majewicz)
 
\qquad \qquad \qquad \qquad 1.1.1.3 (S. Cooper)
 
\qquad \qquad 1.1.2 (V. Kann)
 
\qquad \qquad 1.1.3 (R. Stanley)

\quad 1.2 (X.Y. Sun)
 
\qquad \qquad 1.2.1 (B. Salvy)
 
\qquad \qquad \qquad \qquad 1.2.1.1 (D. Stanton)
 
\qquad \qquad \qquad \qquad 1.2.1.2 (Anon., J. Legrange)
 
\qquad \qquad \qquad \qquad \qquad 1.2.1.2.1 (N. Bergeron)
 
\qquad \qquad \qquad \qquad \qquad \qquad 1.2.1.2.1.1 (R. Simion, R.J. Simpson)
 
\qquad \qquad \qquad \qquad \qquad \qquad \qquad 1.2.1.2.1.1.1 (R. Howe,
R. Graham)
 
\qquad \qquad \qquad \qquad \qquad \qquad \qquad 1.2.1.2.1.1.2 (J. Borwein,
D.E. Knuth)
 
\qquad \qquad \qquad \qquad \qquad 1.2.1.2.2 (D. Robbins)
 
\qquad \qquad \qquad \qquad \qquad 1.2.1.2.3 (L. Zhang)
 
\qquad \qquad \qquad \qquad \qquad \qquad 1.2.1.2.3.1 (M. Bousquet-M\'elou)
 
\qquad \qquad \qquad \qquad \qquad \qquad \qquad 1.2.1.2.3.1.1 (L. Habsieger)
 
\qquad \qquad \qquad \qquad \qquad 1.2.1.2.4 (F. Brenti)
 
\qquad \qquad \qquad \qquad \qquad 1.2.1.2.5 (L. Ehrenpreis, S.B. Ekhad,
G. -C. Rota)
 
\qquad \qquad \qquad \qquad 1.2.1.3 (F. Bergeron, M. Wachs)
 
\quad 1.3 (M. Knopp,  S. Parnes)
 
\qquad \qquad 1.3.1 (M. Knopp)
 
\qquad \qquad \qquad \qquad 1.3.1.1 (R. Canfield)
 
\qquad \qquad \qquad \qquad \qquad 1.3.1.1.1 (W. Chen)
 
\qquad \qquad \qquad \qquad 1.3.1.2 (J. Noonan)
 
\qquad \qquad \qquad \qquad 1.3.1.3 (K. Ding)
 
\qquad \qquad \qquad \qquad 1.3.1.4 (P. Zimmermann)
 
\qquad \qquad \qquad \qquad \qquad  Case I (C. Dunkl, R. Ehrenborg)
 
\qquad \qquad \qquad \qquad \qquad \qquad \qquad Subcase Ia (C. Dunkl)
 
\qquad \qquad \qquad \qquad \qquad \qquad \qquad Subcase Ib (C. Dunkl, 
K. Eriksson, R. Sulanke)
 
\qquad \qquad \qquad \qquad \qquad  Case II (M.  Werman)
 
\quad 1.4 (O. Foda)
 
\qquad \qquad 1.4.1 (J. Haglund)
 
\qquad \qquad \qquad \qquad 1.4.1.1 (M. Readdy, N. Takayama)
 
\qquad \qquad \qquad \qquad \qquad  1.4.1.1.1 (A. Fraenkel, J. Labelle)
 
\qquad \qquad \qquad \qquad 1.4.1.2 (J. Friedman)
 
\qquad \qquad \qquad \qquad 1.4.1.3 (S. Okada, Chu W.)
 
\qquad \qquad \qquad \qquad 1.4.1.4 (F. Garvan, G. Gasper)
 
\qquad \qquad \qquad \qquad \qquad  Case I (V. Strehl)
 
\qquad \qquad \qquad \qquad \qquad \qquad \qquad Subcase Ia 
(A.  Granville)
 
\qquad \qquad \qquad \qquad \qquad \qquad \qquad Subcase Ib 
(E. Grinberg)
 
\qquad \qquad \qquad \qquad \qquad  Case II (G. Kalai, I. Sheftel)
 
\qquad \qquad \qquad \qquad \qquad  1.4.1.4.1 (C. Krattenthaler)
 
\qquad \qquad \qquad \qquad \qquad \qquad \qquad  1.4.1.4.1.1 (G. Labelle)
 
\qquad \qquad \qquad \qquad 1.4.1.5 (G. Almkvist, D. Loeb, V. Strehl)
 
\qquad \qquad \qquad \qquad \qquad  Case I  (P. Leroux, V. Strehl)
 
\qquad \qquad \qquad \qquad \qquad \qquad \qquad Subcase Ia 
(E. Lewis, V. Strehl)
 
\qquad \qquad \qquad \qquad \qquad \qquad \qquad Subcase Ib 
(D. Loeb, V. Strehl)
 
\qquad \qquad \qquad \qquad \qquad  Case II (V. Strehl)
 
\qquad \qquad \qquad \qquad \qquad  1.4.1.5.1 
(G. Bhatnagar, S.  Milne, V. Strehl)
 
\qquad \qquad \qquad \qquad 1.4.1.6 (F. Garvan)
 
\quad 1.5 (R. Supper)
 
\qquad \qquad 1.5.1 (Han G.-N., P. Paule, A. Ram)
 
\qquad \qquad \qquad \qquad 1.5.1.1 (S. B. Ekhad, C. Zeilberger)
 
\qquad \qquad \qquad \qquad 1.5.1.2 (S. B. Ekhad, T. Zeilberger)
 
\qquad \qquad 1.5.2 (A. Regev, J. Remmel)
 
\qquad \qquad \qquad \qquad 1.5.2.1 (C. Reutenauer, B. Reznick)
 
\qquad \qquad \qquad \qquad \qquad  Case I (B. Reznick)
 
\qquad \qquad \qquad \qquad \qquad  Case II (B. Reznick, C. Rousseau)
 
\qquad \qquad \qquad \qquad 1.5.2.2 (B. Sagan)
 
\qquad \qquad \qquad \qquad \qquad  1.5.2.2.1 (S.B. Ekhad, H. Zeilberger)
 
\qquad \qquad \qquad \qquad 1.5.2.3 (W. Stromquist, X.G. Viennot)
\eject
\centerline{\bf The Checkers}
\eightrm
 
The sublemmas checked by each of the checkers is given inside the
square bracket at the end of each entry, where the  periods were dropped.
Thus {\ebf [123456]} would stand (if it existed) for 1.2.3.4.5.6.
If the person also contributed general comments about the rest of the paper,
this is indicated by `GC'. Persons who are only credited with `GC' contributed
important general comments, but did not check the formal correctness of any
statement.
 
The following abbreviations were used: 
AcM:=`Acta Math.', AdM:=`Advances in Math',
AMM:=`Amer. Math. Monthly',
BAMS:=`Bulletin of the Amer. Math. Soc', 
CM:=`Contemporary Math.',
CR:=`Comptes Rendus',
DM:=`Discrete Math.', JCT(A):=`J. of Comb. Theory (Ser. A)', 
ICM:=`Inter. Congress of Mathematicians', 
JSC:=`J. for Symbolic Comp.', LNM:=`Lecture Notes in Math'(Springer),
MI:=`Math. Intell.',
PNAS:=`Proc. Nat. Acad. Sci.',
PAMS:=`Proc. Amer. Math. Soc.',
Inv:=`Inventionae Math.', SB:=`Seminaire Bourbaki',.
Sloane's abbreviated citing scheme has  been followed: 
[volume first\_page year].
 
{\ebf Gert Almkvist},  Lund,
{\et gert@maths.lth.se},  is the `guy who generalized a mistake of Bourbaki'
(CM 143 609 93), which lead to important  work in algebraic K-theory. Now
he became a wizard in asymptotic number theory (CM 143 21 93).
{\ebf [1415]}
 
{\ebf Noga Alon}, Tel-Aviv, {\et noga@math.tau.ac.il }, is one
of the most powerful combinatorialists, graph-theorists, and
complexity-ists in the world, who
developed and applied innovative algebraic techniques.
His many honors include ICM 90', and being granted
a permanent-member-size office when he was member of IAS in 93-94.
{\ebf [11]}
 
{\ebf George Andrews}, Penn State,
{\et andrews@math.psu.edu}, is the king of q-theory, who pioneered
computer-assisted research that lead to astounding results in
combinatorics, number theory, and physics (e.g. CBMS \#66,AMS.)
{\ebf [1]}
 
{\ebf Richard Askey}, Wisconsin, {\et askey@math.wisc.edu},
is the king of special functions, of Askey-Gasper inequality 
(AcM 154 137 85) and Askey-Wilson polynomials  fame.
{\ebf [GC]}
 
{\ebf Dror Bar-Natan}, Harvard, {\et dror@math.harvard.edu},
has done ground-breaking work in topological Quantum Field theory,
and has a novel approach to Vassiliev invariants. His 
beautiful papers
can be anonymously ftp-ed from {\et math.harvard.edu}, and once
he moves to Jerusalem (Sept. 1995), from {\et math.huji.ac.il}.
{\ebf [1112]}
 
{\ebf Bernard Beauzamy}, Paris VII  and Inst. Calc. Math.,
{\et beauzamy@mathp7.jussieu.fr }, is a leading Banach spacer
who has also made significant contributions to computer algebra
(JSC 13 463 92.){\ebf [GC]}
 
{\ebf Francois Bergeron}, UQAM, {\et bergeron@catalan.math.uqam.ca},
is a leading Species-ist who has also made beautiful work on card shuffling
and the descent algebra.
{\ebf [1213]}
 
{\ebf Nantel Bergeron}, Harvard, {\et nantel@math.harvard.edu},
contributed significantly to card-shuffling, knot theory (with Bar-Natan),
the descent algebra, and Schubert polynomials.{\ebf [12121]}
 
{\ebf Gaurav Bhatnagar}, Ohio-State, {\et gaurav@math.ohio-state.edu},
is completing a deep and important thesis under Steve Milne. {\ebf [14151]}
 
{\ebf Anders Bj\"orner}, Stockholm, {\et bjorner@math.kth.se},
is a leading algebraic combinatorialist (e.g. ICM 86), who has
applied the arcane M\"obuius function to mundane, but important,
computational complexity problems. {\ebf  [111]}
 
{\ebf Jonathan Borwein}, Simon Fraser,
{\et jborwein@cecm.sfu.ca},
The Borwein brothers 
showed the world (with the Chudnowsky brothers) that computing
$\pi$ to zillion digits involves very beautiful and deep mathematics.
Their paper (with Dilcher) (AMM 100 274 93), 
according to Andrews(MI 16\#4 16 94),
is a prime example of why humans will never be
replaceable by machines. {\ebf [1212112]}
 
{\ebf Mireille Bousquest-M\'elou}, {\et bousquet@labri.u-bordeaux.fr},
has come a long way since she was the `first woman to have ever ranked first
in the concourse to the `great' schools of France.' She has used
Viennot's heaps, as well as the DSV methodology, to find
very deep, and elegant, results in polyomino-enumeration.
(Publ. LACIM (UQAM), \#8). {\ebf [121231]}
 
{\ebf Francesco Brenti}, Perugia and IAS, {\et brenti@math.ias.edu},
applied total positivity to prove deep combinatorial inequalities.
Currently he is studying the Kazhdan-Lusztig polynomials from a combinatorial
point of view.  {\ebf [12124]}
 
{\ebf David M. Bressoud}, Macalister, {\et dvb@math.psu.edu},
is a leading q-ist (JNT 12 87 80) and combinatorialist(DM 38 313 82),
who is also a master expositor(`Second Year Calculus', `Radical Calculus'.)
{\ebf [So far everything except for 13 and 14 and their descendents.]}
 
{\ebf E. Rodney Canfield}, Georgia, {\et erc@csun1.cs.uga.edu },
disproved the Rota conjecture (AdM 29 1 79), and did many
more things since, especially in asymptotic combinatorics. {\ebf [1311]}
 
{\ebf William Chen}, LANL, {\et chen@monsoon.c3.lanl.gov }, did
breathtaking work on bijections(PNAS 87 9635 90), formal grammars, permutation
statistics and many other things. {\ebf [13111]}
 
{\ebf Chu Wenchang}, Rome, {\et wenchang@mat.uniroma1.it}, has
done beautiful work in systematizing and explaining combinatorial
identities. {\ebf [1413]}
 
{\ebf Shaun Cooper}, Wisconsin, {\et scooper@math.wisc.edu}, is
currently completing his thesis under Dick Askey, and has
proved a conjecture of P. Forrester (DM, to appear.)
{\ebf [1113]}
 
{\ebf Kequan Ding}, IAS and Illinois, {\et ding@math.ias.edu}
has done beautiful work on the crossroads of combinatorics,
geometry, and topology. His approach to the Rook polynomials
is wonderful. {\ebf [1313]}
 
{\ebf Charles Dunkl}, Virginia, {\et cfd5z@virginia.edu},
is of `Dunkl operators' fame, without which Cherednick would
not have been able to prove Macdonald's constant term conjecture.
{\ebf [1314, case I, and descendents]}
 
{\ebf Richard Ehrenborg}, {\et ehrenborg@mipsmath.math.uqam.ca},
is both a great juggling mathematician
and a great mathematical juggler.{\ebf [1314, case I, and descendents]}

{\ebf Leon Ehrenpreis}, Temple, {\et leon@math.temple.edu},
did work on partial differential equations that Halmos et. al.
(AMM 83 503 76) consider to be one of the crowning achievements
of American mathematics since 1940. He is listed by
Dieudonn\'e (`Panorama of Pure Math', p. 70,)
along with Laplace, Fourier, and 13 others,
as one of the originators of the general
theory of linear partial differential equations.  {\ebf [12125,GC]}
 
{\ebf Shalosh B. Ekhad}, Temple, {\et ekhad@math.temple.edu},
has written more than fifteen papers, and intends to write many more.
{\ebf [12125,1511,1512,15221]}
 
{\ebf Kimmo Eriksson}, KTH, {\et kimmo@math.kth.se}
has done beautiful work on Coxeter
groups, the Number game, and Fulton's essential set. 
He is also a great speaker. {\ebf [1314, subcase Ib]}
 
{\ebf Dominique Foata}, Strasbourg, {\et foata@math.u-strasbg.fr },
has revolutionized combinatorics at least twice (LNM \#85 and JCT(A)
24 250 78 (see also ICM 83)). 
The first time it also revolutionized theoretical computer
science, and the second time it also revolutionized special functions.
{\ebf [1]}
 
{\ebf Omar Foda}, Melbroune, {\et foda@maths.mu.OZ.AU}, 
is a leading mathematical physicist and quantum group-ist,
who has proved far-reaching Rogers-Ramanujan-type
identities.{\ebf [14]}
 
{\ebf Aviezri Fraenkel},
Weizmann Inst., {\et MAFRAENKEL@WEIZMANN.WEIZMANN.AC.IL},
has both Erd\"os number 1, and Einstein number 2 (via his advisor Ernst
Strauss.) He is a leader in games, combinatorial number theory, complexity,
and many other things. His uncle is Abraham H.  of Zermelo-F. fame.
{\ebf [14111]}
 
{\ebf Jane Friedman}, {\et  janef@teetot.acusd.edu}, has made lasting
contributions to Rogers-Ramanujan theory (Temple thesis, 90),
to the generalized Odlyzko conjecture (PAMS 118 1013 93), and to
lattice path counting. {\ebf [1412]}
 
{\ebf Frank Garvan}, Florida, {\et frank@math.ufl.edu},
cracked a 1944-conjecture of Dyson by cranking out a new
`crank' for partitions (Dyson, in: `Ramanujan Revisited'). He
also did many other things, including the first proof
of the F4-Macdonald root system conjecture (BAMS 24 343 91).
{\ebf [1416]}
 
{\ebf George Gasper}, Northwestern, {\et george@math.nwu.edu},
is of the fame of the Askey-Gasper inequality and of the Gasper-Rahman classic
on `Basic Hypergeometric Series'. He is the most
serious candidate for proving an inequality that 
was shown by Jensen to imply RH.  {\ebf [1414]}
 
{\ebf Bill Gosper}, Independent, {\et rwg@NEWTON.macsyma.com},
is a 20th-century Ramanujan. His Gosper algorithm (PNAS 75 40 78)
changed summation theory for ever. {\ebf [GC]}
 
{\ebf Ron Graham}, Bell, {\et rlg@research.att.com}, ($\in$NAS),
is present and former president of AMS and AJA resp. Both his
technical contributions and charismatic personality (e.g. DLS II,
Univ. Videos) has helped make combinatorics a house-hold name.
{\ebf [1212111]}
 
{\ebf Andrew Granville}, Georgia, {\et andrew@sophie.math.uga.edu },
proved the Carmichael conjecture (with Pomerance)(ICM 94), need I say more? 
{\ebf [1414, subcase Ia]}
 
{\ebf Eric Grinberg}, Temple,  {\et grinberg@math.temple.edu },
has done important work on the Radon transform, and on the
Bussman-Petit problem.  {\ebf [1414, subcase Ib]}
 
{\ebf Laurent Habsieger}, Bordeaux, 
{\et habsieger@geocub.greco-prog.fr}, proved, in his brilliant
thesis,  Askey's conjectured extension of Selberg's integral
(also done by Kadell), and the G2 case of Macdonald's constant
term conjecture (CR 303 211 86),
(also done by Zeilberger.) He has since diversified,
and did beautiful work in number theory. 
Most recently he improved bounds in the notorious
`football(soccer) pool' problem. {\ebf [1212311]}
 
{\ebf Jim Haglund} , Illinois (starting 8/95), 
{until 8/95: \et jhaglund@kscmail.Kennesaw.edu},
has written a brilliant thesis under Rodney Canfield, that
contains lots of beautiful and deep results about rook polynomials,
KOH-type identities, and the Simon Newcomb problem. Lately he is
studying the combinatorics of the Riemann zeta function.
{\ebf [1414, case I]}
 
{\ebf Han Guo-Niu}, Strasbourg, {\et han@math.u-strasbg.fr },
is the world's greatest expert on the Denert statistics,
and is a very powerful bijectionist and word-enumerator.
{\ebf [151]}
 
{\ebf Roger Howe}, Yale, {\et howe@math.yale.edu}, ($\in$NAS),
is not only one of the deepest and most influential Lie theorists
of our day, but also a great lecturer and expositor (AMM  90 353 83
(L. Ford 84' award.)) {\ebf [1212111]}
 
{\ebf Warren Johnsnon}, Penn State, {\et johnson@math.psu.edu},
has done deep and elegant work both on the q-Bell numbers
and q-Abel identities. {\ebf [1,GC]} 
 
{\ebf Gil Kalai}, Jerusalem, {\et kalai@cs.huji.ac.il}, revolutionized
convexity and polytope theory with his `algebraic shifting' (ICM 94).
With Jeff Kahn, he finished off (negatively) the
famous Borsuk problem (that resisted
the attacks of at least one Fields-medalist). {\ebf [1414, case II]}
 
{\ebf Viggo Kann}, KTH, {\et viggo@nada.kth.se}, wrote a brilliant
thesis on approximate algorithms. His computer, Sol Tre, is a co-author
of S.B. Ekhad. {\ebf [112]}
 
{\ebf Marvin Knopp}, Temple, {\et Never had e-mail (and never will)},
wrote the classic `Modular Functions' (Chelsea), started
the notion of `rational period function' and developed a deep theory
of Eichler cohomology. {\ebf [13,131]}
 
{\ebf Don Knuth}, Stanford,
{\et No more e-mail}, ($\in$ NAS), is both the Diderot
(Art of Comp. Prog.), and Gutenberg(TeX) of the computer age. He is the
unique element of TURING $\cap$ HOPPER. {\ebf [1212112,GC]}
 
{\ebf Christian Krattenthaler}, Vienna, {\et kratt@pap.univie.ac.at },
is a Maestro in piano playing, lattice-path combinatorics,
hypergeometric series, and computer algebra. {\ebf [14141]}
 
{\ebf Gilbert Labelle},UQAM, {\et gilbert@mipsmath.math.uqam.ca },
of Species-fame, has the nicest proof ever, of the Largange
inversion formula (AdM 42 217 81.) {\ebf [141411]}
 
{\ebf Jacques Labelle},UQAM, {\et jacques@mipsmath.math.uqam.ca },
is a master chess player, species-ist and combinatorialist.
{\ebf [14111]}
 
{\ebf Jeff Lagarias}, Bell, {\et jcl@research.att.com },
is not only a high-powered researcher in number theory,
cryptology and complexity theory, but is also a master expositor
(AMM  92 3 85, L. Ford award 86'), with a deep sense of history.
With Shor he disproved a long-standing conjecture of Keller
about tiling (BAMS 27 279  92). This leads me to believe that perhaps
the Jacobian conjecture (also posed by Keller) is false as well,
at least in  high dimensions. {\ebf [GC]}
 
{\ebf Jane Legrange}, Bell, {\et jdl@allwise.att.com },
is a Member of Technical Staff in the Display Research Department,
at AT\&T Bell Labs, Murray Hill, NJ.
Thanks to people like her, very soon we won't need paper or
printers. In 1992-1993, she was Visiting Associate Professor at Princeton
University, as a recipient of the coveted VPW NSF grant.
She published extensively on liquid crystals and DNA, in numerous
journals, including {\eit Physical Review Letters} and {\eit Nature}.
{\ebf [1212, GC]}
 
{\ebf Pierre Leroux}, UQAM, {\et leroux@mipsmath.math.uqam.ca},
did beautiful work both on species and (with Viennot) on
the combinatorics of differential equations. With Foata, he
has met Askey's challenge: to find a combinatorial proof to
the generating function for Jacobi polynomials. (PAMS 87 47 83.)
{\ebf [1415,  case I]}
 
{\ebf Ethan Lewis}, Haverford, {\et e1lewis@haverford.edu},
proved a long-standing conjecture on infinite covering systems,
and extended the Cartier-Foata method.
{\ebf [1415,  subcase Ia]}
 
{\ebf Daniel Loeb}, Bordeaux, {\et labri@u-bordeaux.fr},
has (in part with Rota) largely extended the scope of the Umbral Calculus,
to include formal power series of logarithmic type.
He is a master of `negative thinking' and has started a new era
in the umbral calculus {\ebf [1415,subcase Ib]}
 
{\ebf John Majewicz}, Temple, {\et jmaj@math.temple.edu}, is
currently studying for his Ph.D. under D. Zeilberger.
He has extended the WZ method to Abel-type sums. 
His papers can be anon. ftp-ed from
{\et ftp.math.temple.edu}, directory 
{\et pub/jmaj}. {\ebf [11121]}{\ebf [11122]}
 
{\ebf Steve Milne}, Ohio-State, {\et milne@function.mps.ohio-state.edu},
does extremely deep research in multivariate hypergeometric series.
One of Gelfand's dreams is to relate his geometro-combinatorial 
approach to hypergeometric functions with Milne's (and Gustafson's) work.
{\ebf [14151]}
 
{\ebf John Noonan}, Temple, {\et noonan@math.temple.edu}, is
currently studying for his Ph.D. under D. Zeilberger.
He studies permutations that sin a prescribed number of times. 
His great programs and papers can be downloaded from
{\et http://www.math.temple.edu/\~{\quad}noonan}. {\ebf [11121]}
 
{\ebf Andrew Odlyzko}, Bell, {\et amo@research.att.com}.
Contrary to what Herb Wilf (and probably many others) used to
believe, ODLYZKO is {\eit not} a collective name (like Bourbaki)
for a team of brilliant mathematicians, but consists of a singleton.
His many contributions include work on small discriminants,
the disproof of the Mertens conjecture, and vastly extending
the empirical evidence for RH (ICM 86).  {\ebf [GC]}
 
{\ebf Kathy O'Hara}, Independent (formerly Iowa),{\et ohara@math.uiowa.edu},
has astounded the combinatorial world with her injective proof of
the unimodality of the Gaussian polynomials (AMM 96 590 89.)
{\ebf [1212111,GC]}
 
{\ebf Soichi Okada}, Nagoya,{\et okada@nagoya-u.ac.jp},
is one of the world's greatest determinant and pfaffian
evaluators, and an all-around algebraic combinatorialist. 
{\ebf [1413]}
 
{\ebf Craig Orr}, Miami, {\et craig@paris-gw.cs.miami.edu }, has
completed a brilliant thesis under D. Zeilberger. {\ebf [1111]}
 
{\ebf Sheldon Parnes}, Simon Fraser, {\et sparnes@cecm.sfu.ca  },
completed last year a brilliant thesis under D. Zeilberger. {\ebf [13]}
 
{\ebf Peter Paule}, RISC-Linz, {\et ppaule@risc.uni-linz.ac.at},
is at the cutting edge between combinatorics and computer algebra.
He contributed significantly to the Bailey chain, developed a general
theory of summation, and found short and sweet computer proofs
(much shorter then
the original ones by Ekhad and Tre) of the Rogers-Ramanujan identities
(ElJC 1 R10.) {\ebf [151, 152, and their descendents]}
 
{\ebf Bob Proctor}, North Carolina, {\et rap@math.unc.edu},
interfaces the combinatorics of tableaux with
representation theory very elegantly. His beautiful exposition (AMM 89 721 82)
of Stanley's brilliant proof of the Erd\"os-Moser conjecture,
stripped it of its inessential (and intimidating) high-brow
components, and exposed its bare beauty for the rest of us to savor.
{\ebf [13111]}
 
{\ebf Arun Ram}, Wisconsin, {\et ram@math.wisc.edu}, has
a very healthy schizophrenia between combinatorics
and representation theory, that helped both very much.
He is also a brilliant lecturer. {\ebf [151]}
 
{\ebf Marge Readdy}, UQAM {\et readdy@catalan.math.uqam.ca},
is one of the most brilliant academic grandchildren of Richard
Stanley. She is a very versatile posetist. {\ebf [1411, GC]}
 
{\ebf Jeff Remmel}, UCSD, {\et jremmel@math.ucsd.edu },
has done elegant and deep work on tableaux, permutations,
and many other combinatorial objects. The Garsia-Remmel
rook polynomials attracted great attention. {\ebf [152]}
 
{\ebf Amitai Regev}, Penn State and Weizmann Inst.,
{\et regev@math.psu.edu}, proved, in his thesis, that the tensor product
of PI rings is yet another one. He is of `Regev theory' fame
(see Rowan's book) that uses tableaux to study PI algebras.
{\ebf [152]}
 
{\ebf Christoph Reutenauer}, UQAM, {\et christo@catalan.math.uqam.ca },
has a way with words and plays freely with Lie algebras. Read all
about it in the proc. of FPSAC VI.  {\ebf [1521]}
 
{\ebf Bruce Reznick}, Illinois, {\et reznick@math.uiuc.edu},
is one of the greatest 19th-century algebraists alive. By
the periodicity of mathematical fashions, he is also a great
21st-century mathematician. Hilbert, would have
loved his work. {\ebf [1521]}
 
{\ebf Don Richards}, Virginia, {\et dsr2j@virginia.edu},
is a leader on complex functions of matrix arguments. With
Ken Gross, he is developing a deep theory of hypergeometric
functions on complex matrix space (BAMS 24 349 91.)
{\ebf [GC]}
 
{\ebf Dave Robbins}, IDA-CRD, {\et robbins\%idacrd@princeton.edu},
is responsible (with WM and HR), to my spending (wasting(?)) a big chunk of
my life proving (which was fun) and then debugging (which was no fun!)
the ASM conjecture. Before coming up with his notorious conjectures,
he proved his mettle by proving the `most interesting open problem
in enumeration' (Stanley): Macdonald's formula for CSPPs (ref. [MRR1].)
{\ebf [12122]}
 
{\ebf Gian-Carlo Rota}, MIT, {\et rota@math.mit.edu}, ($\in$ NAS), has laid
the foundation to Combinatorial Theory (see `Classic Papers in Comb.',
ed. by Gessel and Rota.)  {\ebf [12125]}
 
{\ebf Cecil Rousseau}, Memphis State, {\et rousseau@hermes.msci.memst.edu}
is a great graph theorist and problem solver. He served as coach of the
US math Olympic team for many years. {\ebf [1521, case II]}
 
{\ebf Bruce Sagan}, Michigan State, {\et sagan@nsf1.mth.msu.edu },
is a great bijectionist and master expositor (`The symmetric group').
{\ebf [1522]}
 
{\ebf Bruno Salvy }, INRIA, {\et Bruno.Salvy@inria.fr },
is pioneering the use of computer algebra for combinatorics
and theoretical computer science. With Paul Zimermann, he
has developed the versatile Maple package {\et gfun}
for handling generating functions.{\ebf [121]}
 
{\ebf Isabel Sheftel}, Jerusalem, {\et sheftel@cs.huji.ac.il}, 
is writing a thesis under the direction of Gil Kalai.
{\ebf [1414, case II]}
 
{\ebf Rodica Simion}, George Washington, {\et simion@gwuvm.gwu.edu},
has done very original and seminal work both on bijections,
posets, enumeration, and (multi-q) theory.
{\ebf [121211]}
 
{\ebf R. Jamie Simpson}, Curtin, {\et TSIMPSONR@cc.curtin.edu.au},
is one to the leaders in covering sequences. He had
a beautiful (and hence elementary)
proof (found independently by Berger, Felzenbaum, 
and Fraenkel) that the top moduli of an exact  covering set must
be equal. {\ebf  [121211]}
 
{\ebf Richard Stanley}, MIT, {\et rstan@math.mit.edu}, ($\in$ NAS), is perhaps
more responsible than anyone else to the creation of a new
AMS classification (05E.) The totality of his influence
is even more than the sum of its parts, which include
having so many brilliant academic offsprings, some of whom
are listed here. {\ebf [113]}
 
{\ebf Dennis Stanton}, Minnesota, {\et stanton@s2.math.umn.edu },
is equally comfortable with block designs (from a fancy, `algebraic'
point of view.), bijections (`Constructive combinatorics', with D. White),
and hypergeometric analysis, to all of which he contributed
deeply. Last but not least he is co-responsible for the `Stanton-Stembridge
trick'. {\ebf [1211]}
 
{\ebf Volker Strehl}, {\et strehl@immd1.informatik.uni-erlangen.de },
has carried the Foata approach to special functions to 
breathtaking  heights. His {\et Habilationschrfft} is a true classic,
that I hope would be published as a book. {\ebf [1414, case I;
1415 and its descendents]}
 
{\ebf Walt Stromquist}, Wagner Assoc., {\et  74730.2425@compuserve.com},
was, for a few months in 1975, the holder of the record for the
minimal number of  countries needed for a planar map not to be
four-colorable. He retains this title, to this day, amongst humans.
Since then he did a lot of things, including instructions
on how to divide a cake fairly. {\ebf [1523]}
 
{\ebf Bob Sulanke}, Boise State, {\et sulanke@math.idbsu.edu},
did beautiful and elegant work on the combinatorics of lattice
paths. He is a master of elegant bijections. {\ebf [1314, subcase Ib]}
 
{\ebf X.Y.Sun}, Temple, {\et xysun@math.temple.edu},
hopes to start his dissertation work soon, under D. Zeilberger.
He is a true Maple whiz. {\ebf [12]}
 
{\ebf Sheila Sundaram}, Miami, {\et sheila@paris-gw.cs.miami.edu}, is a leading
symplectic tableau-ist and plethysm-ist. See FPSAC VI.
{\ebf [GC]}
 
{\ebf Rapha\"ele Supper}, Strasbourg, {\et c/o avanissian@math.u-strasbg.fr},
is a very promising young analyst, who has already contributed
significantly to the theory of entire arithmetic functions (CR 312 781 91,)
and to harmonic functions of exponential type. {\ebf [15]}
 
{\ebf Nobuki Takayama}, Kobe, {\et taka@math.s.kobe-u.ac.jp},
is one of the greatest Gr\"obner bases algorithmicians in
the world, as applied to special functions. He developed
the package KAN. {\ebf [1411]}
 
{\ebf Xavier G. Viennot}, Bordeaux, {\et viennot@geocub.greco-prog.fr }
revolutionized combinatorial statistical physics (SB, Ex. \# 626).
His charismatic research and lecture-persona (`Viennotique'=
transparency wizardry,) helped make combinatorics mainstream, even
in (formerly) snobbish France. He won a CNRS silver medal.
{\ebf [1523,GC]}
 
{\ebf Michelle Wachs}, Miami, {\et wachs@math.miami.edu },
is a leading algebraic combinatorialist. She worked on posets,
permutation statistics, and vastly extended the scope of
{\eit shellabilty}, that, according to Gil Kalai, should from now
on be called {\eit michelle-ability}. {\ebf [1213]}
 
{\ebf Michael Werman}, Jerusalem, {\et werman@cs.huji.ac.il},
is a master (mathematical-) imager (e.g. PRL 5 87 207).
{\ebf [1314, case II]}
 
{\ebf Herb Wilf}, Penn, {\et wilf@central.cis.upenn.edu},
went through the three analyses: numerical, mathematical, and
(happily ever after,) combinatorial. With Nijenhuis and Greene
he designed the best mathematical proof ever. He was a pioneer
in combinatorial algorithms, in particular for random selection and
listings of combinatorial objects. (See  CBMS \#55, SIAM.)
{\ebf [11]}
 
{\ebf Celia Zeilberger}, {\et 71324.3565@compuserve.com},
attends sixth grade at John Witherspoon Middle School,
Princeton, NJ. {\ebf [1511]}
 
{\ebf  Hadas Zeilberger}, {\et 71324.3565@compuserve.com},
attends `A-class' at the U-NOW Day Nursery, Princeton,
that is located on the same street (a couple of blocks away)
from the residence of the (tentative) prover of FLT. {\ebf [15221]}
 
{\ebf Tamar Zeilberger}, {\et 71324.3565@compuserve.com},
attends third grade at Johnson Park Elementary School,
Princeton, NJ. {\ebf [1512]}
 
{\ebf Li Zhang}, Temple, {\et lzhang@math.temple.edu},
is working towards his dissertation, under D. Zeilberger.
He is a great Maple whiz. {\ebf [12123]}
 
{\ebf Paul Zimmermann }, INRIA, {\et Paul.Zimmermann@inria.fr },
is one of the pioneers in the use of computer alegbra for combinatorics
and theoretical computer science. See also under Salvy.{\ebf [1314]}
 
\bye